\documentclass[12pt,a4paper]{article}
\usepackage{amsmath,amsfonts,amssymb,amscd}

\usepackage{amsmath}
\usepackage{tikz-cd}
\usetikzlibrary{matrix, calc, arrows}

\def\Ker{\operatorname{Ker}}

\def\id{\operatorname{id}}

\def\d{\operatorname{d}}
\def\D{\operatorname{D}}

\def\Im{\operatorname{Im}}

\def\mod{\operatorname{mod}}

\newcounter{th}
\def\t{\refstepcounter{th}{\bf \noindent{Theorem} \arabic{th}. }}

\newcounter{prop}
\def\prop{\refstepcounter{prop}{\bf \noindent{Proposition} \arabic{prop}. }}

\newcounter{lem}
\def\lem{\refstepcounter{lem}{\bf \noindent{Lemma} \arabic{lem}. }}

\newcounter{de}
\def\de{\refstepcounter{de}{\bf \noindent{Definition} \arabic{de}. }}

\newcounter{ex}

\begin{document}

\begin{center}
   \Large {\bf\LARGE Graded manifolds of type $\Delta$ and $n$-fold vector bundles}
\end{center}

\begin{center}
	\large {\it Elizaveta Vishnyakova}
\end{center}

\begin{abstract}
	Vector bundles and double vector bundles, or $2$-fold vector bundles, arise naturally for instance as base spaces for  algebraic structures such as  Lie algebroids, Courant algebroids and double Lie algebroids. It is known that all these structures possess a unified description using the language of super\-geometry and $\mathbb{Z}$-graded manifolds of degree $\leq 2$. Indeed, a link has been established between the super and classical pictures
 by the geometrization process, leading to an equivalence of  the category of $\mathbb{Z}$-graded manifolds of degree $\leq 2$ and the category of (double) vector bundles with additional structures. 
 
 In this paper we study  the geometrization process in the case of  $\mathbb Z^r$-graded manifolds of type $\Delta$, where $\Delta$ is a certain weight system and $r$ is the rank of $\Delta$. We establish an equivalence between a subcategory of the category of $n$-fold vector bundles and the category of graded manifolds of type $\Delta$.

\end{abstract}

\section{Introduction}

\textsf{Graded manifolds of type $\Delta$.} A graded manifold of type $\Delta$, a notion that we introduce here, is a natural generalization of the notion of a non-negatively $\mathbb{Z}$-graded manifold of degree $n$. We work in the category of smooth or complex-analytic graded manifolds and we use the language of sheaves and ringed spaces as in the theory of supermanifolds \cite{Leites,Man}.

$\mathbb{Z}$-graded manifolds of degree $n$ were studied by various authors in for instance the context of the theories of Lie algebroids, Courant algebroids, double Lie algebroids and their higher generalizations \cite{Vor,Roytenberg,Bursztyn,LS,Luca,Fernando,JL,BruceGrab SIGMA}.  We can define a non-negatively $\mathbb Z$-graded manifold of degree $n$ as a ringed space which possesses  an atlas with homogeneous coordinates with weights (or degrees) labeled by integers $0,1,\ldots,n$, see \cite{Voronov graded}. In this paper we study non-negatively $\mathbb{Z}^r$-graded manifolds, where $r\geq 1$.  Everywhere {\it graded manifold of type $\Delta$} means $\mathbb{Z}^r$-graded manifold of type $\Delta$, where $r$ is the rank of $\Delta$.

In the case $r>1$ the notion of a degree for graded manifolds is not sufficient to characterize the corresponding category. For example, consider the iterated tangent bundle $T(T(M))$ of a manifold $M$. The structure sheaf of $T(T(M))$ is naturally $\mathbb Z^2$-graded. Indeed, on $T(T(M))$  we can choose  local charts with coordinates in the following form: 
$$
x_i, \,\,\, \d_1(x_j),\,\,\, \d_2(x_s), \,\,\, \d_2(\d_1(x_t)),
$$
where $(x_i)$ are local coordinates on the manifold $M$ which we assume have weight $(0,0)$. Here $\d_1$ and $\d_2$ are the first and second de Rham differentials. Let the local coordinates $\d_1(x_j)$ and $\d_2(x_s)$ and $\d_2(\d_1(x_t))$ have weights  $(1,0)$ and $(0,1)$ and $(1,1)$, respectively. We see that this is a $\mathbb{Z}$-graded manifold of degree $2=1+1$ with respect to the total degree. However in the $\mathbb Z^2$-graded case we can be more precise and consider graded manifolds of multi-degree $(n_1,n_2)$, in the $\mathbb Z^3$-graded case we should consider graded manifolds of multi-degree $(n_1,n_2,n_3)$ and so on. From this point of view $T(T(M))$ is a $\mathbb{Z}^2$-graded manifold of degree $(1,1)$.

 Another observation here is that the numbers $(n_i)$ are also not sufficient to describe the whole picture. For example we can consider a category of $\mathbb{Z}^2$-graded manifolds of degree $(2,2)$ such that any object in this category possesses an atlas with  local coordinates of degrees $(0,0)$, $(2,0)$ and $(0,2)$. We see that in this case we can specify the definition of a $\mathbb{Z}^2$-graded manifold of degree $(2,2)$ and consider the category of graded manifolds of type $\Delta=\{(0,0), (2,0), (0,2)\}$. 
 In addition we can think about $2$-tuples $(0,0)$, $(2,0)$ and $(0,2)$ as vectors in $\mathbb K^2$,  where for convenience we assume that $\mathbb K= \mathbb R$ or $\mathbb C$. We introduce a monoid or a  weight system $\Delta\subset \mathbb K^2$ that parametrizes degrees of local coordinates, see Definition \ref{de weight system} for details.

 Summing up in this paper we study a more precise notion of a non-negatively $\mathbb{Z}^r$-graded manifold, i.e. the notion of a graded manifold of type $\Delta$. In addition we assume as in \cite{Voronov graded} that local coordinates have parities that are related but not determined by the weights. This approach suggests a reduction of some questions about graded manifolds of type $\Delta$ to the study of the combinatorics of the monoid $\Delta$, as it is done for instance in this paper. A further example here is the following.  We consider the root system $A_2$ of the Lie algebra $\mathfrak{sl}_3$. Then certain reflections in this root system correspond to the dualizations of double vector bundles \cite{GraciaMackenzie}, see Section $3.5$ for details. Another benefit of this approach is the possibility to give a  more precise definition of an $n$-fold vector bundle. For instance this  new definition distinguishes between the categories of double vector bundles and of double vector bundles with the trivial core. We introduce the notion of multiplicity free weight system $\Delta$ and study $n$-fold vector bundles of type $\Delta$, where $n$ is the rank of $\Delta$, see Definitions \ref{de n-fold vector bundle}, \ref{de multiplicity free weight system} and \ref{de n-fold vector bundle general}.

\medskip

\noindent\textsf{Geometrization process.} A geometrization process is a functor from the category of graded manifolds to the category of smooth (or holomorphic) manifolds. Such functors are well-known, for example the functor of points for graded or supermanifolds and the linearisation functor \cite{BruceGrGr}. Often it is interesting to ask  which 
  classical manifolds arise from graded manifolds. The goal of this paper is to answer this question for graded manifolds of type $\Delta$.  For motivation, consider the following table of correspondences.

	\begin{center}
	\begin{tabular}{|c|c|}
		\hline
		\quad \quad	{\bf Geometric structures}\quad \quad\quad & {\bf Supergeometric structures}  \\
			\textsf{ base space} / structure& \textsf{ base space} / structure\\
		\hline
			&\\
		\textsf{ vector bundles} & \textsf{ $\mathbb{Z}$-graded manifolds of degree $1$} \\
		Lie algebroids & $[Q,Q]=0$\\
		\hline
		&\\
		\textsf{ metric vector bundles}	 & \textsf{ symplectic $\mathbb{Z}$-graded manifolds} \\
		& \textsf{ of degree $2$}  \\	
		Courant algebroids	& $[Q,Q]=0$\\	
		\hline
		&\\
		\textsf{ metric double vector bundle} & \textsf{ $\mathbb{Z}$-graded manifolds of degree $2$} \\
		VB-Courant algebroids	& $[Q,Q]=0$\\	
		&\\	
		\hline
		& \textsf{ $\mathbb{Z}$-graded manifolds of degree $n>2$ and}\\
	\textsf{ ?} & \textsf{ graded manifolds of type $\Delta$} \\
			& $[Q,Q]=0$\\
		\hline 
	\end{tabular}
\end{center}
Here $Q$ is an odd homological vector field.

The first column of this table represents the world of the classical or commutative geometry, while the second column characterizes the world of super\-geometry. In supergeometry together with usual ``boson-type'' or commuting or even variables one considers also  ``fermion-type'' or anticommuting or odd or ghost variables. 
 Recently it was discovered that the language of supergeometry is very useful for example in the theory of different types of Lie algebroids.  For instance A.~Vaintrob \cite{Va} established an equivalence between categories of Lie algebroids and $\mathbb{Z}$-graded manifolds of degree $1$  with a homological vector field $Q$ of degree $+1$. This equivalence is represented in the first line of our table. 
 Sometimes the super point of view leads to an essential  simplification of the theory of Lie algebroids. For example, the definition of a Lie algebroid
 morphism is quite non-trivial, but this definition has a natural reformulation in terms of homological vector fields.

  The geometrization process is a ``map'' from the second column to the first one. In the case of Lie algebroids however this process is trivial since $\mathbb{Z}$-graded manifolds of degree $1$ are in one-to-one correspondence with vector bundles. In this case the important part is the additional structure, i.e. the homological vector field $Q$.

   The second line of the table represents an equivalence between the categories of Courant algebroids and symplectic $\mathbb{Z}$-graded manifolds of degree $2$ with a certain homological vector field $Q$. The result is due to   P.~\v{S}evera \cite{Severa}  and  D.~Roytenberg \cite{Roytenberg}, independently.  In this case the geometrization process was used not directly.  
   The category of $\mathbb{Z}$-graded manifolds of degree $2$ (not necessary symplectic) was studied by D.~Li-Bland in \cite{Li-Bland}. His result corresponds to the third line of the table. This is an equivalence between the category of $\mathbb{Z}$-graded manifolds of degree $2$ with a certain homological vector field $Q$ and the category of metric double vector bundles with the structure of  VB-Courant algebroid (VB means ``vector bundle''). In other words, to any $\mathbb{Z}$-graded manifold of degree $2$,  Li-Bland assigned a usual manifold, i.e. a metric double vector bundle, and he determined an additional structure that corresponds to the homological vector field $Q$. For  more about applications of Courant algebroids and VB-algebroids, see \cite{BursztynAdv,Courant,GraciaMetha,Gualteri,Kosmann,Mackenzie Crelles}.

    Other results in the direction of the third row of our table were obtained in  \cite{Bursztyn,Fernando,JL}. In these papers the authors assigned to a $\mathbb{Z}$-graded manifold of degree $2$ a usual manifold, or more precisely a double vector bundle with different types of additional structures. They also studied structures determined by a homological vector field $Q$.

   A natural question is to investigate the last line of this table. Thus, in this paper we study more generally     
   non-negatively $\mathbb Z^r$-graded manifolds of type $\Delta$. Due to the complexity we consider graded manifolds without any additional vector fields $Q$. This question is left for the future.

While this paper was in preparation, there appeared another result in this direction \cite{BruceGrab SIGMA}, in which  the authors study  graded bundles of degree $k$ that are  the special case $r=1$ of our graded manifolds of type $\Delta$.  More precisely, in \cite{BruceGrab SIGMA} the authors constructed a functor, which they called the full linearization functor, from the category of graded bundles of degree $k$ ($\mathbb{Z}$-graded manifolds of degree $k$ in our sense, i.e. the parities of coordinates do not necessary coincide with the parities of degrees) to the category of symmetric $k$-fold vector bundles with a family of morphisms that are parametrized by the symmetric group $S_k$. They showed that this functor is an equivalence of categories. Note that in the present paper we consider a different category of  $k$-fold vector bundles, i.e. $k$-fold vector bundles of type $\Delta$ with a family of odd commuting vector fields.

 \medskip
 
\noindent\textsf{Main result.}  Our results can be described as follows. We fix a weight system $\Delta\subset \mathbb K^r$ of rank $r$. Further we choose the parities of the basic weights, see Definition \ref{de weight system}. Then we construct the corresponding multiplicity free weight system $\Delta'= \Delta'(\Delta)$ of rank $r'$, which is in general different from $r$. These two weight systems determine the category $\Delta$\textsf{Man} of graded manifolds of type $\Delta$ and the category $\Delta'$\textsf{VB} of  $r'$-fold vector bundles of type $\Delta'$. 
Further we construct a functor $\mathbb F:$ $\Delta$\textsf{Man}$\to$ $\Delta'$\textsf{VB}, where the main idea is to use the $(r'-r)$-iterated tangent bundle $T\cdots T(\mathcal N)$ of a graded manifold $\mathcal N$. 
(The authors in \cite{BruceGrab SIGMA} introduced  independently a similar construction of the functor $\mathbb F$, the linearisation functor, for the case $r=1$.)
 We finally define the subcategory $\Delta'$\textsf{VBVect} of the category $\Delta$\textsf{Man}, consisting  of $r'$-fold vector bundles of type $\Delta'$ with $(r'-r)$ odd commuting homological vector fields.  These vector fields arise from the iterated de Rham differentials on the structure sheaf of $T\cdots T(\mathcal N)$. We prove that the image of $\mathbb F$ coincides with $\Delta'$\textsf{VBVect}, and moreover that $\mathbb F$ determines an equivalence of the categories $\Delta$\textsf{Man}  and $\Delta'$\textsf{VBVect}.

 \bigskip
 
 \textbf{Acknowledgements:}
E.~V. was partially partially supported by Max Planck Institute for Mathematics, Bonn, by the University of S\~{a}o Paulo, Brazil, FAPESP, grant 2015/15901-9, by SFB TR 191, Germany, and by the Universidade Federal de Minas Gerais.

\section{Graded manifolds of type $\Delta$}

About $\mathbb Z$-graded manifolds of degree $n$ see for instance in \cite{BruceGrab SIGMA,Fernando,GR,JL,LS,Roytenberg,Vor}.

\subsection{A weight system}   Roughly speaking a {\it weight system} $\Delta$  is a monoid that parametrizes weights of local coordinates of a $\mathbb Z^r$-graded manifold. Let us explain this notion in details.

We choose $r$ formal parameters $\alpha_1,\ldots, \alpha_r$, which we will call {\it basic weights}. It is convenient to think about $\alpha_i$ as  vectors in $\mathbb R^r$ or $\mathbb C^r$.  

\medskip
\de\label{de weight system} A {\it weight system} is a subset
\begin{equation}\label{eq Delta as a subset}
\Delta\subset \mathbb Z\alpha_1\oplus \cdots \oplus \mathbb Z\alpha_r
\end{equation}
 satisfying the following properties:
 \begin{enumerate}
 \item $\Delta$ is finite;
 \item $\{0\} \in \Delta$ and $\alpha_i \in \Delta$, where $i=1,\ldots, r$;
 \item if $\delta\in \Delta$ and $\delta = \sum a_i \alpha_i$, where $a_i\in  \mathbb Z$, then $a_i\geqslant 0$.
 \end{enumerate}
The number $r$ is called the {\it rank of $\Delta$}.   

\medskip 

We also will assign a parity $\bar{\alpha}_i\in \{\bar 0, \bar 1\}$  to any basic weight $\alpha_i$. If the parities of $\alpha_i$ are fixed for any $i$, the parities of all other elements from $\Delta$ are determined by the rule $\overline{\delta_1+ \delta_2} = \bar{\delta}_1 + \bar{\delta}_2$. 

\medskip 
 
 In the next section we will introduce {\it graded manifolds of type $\Delta$}. This is graded manifolds with an atlas such that local coordinates are parameterized by elements from $\Delta$.
 The first condition of Definition \ref{de weight system} means that local coordinates of a graded manifold of type $\Delta$ may have only finite number of different weights. It is  a natural agreement for a finite dimensional graded manifold. Further, the first part of the condition $2$ means that, as in the theory of $\mathbb{Z}$-graded manifolds of degree $n$, we have an underlying manifold which structure sheaf is indicated by $0\in \Delta$. The second part of the condition $2$ is technical. The last condition shows that the structure sheaf of our graded manifold is non-negatively $\mathbb Z^r$-graded.

 Examples of weight systems are: 
\begin{equation}\label{eq examples of weight systems}
\Delta_{\mathbf D_2}:= \{0,\,\alpha_1,\,\alpha_2,\,\alpha_1+\alpha_2\},\quad \Delta_{\mathcal M_3}:= \{0,\,\alpha_1, \,2\alpha_1,\, 3\alpha_1\}.
\end{equation}
The weight system $\Delta_{\mathbf D_2}$ corresponds to a double vector bundle $\mathbf D_2$ and  
the weight system $\Delta_{\mathcal M_3}$ corresponds to a $\mathbb Z$-graded manifold $\mathcal M_3$ of degree $3$, see Section $3$. 

Let us give two examples of the parity agreement which we will use in this paper. Consider the weight system $\Delta_{\mathbf D_2}$. In this case we have four possibilities to assign parities for the basic weights $\alpha_1$ and $\alpha_2$. Indeed, we can assume that {\bf (1)} $\bar \alpha_1= \bar \alpha_2=\bar 0$; 
	{\bf (2)} $\bar \alpha_1= \bar \alpha_2=\bar 1$;
	{\bf (3)} $\bar \alpha_1=\bar 0 $ and $\bar \alpha_2=\bar 1 $;
	{\bf (4)} $\bar \alpha_1=\bar 1 $ and $\bar \alpha_2=\bar 0 $. 
All these cases will lead to different categories of graded manifolds of type $\Delta_{\mathbf D_2}$. For example the first case corresponds to the category of pure even double vector bundles,  while in the second case we deal with the category of double vector bundles such that both side bundles of a double vector bundle are odd. Note that the third and the fourth cases lead to isomorphic categories of graded manifolds. More information about all these categories can be found in \cite{Voronov graded}.

In case of the weight system $\Delta_{\mathcal M_3}$ we have two possibilities for the assignment of parity 
$\bar\alpha_1$: $\bar \alpha_1 =\bar 0$ and $\bar \alpha_1 =\bar 1$. Usually in the literature one considers only the second case and the corresponding graded manifolds are called $\mathbb Z$-graded manifold  of degree $3$. The first case corresponds to the category of pure even $\mathbb Z$-graded manifold  of degree $3$, which is less interesting and it is usually omitted.

The construction of the functor, that we study in this paper, works for all possible choices of parities for the basic weights of $\Delta$.
 Note that a different choice of parities for the basic weights of $\Delta$ leads to different categories of graded manifolds of type $\Delta$.

\subsection{Definition of a graded manifold of type $\Delta$}

Let us take a weight system $\Delta$ as in Definition \ref{de weight system} and let us fix parities of the basic weights. In other words, an $r$-tuple $(\bar \alpha_1, \ldots, \bar \alpha_r) \in \mathbb Z_2^r$ is fixed.   Consider a finite dimensional vector space $V$ over $\mathbb K$, where $\mathbb K = \mathbb R$ or $\mathbb C$, with a decomposition into a direct sum of vector subspaces 
 $V_{\delta}$, where ${\delta \in \Delta}$. In other words,
 $$
  V= \bigoplus_{\delta \in \Delta} V_{\delta}.
 $$
  We say that elements from $V_{\delta}\setminus \{0\}$ have weight $\delta$ and have the same parity as $\delta$. In other words, $V_{\delta}$ is a vector subspace in $V$ of parity  $\bar{\delta}$ and of weight $\delta$.  
 Further, we denote by $S^*( V)$ the super-symmetric power of $V$. As usual the weight of a product is the sum of weights of factors. For example the weight of the following product
$$
v_{\delta_1}\cdot v_{\delta_2}\in V_{\delta_1} \cdot V_{\delta_2} \subset S^*( V),
$$
 where $v_{\delta_i}\in  V_{\delta_i}$, is equal to $\delta_1+\delta_2$. 
The same agreement holds for parities.

Consider the $\mathbb Z^r$-graded ringed space $\mathcal U = (\mathcal U_0,\mathcal O_{\mathcal U})$, where $\mathcal U_0 =V^*_0$, and the sheaf $\mathcal O_{\mathcal U}$ is given by the following formula:
\begin{equation}\label{eq structure sheaf of delta domain}
\mathcal O_{\mathcal U}: = \mathcal F_{\mathcal U_0}\otimes_{\mathbb K} S^*\Big(\bigoplus_{\delta \in \Delta\setminus \{0\}} V_{\delta}\Big).
\end{equation}
Here $\mathcal F_{\mathcal U_0}$ is the sheaf of smooth (the case $\mathbb K= \mathbb R$) or holomorphic (the case $\mathbb K= \mathbb C$) functions on $\mathcal U_0 = V^*_0$. 
 The ringed space $\mathcal U$ is a non-negatively  $\mathbb Z^r$-graded ringed space and $\Delta$ is the set of weights of its local coordinates. More precisely, let us choose a basis  $(x_i)$ in $V_0$ and a basis $(\xi_j^{\delta})$ in any $ V_{\delta}$. Then we can consider the set $(x_i, \xi_j^{\delta})_{\delta \in \Delta\setminus \{0\}}$ as the set of local coordinates on $\mathcal U$. We assign the weight $0$ and the parity $\bar 0$ to any $x_i$ and the weight $\delta$ and the parity $\bar \delta$ to any $\xi_j^{\delta}$. We see that the weight system $\Delta$ parametrizes the weights of local coordinates in $\mathcal U$. 
We will call the ringed space $\mathcal U$ a {\it graded domain of type $\Delta$, of parities $(\bar \alpha_1, \ldots, \bar \alpha_r) \in \mathbb Z_2^r$ and of dimension $\{\dim V_{\delta}\}_{\delta \in \Delta}$}. Note that in this case the dimension is a set of numbers parametrized by the elements from $\Delta$.

\medskip

\de\label{de weightd manifold} $\bullet$ A {\it graded manifold of type $\Delta$, of parities $(\bar \alpha_1, \ldots, \bar \alpha_r) \in \mathbb Z_2^r$ and of dimension $\{\dim V_{\delta}\}_{\delta \in \Delta}$} is a $\mathbb Z^r$-graded ringed space $\mathcal N = (\mathcal N_0, \mathcal O_{\mathcal N})$, that is locally isomorphic to a graded domain of type $\Delta$, of parities $(\bar \alpha_1, \ldots, \bar \alpha_r)$ and of dimension 
$
\{\dim V_{\delta}\}_{\delta \in \Delta}.
$

$\bullet$ A {\it morphism of graded manifolds of type $\Delta$ and of parities $(\bar \alpha_1, \ldots, \bar \alpha_r)$} is a morphism of the corresponding $\mathbb Z^r$-graded ringed spaces.     

\medskip

 We will denote the {\it category of graded manifolds of type $\Delta$ and of parities $(\bar \alpha_1, \ldots, \bar \alpha_r)$} by  $\Delta_{(\bar \alpha_1, \ldots, \bar \alpha_r)}$\textsf{Man} or just by $\Delta$\textsf{Man}, when the parity agreement is clear from the context. Note that a graded manifold of type $\Delta$ is defined only if we fix parities $(\bar \alpha_1, \ldots, \bar \alpha_r) \in \mathbb Z_2^r$ of the basic weights in $\Delta$; and to different elements in $(\bar \alpha_1, \ldots, \bar \alpha_r) \in \mathbb Z_2^r$ we assign different categories of graded manifolds of type $\Delta$. If a particular choice of parities $(\bar \alpha_1, \ldots, \bar \alpha_r)$ is not mentioned explicitely, this means that a statement or a construction hold for any choice of the parities.

We can describe a graded manifold of type $\Delta$ in terms of atlases and local coordinates. On a graded manifold $\mathcal N$ of type $\Delta$  there exists an atlas such that in any local chart we can chose local coordinates of weights $\delta\in \Delta$ and we require that transition functions between any two charts preserve all weights.  Note that the structure sheaf of the underlying manifold $\mathcal N_0$ of $\mathcal N$ is equal to  $(\mathcal O_{\mathcal N})_0$ and any homogeneous subsheaf $(\mathcal O_{\mathcal N})_{\delta}$ in $\mathcal O_{\mathcal N}$, where $\delta\in \Delta$, is a $(\mathcal O_{\mathcal N})_0$-locally free sheaf on $\mathcal N_0$.

\section{Examples of graded manifolds of different types}

\subsection{Example 1}
\noindent{\it $\mathbb{Z}$-graded manifolds of degree $n$.} An example of a graded manifold of type $\Delta$ is a $\mathbb Z$-graded manifold of degree $n$.  

\medskip
\de\label{de graded manifold of degree n} 
A {\it $\mathbb{Z}$-graded manifold of degree $n$ and of parity $\bar \alpha_1\in \mathbb Z_2$} is a graded manifold $\mathcal M_n$ of type $\Delta_{\mathcal M_n}$ and of parity $\bar \alpha_1\in \mathbb Z_2$, where
	\begin{equation}\label{eq weight system of M_n}
	\Delta_{\mathcal M_n} = \{0,\alpha_1,\ldots, n \alpha_1\} \subset \mathbb Z \alpha_1.
	\end{equation}
 The number $n$ is called the {\it degree} of $\mathcal M_n$. 
 
\medskip

\subsection{Example 2}
	
\noindent{\it Double and $r$-fold vector bundles.} Another example of graded manifolds of type $\Delta$ is a double and more general an $r$-fold vector bundle. For instance a double vector bundle is a graded manifold of type $\Delta_{\mathbf D_2}$, see (\ref{eq examples of weight systems}).  A  triple vector bundle $\mathbf D_3$ has the weight system 
$$
\Delta_{\mathbf D_3}:= \{0,\alpha_1,\alpha_2,\alpha_3,\alpha_1+\alpha_2, \alpha_1+\alpha_3,\alpha_2+\alpha_3, \alpha_1+\alpha_2+ \alpha_3\}
$$
and so on. We can characterize the weight system $\Delta_{\mathbf D_r}$ in the following way:

\smallskip

{\it  the weight system $\Delta_{\mathbf D_r}$ has rank $r$ and it contains all linear combinations of $\alpha_i$ with coefficients $0$ or $1$, i.e. all linear combinations of $\alpha_i$ without multiplicities. 
}
\smallskip

\medskip

\de\label{de n-fold vector bundle} An {\it $r$-fold vector bundle  of parities $(\bar \alpha_1, \ldots, \bar \alpha_r) \in \mathbb Z_2^r$} is a graded manifold of type $\Delta_{\mathbf D_r}$ and of parities $(\bar \alpha_1, \ldots, \bar \alpha_r)$. 

\medskip

\noindent{\bf Remark.}  This definition of an $r$-fold vector bundle is equivalent to a classical one as was shown in \cite[Theorem 4.1]{GR}, see also \cite{Vor}.

\medskip

In this paper we also will use a more general notion of an $r$-fold vector bundle: an $r$-fold vector bundle of type $\Delta$. For a motivation let us consider an example. Let us take a double vector bundle $\mathbf D_2$ with trivial core (see \cite{Mackenzie} for definitions). In our notations this means that this double vector bundle does not have local coordinates of weight $\alpha_1+\alpha_2$  in a certain atlas. Hence we can assume that the weight system of any double vector bundle with trivial core has the form $\{0,\alpha_1,\alpha_2\}$. For our purpose it is convenient to distinguish these two categories: the category of double vector bundles and the category of  double vector bundles with trivial core. So we will speak about 
the category of double vector bundles of type $\Delta_{\mathbf D_2}$ and the category of double vector bundles of type $\{0,\alpha_1,\alpha_2\}$. Summing up, in this paper we will use the following definitions.

\medskip
\de\label{de multiplicity free weight system} A weight system $\Delta$ is called {\it multiplicity free} if $\Delta$ contains only linear combinations of $\alpha_i$, where $i=1,\ldots,r,$ with coefficients $0$ or $1$. In other words it contains only linear combinations of $\alpha_i$ without multiplicities.

\medskip

Clearly any multiplicity free weight system is contained in some $\Delta_{\mathbf D_r}$ and $\Delta_{\mathbf D_r}$ is the maximal multiplicity free system of rank $r$.

\medskip

\de\label{de n-fold vector bundle general} $\bullet$ An {\it $r$-fold vector bundle of type $\Delta$ and of parities $(\bar \alpha_1, \ldots, \bar \alpha_r)\linebreak \in \mathbb Z_2^r$} is a graded manifold of type $\Delta$ and  of parities $(\bar \alpha_1, \ldots, \bar \alpha_r)$, where $\Delta$ is a multiplicity free weight system of rank $r$.

$\bullet$ A {\it morphism of $r$-fold vector bundles of type $\Delta$ and  of parities $(\bar \alpha_1, \ldots, \bar \alpha_r)\linebreak \in \mathbb Z_2^r$} is a morphism of the corresponding graded manifolds. 

\medskip

Let $\Delta$ be a multiplicity free weight system of rank $r$. 
We will denote the category of  $r$-fold vector bundles of type $\Delta$ and  of parities $(\bar \alpha_1, \ldots, \bar \alpha_r) \in \mathbb Z_2^r$ by $\Delta_{(\bar \alpha_1, \ldots, \bar \alpha_r)}$\textsf{VB} or just by $\Delta$\textsf{VB}.

\medskip

\noindent{\bf Remark.}  In \cite{Vor} the parity reversion functor is defined for $r$-fold vector bundles of type $\Delta_{\mathbf D_r}$. This functor establishes equivalences between all categories of $r$-fold vector bundles of parities $(\bar \alpha_1, \ldots, \bar \alpha_r) \in \mathbb Z_2^r$. In other words all categories of graded manifolds of type $\Delta_{\mathbf D_r}$ and of parities $(\bar \alpha_1, \ldots, \bar \alpha_r) \in \mathbb Z_2^r$ are equivalent. The same holds for 
the category of $r$-fold vector bundles of type $\Delta$, where $\Delta$ is multiplicity free weight system, since the category $\Delta$\textsf{VB} is a subcategory of $\Delta_{\mathbf D_r}$\textsf{VB}. 

\medskip

\noindent{\bf Example.} For more motivation for  Definition \ref{de n-fold vector bundle general}, consider the category of triple vector bundles, that is the category of graded manifolds of type $\Delta_{\mathbf{D}_3}$, see \cite{GraciaMackenzie} and \cite[Example 3]{Vor}. We can visualise a triple vector bundle $\mathbf D$ in the following way
\begin{center}
	\begin{tikzpicture}
	\matrix (m) [matrix of math nodes, row sep=1em,
	column sep=1em]{
		& \mathbf D& &\mathbf D_{23} \\
		\mathbf D_{13} & &\mathbf D_3 & \\
		&\mathbf D_{12} & &\mathbf D_2 \\
	\mathbf	D_1 & & M & \\};
	\path[-stealth]
	(m-1-2) edge (m-1-4) edge (m-2-1)
	edge [densely dotted] (m-3-2)
	(m-1-4) edge (m-3-4) edge (m-2-3)
	(m-2-1) edge [-,line width=6pt,draw=white] (m-2-3)
	edge (m-2-3) edge (m-4-1)
	(m-3-2) edge [densely dotted] (m-3-4)
	edge [densely dotted] (m-4-1)
	(m-4-1) edge (m-4-3)
	(m-3-4) edge (m-4-3)
	(m-2-3) edge [-,line width=6pt,draw=white] (m-4-3)
	edge (m-4-3);
	\end{tikzpicture}
\end{center}
Here all sides of this cube are double vector bundles with certain compatibi\-lity conditions. According to \cite{GR} we can cover $\mathbf D$, $\mathbf D_{ij}$, $\mathbf D_s$ and $M$ with local charts that have the following local coordinates
\begin{align*}
\mathbf D:&\quad (x, \xi^{\alpha_1}, \xi^{\alpha_2}, \xi^{\alpha_3},\xi^{\alpha_1+\alpha_2}, \xi^{\alpha_1+\alpha_3}, \xi^{\alpha_2+\alpha_3}, \xi^{\alpha_1+\alpha_2 + \alpha_2});\\
\mathbf D_{12}: &\quad (x, \xi^{\alpha_1}, \xi^{\alpha_2}, \xi^{\alpha_1+\alpha_2});\\
\mathbf D_{13}: &\quad (x, \xi^{\alpha_1}, \xi^{\alpha_3}, \xi^{\alpha_1+\alpha_3});\\
\mathbf D_{23}: &\quad (x, \xi^{\alpha_2}, \xi^{\alpha_3}, \xi^{\alpha_2+\alpha_3});\\
\mathbf D_{1}: &\quad (x, \xi^{\alpha_1}); \quad \mathbf D_{2}: \quad (x, \xi^{\alpha_2}); \quad \mathbf D_{3}: \quad (x, \xi^{\alpha_3}); \quad M: (x).
\end{align*}
As usual superscript indicates the weight of a coordinate. We have omitted subscripts in order to simplify notation.
 In these charts all projections are given in a  natural way. For example the map $\mathbf D\to \mathbf D_{12}$ in our coordinates is given by 
$$
(x, \xi^{\alpha_1}, \xi^{\alpha_2}, \xi^{\alpha_3},\xi^{\alpha_1+\alpha_2}, \xi^{\alpha_1+\alpha_3}, \xi^{\alpha_2+\alpha_3}, \xi^{\alpha_1+\alpha_2 + \alpha_2})\longmapsto   (x, \xi^{\alpha_1}, \xi^{\alpha_2}, \xi^{\alpha_1+\alpha_2}).
$$

Consider for example the double vector bundle $\mathbf D_{12}$. It may happen that it has trivial core. This means that the intersection of kernels of the projections $\mathbf D_{12}\to \mathbf D_1$ and $\mathbf D_{12}\to \mathbf D_2$ is trivial. In our coordinates this means that we do not have coordinates $\xi^{\alpha_1+\alpha_2}$ of weight $\alpha_1+\alpha_2$. Hence if $\mathbf D_{12}$ has trivial core, then $\mathbf D$ is a graded manifold of type
$$
\Delta_{12}= \{0,\alpha_1,\alpha_2,\alpha_3, \alpha_1+\alpha_3,\alpha_2+\alpha_3, \alpha_1+\alpha_2+ \alpha_3\}.
$$
Here the weight $\alpha_1+\alpha_2$ is omitted. The category of triple vector bundles such that $\mathbf D_{12}$ has trivial core is the category of graded manifold of type $\Delta_{12}$.

Another example is the category of triple vector bundles with trivial ultracore. Recall that the ultracore of $\mathbf D$ is the intersection of kernels of the projections $\mathbf D\to \mathbf D_{ij}$. The category of triple vector bundles with trivial ultracore is the category of graded manifolds of type
$$
\Delta_{123} = \{0,\alpha_1,\alpha_2,\alpha_3, \alpha_1+\alpha_2, \alpha_1+\alpha_3,\alpha_2+\alpha_3\}.
$$

\subsection{Example 3}

\noindent{\it Vector bundles over a graded manifold of type $\Delta$ shifted by a weight.} Let $\mathcal N$ be a graded manifold  of type $\Delta$.  In the category $\Delta$\textsf{Man} we can define a vector bundle in the usual way: {\it a vector bundle over $\mathcal N$  is a graded manifold $\mathbb E$ of type $\Delta$ with a morphism $\mathbb E\to \mathcal N$ that satisfies the usual condition of local triviality and that  has $\mathcal O_{\mathcal N}$-linear transition functions between the trivial pieces.} However sometimes it is more convenient to introduce an additional formal basic weight, say $\beta$, and to think about $\mathbb E$  as about a graded manifold of type $\Delta_{\mathbb E}$, where
\begin{equation}\label{eq weight of the tangent space}
\Delta_{\mathbb E}:= \Delta\cup \{\beta+\delta \,\, |\,\, \delta\in \Delta\}.
\end{equation}

In more details, let us choose local sections $(e_j^{\delta})$ of $\mathbb E$ of weights $\delta \in \Delta$. (The weights of local sections are indicated by the superscript.) Now in any chart we replace the weight $\delta$ of $e_j^{\delta}$ by $\beta+ \delta$.   Clearly, this operation is well-defined and the $\mathcal O_{\mathcal N}$-linearity of transition functions means that the total weight of sections of $\mathbb E$ is preserved. 
In the literature this operation is called the {\bf shift of a vector bundle by a weight $\beta$} and is denoted by $\mathbb E[\beta]$. Below in Section $3.4$ we consider the shift of a tangent bundle in more details.

Summing up, if $\mathbb E$ is a vector bundle over a manifold $\mathcal N$ of type $\Delta$, we will assume that $\mathbb E$ is of type (\ref{eq weight of the tangent space}) for some additional weight $\beta$. We see that the weight $\beta$ has no multiplicity in (\ref{eq weight of the tangent space}), i.e. $\beta$ is contained in weights from (\ref{eq weight of the tangent space}) with coefficients $0$ or $1$. Note that the converse statement is also true. Indeed, let us take a graded manifold $\mathbb E$ of type $\Delta$ such that a certain basic weight $\alpha_i$ has no multiplicity in $\Delta$. Then $\mathbb E$ is a vector bundle over a graded manifold $\mathcal N$ of type 
$$
\Delta':= \Delta \cap \bigoplus_{j\ne i} \mathbb Z \alpha_j. 
$$
The weight system $\Delta'$ satisfies conditions of Lemma \ref{lem Delta' subset Delta denote a smf}, see below. Hence  
$\mathcal N$ is a well-defined graded manifold of type $\Delta'$.

\subsection{Example 4}

\noindent{\it The tangent bundle of a graded manifold of type $\Delta$.} 
Let us remind the construction of the shift of a tangent bundle by a number or by a weight for $\mathbb Z$-graded manifolds of degree $n$ in more details,  see also Section $3.3$. Let $\alpha$ be the basic weight of $\Delta_{\mathcal M_n}$. Consider a $\mathbb Z$-graded manifold $\mathcal M$ of degree $n$ and its tangent bundle $T\mathcal M$. The shift of $T\mathcal M$ by $k\in \mathbb Z$ is denoted usually by $T[k]\mathcal M$. Note that in the present notation it would be more precise to write $T[k\alpha]\mathcal M$ for the shift of $T\mathcal M$ by $k$.

Let $(x_i,\xi_j^{s\alpha})$, where $s=1,\ldots, n$, be local coordinates in $\mathcal M$, where the weights of $x_i$ and $\xi_j^{s\alpha}$ are equal to $0$ and $s\alpha$, respectively. The corresponding local coordinates in $T[k]\mathcal M$ are $(x_i,\xi_j^{\alpha}, \d\! x_i,\d\! \xi_j^{\alpha})$. The notion ''$T\mathcal M$ is shifted by $k\alpha$`` means that the weight of $\d\! x_i $ (and of $\d\! \xi_j^{s\alpha}$) is shifted by $k\alpha$ and it is equal to the weight of $x_i$ plus $k\alpha$ (or it is equal to the weight of $\xi_j^{\alpha}$ plus $k\alpha$, respectively). In other words we assume that the weight of $\d\! x_i $ is equal to $k\alpha$ and the weight of $\d\! \xi_j^{s\alpha}$ is equal to $(s+k)\alpha$. We can easily verify that this definition of the weight does not depend on the choice of local coordinates. Another observation is that instead of the shift by $k\alpha$, we can shift the tangent bundle $T\mathcal M$ by an additional formal weight, say $\beta$. In this case we get  $T[\beta]\mathcal M$. This means that we assume that the weight of $\d\! \xi_j^{s\alpha}$ is equal to $s\alpha + \beta$ and the weight of $\d\! x_i $ is equal to $\beta$. In particular it follows that the weight of the de Rham differential $\d =\d_{dR}$ is equal to $\beta$ in this case since $\d$ sends a coordinate $\xi_j^{s\alpha}$ of weight $s\alpha$ to the coordinate $\d \!\xi_j^{s\alpha}$ of weight $s\alpha + \beta$. Throughout this paper we also assume that the {\bf de Rham differential $\d =\d_{dR}$ is odd}.

In this paper we study iterated tangent bundles shifted by different weights. In fact the order of these shifts is not important and leads to {\it isomorphic} graded manifolds. In more details consider the iterated tangent bundles $\mathcal N^1 :=T[\beta_2] (T[\beta_1] \mathcal M)$ and $\mathcal N^2 := T[\beta_1] (T[\beta_2] \mathcal M)$, where $\beta_i$ are additional weights. These two graded manifolds are naturally isomorphic. Indeed, let $\eta$ be $x_i$ or $\xi_{j}^{s\alpha}$. Then the standard local coordinates in $\mathcal N^1$ and in $\mathcal N^2$ have the following form respectively 
$$
(\eta, \d^1_1 \eta, \d^1_2 \eta, \d^1_2 (\d^1_1 \eta))\quad \text{and} \quad (\eta, \d^2_1 \eta, \d^2_2 \eta, \d^2_2 (\d^2_1 \eta)),
$$
where $\d^i_1$ and $\d^i_2$ are the first and the second  de Rham differentials in $\mathcal N^i$. If $\gamma$ is the weight of $\eta$, then the weights of these coordinates respectively are 
$$
\gamma,\,
\gamma+ \beta_1, \, \gamma+ \beta_2, \,\gamma+ \beta_1+\beta_2 \quad \text{and} \quad \gamma,\,
\gamma+ \beta_2, \, \gamma+ \beta_1, \,\gamma+ \beta_2+\beta_1.
$$
An isomorphism $\mathcal N^1 \simeq \mathcal N^2$ is given by the following formula
$$
 \eta = \eta, \,\, \d^2_1 \eta = \d^1_2 \eta, \,\, \d^2_2 \eta = \d^1_1 \eta, \,\, \d^2_2 (\d^2_1 \eta) =  -\d^1_2 (\d^1_1 \eta).
$$

We will use a similar procedure in the case of graded manifolds of type $\Delta$.
 To simlify notation, we will omit $\beta$'s in itarated tangent bundles, since the idea of the shifts will be clear from the context and all graded manifolds that we will obtain are isomorphic. For the tangent bundle $T\mathcal N$ of a graded manifold $\mathcal N$ of type $\Delta$ we will use a special agreement of shifts. 
 More precisely, let $\mathcal{N}$ be a graded manifold of type $\Delta$ with tangent bundle $T\mathcal{N}$, and let $\alpha_i\in\Delta$ be fixed.  In what follows, we shall consider shifts of the tangent bundle $T\mathcal N$ by $\beta-\alpha_i$, where $\beta$ is an additional weight.
 The shift by $\beta-\alpha_i$ can be understood as the composition of a shift by $-\alpha_i$ and by $\beta$. Explicitly, we have $T[\beta-\alpha_i]\mathcal N = T[\beta][-\alpha_i]\mathcal N$.  Note that in this case the de Rham differential  $\d_{dR}$ must have the weight $\beta-\alpha_i$ as it was pointed above. Since the de Rham differential is odd, the weight $\beta$ has the opposite parity to $\alpha_i$.

To illustrate our construction, let us consider a $\mathbb Z$-graded manifold $\mathcal M$ of degree $n$, see (\ref{eq weight system of M_n}). Here we have only one basic weight $\alpha$ and consider the shift by $\beta - \alpha$, where $\beta$ is an additional weight. The de Rham differential  $\d_{dR}$ has the weight $\beta-\alpha$ and the weight system $\Delta_{T\mathcal M}$ of $T\mathcal M$ is given by the following formula:
$$
\Delta_{T\mathcal M}= \{0,\alpha,\ldots, n \alpha,\, \beta-\alpha,\, \beta,\, \beta+\alpha,\ldots, \beta+ (n-1)\alpha\}.
$$
We see that $T\mathcal M$ is not a non-negatively graded manifold anymore. Indeed, the weight system $\Delta_{T\mathcal M}$ contains weights with negative coefficients, for instance, $\beta-\alpha$. Such shifts we will need further to construct a functor from the category of graded manifolds to the category of $r$-fold vector bundles.

\subsection{Example 5}

 \noindent{\it Root systems of rank $2$ and the corresponding graded manifolds.} There are the following types of 
 root systems of rank $2$: 
 \begin{center}
 	 $A_1\times A_1\simeq D_2$,\,\, $A_2$,\,\, $B_2\simeq C_2$\,\, and \,\, $G_2$.
 \end{center} 
 \begin{center}
	\includegraphics[width=130pt]{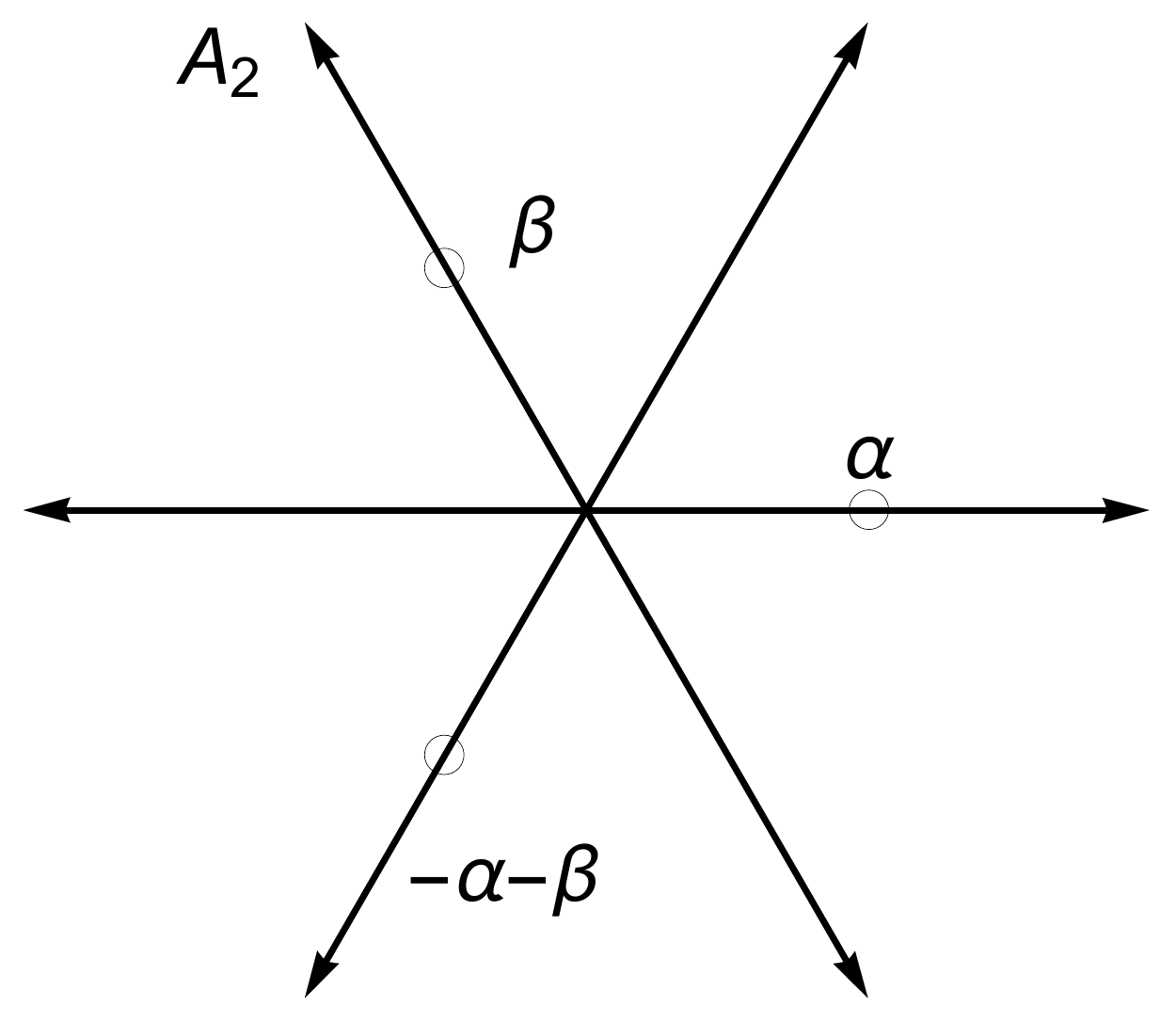} \,\,\,\,\,\,\,\,
	\includegraphics[width=110pt]{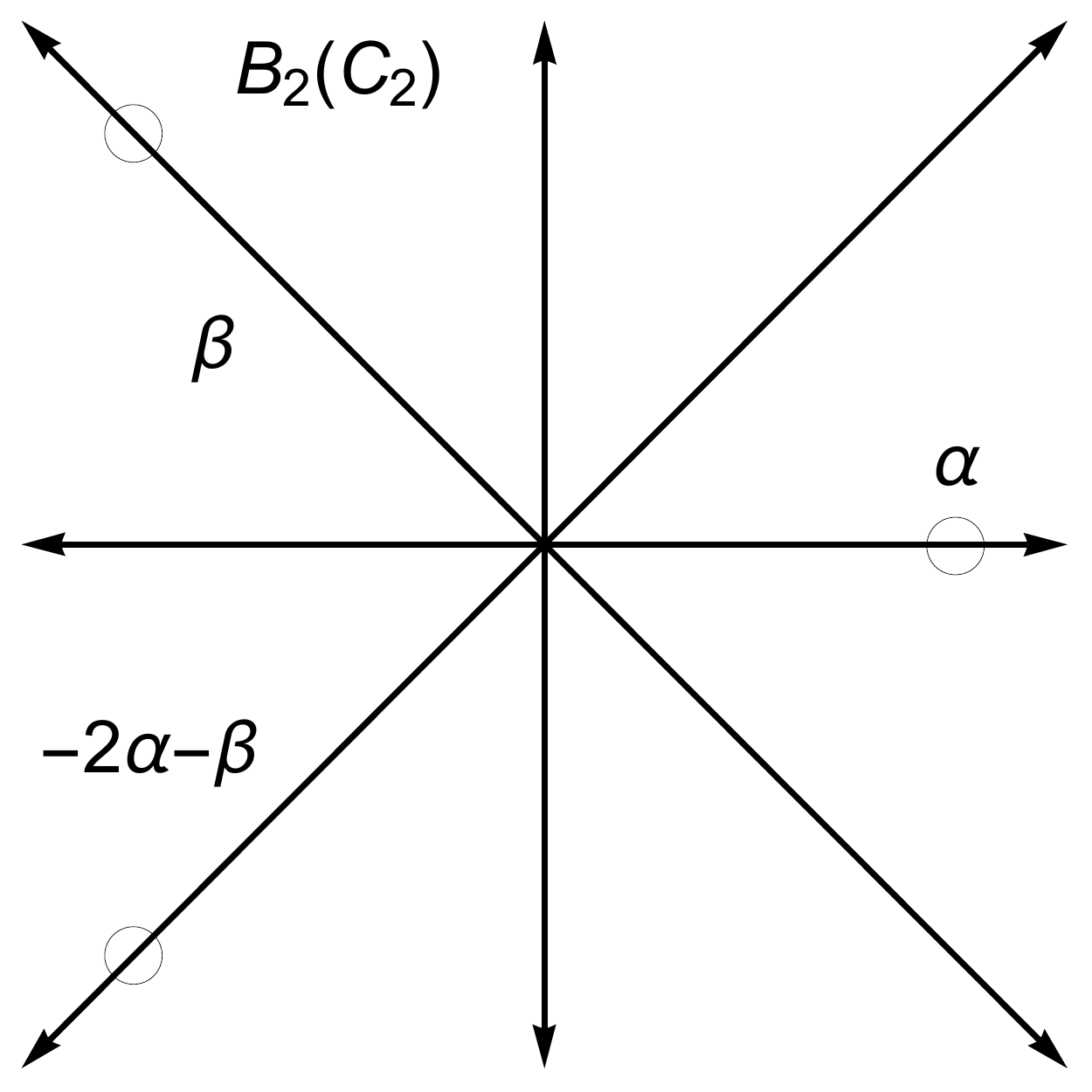} 
\end{center}

\begin{center}
	\includegraphics[width=130pt]{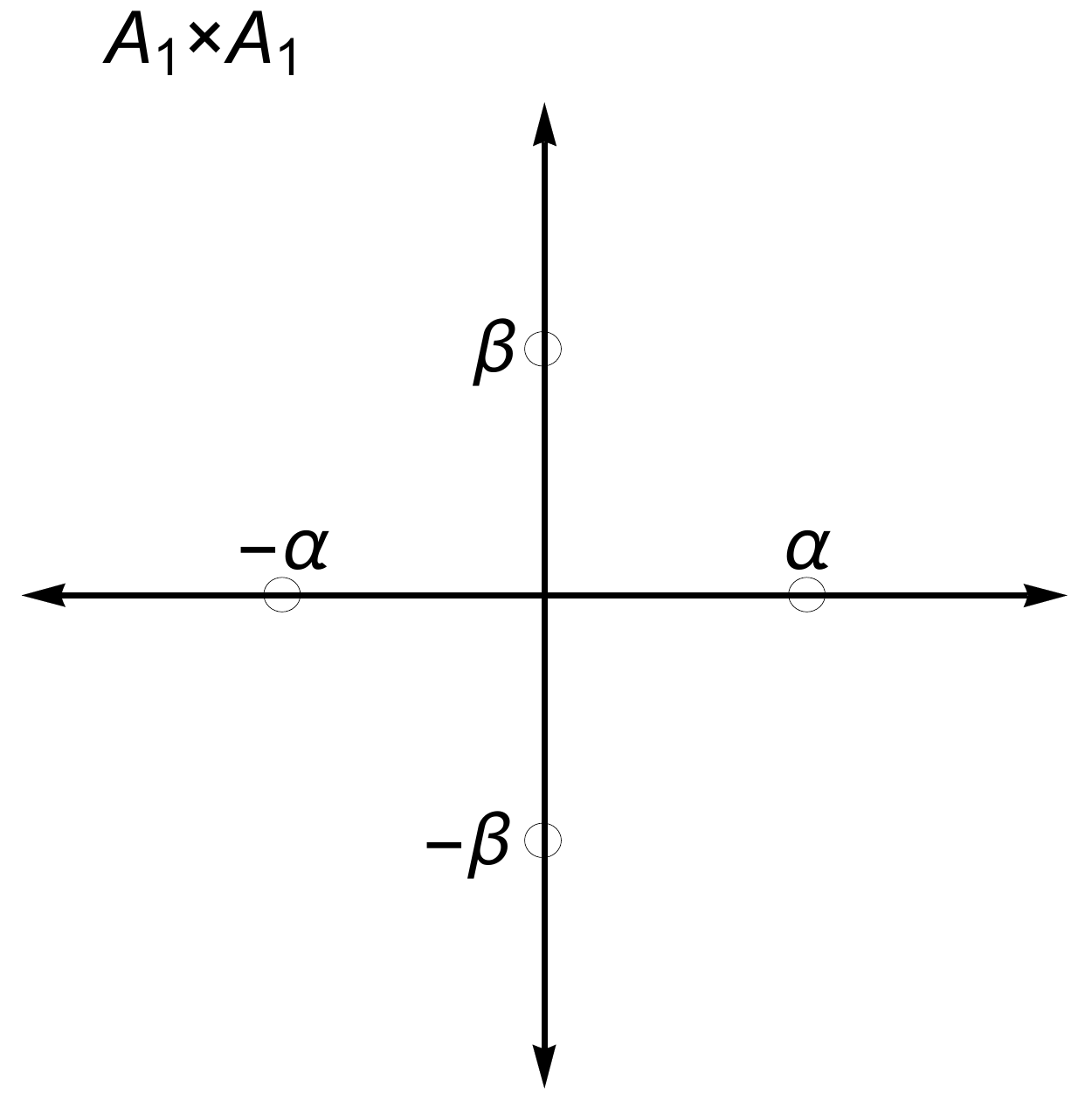} \,\,\,\,\,\,\,\,
	\includegraphics[width=130pt]{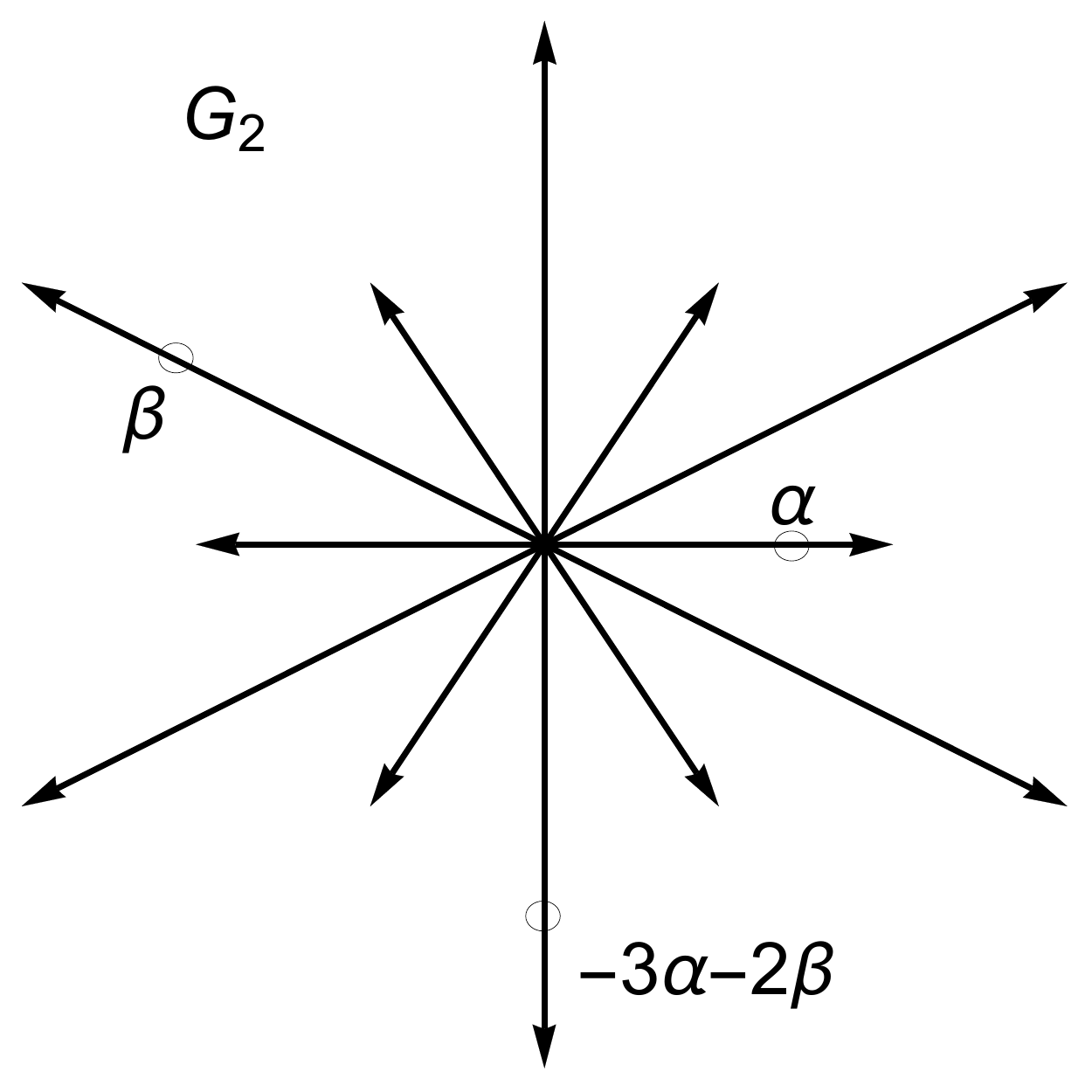} 
\end{center}

 Denote by $\Delta$ a system of positive roots in any of these root systems. Then we can consider the corresponding category of graded manifolds of type $\Delta\cup \{0\}$. Let us characterize these categories.

 \begin{itemize}
 	\item {\bf Case  $A_2$.} Let us choose a system of positive roots $\Delta = \{\alpha,\beta, \alpha+\beta \}$, see picture above. By our definition, see Example $2$, a graded manifold of type $\Delta\cup \{0\}$ is a double vector bundle $\mathbf D_2$:
 	$$
 	\begin{array}{ccc}
 	\mathbf D_2& \longrightarrow& \mathbf A\\
 	\downarrow&&\downarrow\\
 	\mathbf B& \longrightarrow& M\\
 	\end{array}.
 	$$
 	Here $\mathbf D_2\to \mathbf A$, $\mathbf D_2\to \mathbf B$, $\mathbf A\to \mathbf M$ and $\mathbf B\to \mathbf M$ are vector bundles, see \cite{Mackenzie} for precise definition. 
 	
The root system $A_2 = \Delta \cup -\Delta\cup \{0\}$ has also a natural geometric interpretation. Let us choose a chart on $\mathbf D_2$ with the following local coordinates:
 	$$
 	x_i,\,\, \xi_j^{\alpha}, \,\, \xi_s^{\beta}, \,\, \xi_t^{\alpha+\beta}.
 	$$
 	Here $x_i$ are local coordinates of weight $0$ and $\xi_j^{\delta}$ are local coordinates of weight $\delta$, where $\delta\in \Delta$. We use here the standard agreement that $(x_i)$ are local coordinates on $M$, $(x_i,\,\xi_j^{\alpha})$ are local coordinates on $\mathbf A$, $(x_i,\,\xi_s^{\beta})$ are local coordinates on $\mathbf B$ and $(x_i,\,\xi_j^{\alpha}, \,\xi_s^{\beta}, \, \xi_t^{\alpha+\beta})$ are local coordinates on $\mathbf D_2$. Consider the cotangent space $T^*\mathbf D_2$. It has the following local coordinates in the corresponding chart on $T^*\mathbf D_2$:
 	$$
 	\Big\{	x_i,\,\, \xi_j^{\delta}, \,\, \frac{\partial}{\partial 	x_i}, \,\, \frac{\partial}{\partial 	\xi_j^{\delta}}\Big\}_{\delta\in \Delta}.
 	$$
 We assign the weight $-\delta$ to the element $\frac{\partial}{\partial 	\xi_j^{\delta}}$ and the weight $0$ to $\frac{\partial}{\partial 	x_i}$. Hence, $T^*\mathbf D_2$ is a graded manifold of type $A_2$.

Denote by $\hat{T}^*\mathbf D_2$ the graded manifold of type $A_2$ with the structure sheaf $\mathcal O_{\hat{T}^*\mathbf D_2}$ which is locally generated by the following elements: 
 	$$
 	\Big\{	x_i,\,\, \xi_j^{\delta}, \,\, \frac{\partial}{\partial 	\xi_j^{\delta}}\Big\}_{\delta\in \Delta}
 	$$
 and with base $M$. (Clearly, $\mathcal O_{\hat{T}^*\mathbf D_2}$ is a well-defined subsheaf in $\mathcal O_{T^*\mathbf D_2}$.)
 	
 \medskip	
 
 	Further, the double vector bundle $\mathbf D_2$ possesses two dualization operations in the direction $\mathbf A$ and the direction $\mathbf B$. We denote by $\mathbf D_2^{*\mathbf A}$ the dual vector bundle in the direction $\mathbf A$, i.e. $\mathbf D_2^{*\mathbf A}\to \mathbf A$ is the dual of the vector bundle $\mathbf D_2\to \mathbf A$. Similarly we obtain $\mathbf D_2^{*\mathbf B}$. In fact,  $\mathbf D_2^{*\mathbf B}\to \mathbf B$ is the dual of the vector bundle $\mathbf D_2\to \mathbf B$.

It is well-known (see \cite{Mackenzie} and also \cite{GraciaMackenzie}) that  $\mathbf D_2^{*\mathbf A}$ and  $\mathbf D_2^{*\mathbf B}$ are again double vector bundles, hence graded manifolds. We can describe these graded manifolds  using the root system $A_2$. Indeed, the weight system of $\mathbf D_2^{*\mathbf A}$ is 
 	$$
 	\Delta^{*\mathbf A}\cup\{0\},\quad \text{where} \quad \Delta^{*\mathbf A} = \{\alpha, -\alpha-\beta, -\beta \}.
 	 $$
 	Here we can take the weights $\alpha$ and $-\alpha-\beta$ as basic weights, then $-\beta = \alpha +(-\alpha-\beta)$. Similar picture we have for  $\mathbf D_2^{*\mathbf B}$. The weight system of $\mathbf D_2^{*\mathbf B}$ is 
 	$$
 \Delta^{*\mathbf B}\cup\{0\},\quad \text{where} \quad	\Delta^{*\mathbf B} = \{\beta,  -\alpha-\beta, -\alpha\}.
 	$$
 We can take the weights $\beta$ and $-\alpha-\beta$ as basic weights. The structure sheaves of $\mathbf D_2^{*\mathbf A}$ and $\mathbf D_2^{*\mathbf B}$ are subsheaves in the structure sheaf of $\hat{T}^*\mathbf D_2$.

 	Summing up, consider the picture for $A_2$ above, where we  marked the roots $\alpha$, $\beta$ and $-\alpha-\beta$ by a cycle. We see that 
 	
 	\smallskip
 	
 	{\it all double vector bundles that we can obtain from $\mathbf D_2$ using dualizations up to isomorphism correspond to systems of positive roots in $A_2$ such that any of these systems contain exactly two marked roots.}

 	\item {\bf Case  $B_2$.} Consider the following system of positive roots 
 	$$
 	\Delta = \{\alpha,\beta, \alpha+\beta, 2\alpha+\beta \}.
 	$$
 	and the category of graded manifolds of type $\Delta\cup \{0\}$. We see that the weight $\beta$ has no multiplicity in  $\Delta\cup \{0\}$, therefore any graded manifold $\mathcal E$ of this type is a graded vector bundle over a graded manifold $\mathcal M$ of type $\{0,\alpha\}$, see Example $3$.

 	Consider the graded manifold $\hat{T}^*\mathcal E$ of type $B_2 = \Delta\cup -\Delta\cup\{0\}$ that is constructed as in case $A_2$. More precisely the structure sheaf of $\hat{T}^*\mathcal E$ is locally generated by
 	$$
 	\Big\{	x_i,\,\, \xi_j^{\delta}, \,\, \frac{\partial}{\partial 	\xi_j^{\delta}}\Big\}_{\delta\in \Delta}.
 	$$
 	
 Further, we can take the dual vector bundle $\mathcal E^*$ of $\mathcal E$. It is a graded manifold of type $\Delta^*\cup \{0\}$, where 
  	$$
 	\Delta^* = \{\alpha, -2\alpha-\beta,\, -\beta,\, -\alpha-\beta \}.
 	$$
The basic weights here are $\alpha$ and $-2\alpha-\beta$. Again we see that the weight system $\Delta^*$ can be obtained from $\Delta$ by a reflection in the root system $B_2$. Note that the structure sheaves of $\mathcal E$ and $\mathcal E^*$ are subsheaves in the structure sheaf of $\hat{T}^*\mathcal E$.

 	\item {\bf Case  $A_1\times A_1$.} Consider the system of positive roots $\Delta=\{\alpha,\beta\}$. This weight system is multiplicity free, hence it determines the category of certain double vector bundles.
 	 Such double vector bundles are called double vector bundles with trivial core, see  \cite{Mackenzie} for definitions. In other words it is just a sum of two vector bundles. Summing up, the category of graded manifolds of type  $\Delta\cup \{0\}$ is the category of  double vector bundles with trivial core. Again reflections $\alpha \mapsto -\alpha$ and $\beta \mapsto -\beta$ in the weight system $A_1\times A_1$ correspond to dualizations of double vector bundles in different directions.   	 
 	
 	\item {\bf Case  $G_2$.} The category of graded manifolds of type $\Delta\cup \{0\}$, where
 	$$
 	\Delta = \{\alpha,\beta, \alpha+\beta, 2\alpha+\beta, 3\alpha+\beta, 3\alpha+2\beta \}.
 	$$
	We do not know any geometric interpretation in this case.

 \end{itemize}

\section{Constructions of graded manifolds of\\ different types and useful observations}

\subsection{Construction 1}

 Let us take a graded manifold $\mathcal N$ of type $\Delta$. We can associate to $\mathcal N$ a family of graded manifolds of different types. Let $\mathcal N_0$ be the underlying manifold of $\mathcal N$. In our notations, $(\mathcal O_{\mathcal N})_0$ is the structure sheaf of $\mathcal N_0$ and $(\mathcal O_{\mathcal N})_{\delta}$ are $(\mathcal O_{\mathcal N})_0$-locally free sheaves on $\mathcal N_0$, where $\delta\in \Delta$. Let us choose a subset $\Delta'\subset \Delta$ that satisfies the following property: 

\begin{equation}\label{property delta'}
\text{{\it if $\delta \in \Delta'$ and $\delta= \sum\limits_i \delta_i$ for some $\delta_i \in \Delta$, then $\delta_i\in \Delta'$ for any $i$.} }
\end{equation}

\medskip
\lem\label{lem Delta' subset Delta denote a smf}
{\sl Let us take a graded manifold $\mathcal N$ of type $\Delta$.  To any $\Delta'$ satisfying (\ref{property delta'}) we may assign the graded manifold $\mathcal N_{\Delta'}$ of type $\Delta'$.
	
}
\medskip

\noindent{\it Proof.} Consider a local chart on $\mathcal N$ with the structure sheaf in the form  (\ref{eq structure sheaf of delta domain}). Clearly,  we have the following inclusion of the sheaves
$$
 \mathcal F_{\mathcal U_0}\otimes_{\mathbb K} S^*\Big( \bigoplus_{\delta \in \Delta'\setminus \{0\}} V_{\delta}\Big) \hookrightarrow \mathcal F_{\mathcal U_0}\otimes_{\mathbb K} S^*\Big( \bigoplus_{\delta \in \Delta\setminus \{0\}} V_{\delta}\Big).
$$
By (\ref{property delta'}) the transition functions between any such charts preserve this inclusion. Gluing these charts together we get $\mathcal N_{\Delta'}$.$\Box$

\bigskip

\lem\label{lem union and intersection of two Deltas}
{\sl Let as take two weight subsystem $\Delta'$ and $\Delta''$ in $\Delta$  satisfying (\ref{property delta'}). Then 
$$
\Delta'\cap \Delta'' \quad \text{and} \quad \Delta'\cup \Delta'' 
$$ 
also satisfy  (\ref{property delta'}) and determine graded manifolds of type $\Delta'\cap \Delta''$ and $\Delta'\cup \Delta''$. 
}
\medskip

\noindent{\it Proof} we leave for a reader.$\Box$

\subsection{Construction 2}

 Another construction is the following. Let $\Delta$ be  a weight system given by (\ref{eq Delta as a subset}), $\alpha_i$ are basic weights  
 and $\delta=\sum\limits_i a_i\alpha_i\notin \Delta$, where $a_i\geq 0$, be a certain element from lattice (\ref{eq Delta as a subset}).  We set 
 $$
 \Delta':= \Delta\cup \{\delta\}.
 $$
 Assume that a graded manifold $\mathcal N_{\Delta}$ of type $\Delta$ is given.  Our goal now is to add the weight $\delta$ and to  construct another graded manifold $\mathcal N_{\Delta'}$ of type $\Delta'$.

First of all we work with the graded manifold $\mathcal N_{\Delta}$ of type $\Delta$. The structure sheaf of this graded manifold is $\mathbb Z^r$-graded, see Definition \ref{de weightd manifold}. 
This means that we have the decomposition
$$
\mathcal O_{\mathcal N_{\Delta}} = \bigoplus_{\rho\in \mathbb Z \alpha_1\oplus \cdots \oplus \mathbb Z \alpha_r} (\mathcal O_{\mathcal N_{\Delta}})_{\rho}.
$$
Here the sheaf $(\mathcal O_{\mathcal N_{\Delta}})_{\rho}$ is the subsheaf of sections of $\mathcal O_{\mathcal N_{\Delta}}$ of weight $\rho$. Note that some subsheaves $(\mathcal O_{\mathcal N_{\Delta}})_{\rho}$ are equal to $0$. For instance $(\mathcal O_{\mathcal N_{\Delta}})_{\rho} = 0$ if $\rho = \sum\limits_i c_i\alpha_i$, where $c_j<0$ for some $j$. At the same time, it may happen that $\rho \notin \Delta$, however $(\mathcal O_{\mathcal N_{\Delta}})_{\rho}$ is not trivial.

To construct a graded manifold $\mathcal N_{\Delta'}$, we take two $(\mathcal O_{\mathcal N_{\Delta}})_0$-locally free sheaves $\mathcal O_{\delta}$ and $\mathcal E_{\delta}$ on the underlying space $(\mathcal N_{\Delta})_0$ of $\mathcal N_{\Delta}$ such that the sequence 
\begin{equation}\label{eq add a weight exact sequence}
0\to \big(\mathcal O_{\mathcal N_{\Delta}}\big)_{\delta} \to \mathcal O_{\delta} \to \mathcal E_{\delta}\to 0
\end{equation}
is exact.  Note that we always can find such two sheaves. For existence, note that we can take any locally free sheaf $\mathcal E_{\delta}$ on $(\mathcal N_{\Delta})_0$ and define $\mathcal O_{\delta} := \big(\mathcal O_{\mathcal N_{\Delta}}\big)_{\delta} \oplus \mathcal E_{\delta}$. The choice of  $\mathcal O_{\delta}$ and $\mathcal E_{\delta}$  is not unique in general. Therefore the construction of $\mathcal N_{\Delta'}$ depends on the graded manifold $\mathcal N_{\Delta}$ and on the choice of the sequence (\ref{eq add a weight exact sequence}).  Note that any short exact sequence of locally free sheaves is always locally split. Moreover in the category of smooth locally free sheaves a global splitting always exists.

Let us take an atlas $\{U_i\}$ on $(\mathcal N_{\Delta})_0$ such that (\ref{eq add a weight exact sequence}) is split over each $U_i$, i.e. 
$$
\mathcal O_{\delta}|_{U_i} \simeq \big(\mathcal O_{\mathcal N_{\Delta}}\big)_{\delta}|_{U_i} \oplus \mathcal E_{\delta}|_{U_i}.
$$
By definition of locally free sheaves, we may assume that any $U_i$ is small enough such that the sheaf $\mathcal E_{\delta}|_{U_i}$ is free.  Hence, 
$$
\mathcal U_i:= (U_i, \mathcal O_{\mathcal N_{\Delta}}|_{U_i})
$$
is a local chart on $\mathcal N_{\Delta}$. 
Let us choose local coordinates $(x_i)$ in $U_i$, local homogeneous coordinates with non-trivial weights  $(\xi_a^i)$ in $\mathcal U_i$ and a basis of sections $(\eta^i_b)$ in $\mathcal E_{\delta}|_{U_i}\subset \mathcal O_{\delta}|_{U_i}$.  To each $\eta^i_b$ we assign the weight $\delta$. So we are ready to construct a local chart on $\mathcal N_{\Delta'}$. We set
$$
\mathcal V_i:= \big(U_i, \mathcal F_{U_i}\otimes_{\mathbb K} S^*(\xi^i_a,\eta^i_b)\big).
$$
Now our goal is to define transition functions in any intersection $\mathcal V_i\cap \mathcal V_j$.
We take the transition functions $x_j = x_j(x_i)$ and $\xi^j_a = \xi^j_a (x_i, \xi^i_a)$  of $\mathcal N_{\Delta}$ in the intersection of charts $\mathcal U_i \cap \mathcal U_j$  together with the transition functions  $\eta^j_b = \eta^j_b(\xi^i_a, \eta^i_b)$ of $\mathcal O_{\delta}|_{U_i\cap U_j}$. To be more precise, here $\eta^j_b = \eta^j_b(\xi^i_a, \eta^i_b)$ is just the expression of the section $\eta^j_b$ of $\mathcal E_{\delta}|_{U_j}\subset \mathcal O_{\delta}|_{U_j}$ as a linear combination of sections of $\big(\mathcal O_{\mathcal N_{\Delta}}\big)_{\delta}|_{U_i} \oplus \mathcal E_{\delta}|_{U_i} \simeq \mathcal O_{\delta}|_{U_i}$. Summing up, the transition functions in $\mathcal V_i\cap \mathcal V_j$ are $x_j = x_j(x_i)$, $\xi^j_a = \xi^j_a (x_i, \xi^i_a)$ and $\eta^j_b = \eta^j_b(\xi^i_a, \eta^i_b)$. Clearly all total weights are preserved and we get a graded manifold $\mathcal N_{\Delta'}$ of type $\Delta'$. We have completed Construction $2$.

\subsection{A useful observation}

Further we will need the following observation. Let us take two graded manifolds $\mathcal N$ and $\mathcal N'$ of type $\Delta$ with the same underlying space, i.e. $\mathcal N_0 = \mathcal N'_0$. For simplicity of notations we denote the structure sheaves of  $\mathcal N$ and $\mathcal N'$ by $\mathcal O$ and $\mathcal O'$, respectively.  

\medskip

\prop\label{prop isomorphism of graded mnf} {\sl Let us take two graded manifolds $\mathcal N$ and $\mathcal N'$ of type $\Delta$ with the structure sheaves $\mathcal O$ and $\mathcal O'$, respectively, and with  $\mathcal N_0 = \mathcal N'_0$. Assume that the bundle isomorphisms are given $	\varphi_{\delta}: \mathcal O_{\delta} \to \mathcal O'_{\delta}$ for any $\delta \in \Delta.$
If for any 	$\delta \in \Delta$ the following diagram is commutative:
\begin{equation}\label{diagr}
\begin{CD}
\bigoplus\limits_{\delta_1+\delta_2=\delta}
\mathcal O_{\delta_1} \cdot \mathcal O_{\delta_2} @>>> \mathcal O_{\delta}\\
@V{\varphi_{\delta_1} \cdot \varphi_{\delta_2}}VV @VV{\varphi_{\delta}}V\\
\bigoplus\limits_{\delta_1+\delta_2=\delta}
\mathcal O'_{\delta_1}\cdot \mathcal O'_{\delta_2} @>>> \mathcal O'_{\delta}\\
\end{CD}\,\,\,,
\end{equation}
where the horizontal maps are natural inclusions and sums are taken over all $\delta_1,\delta_2\in \Delta\setminus 0$, then $\mathcal N$ and $\mathcal N'$ are isomorphic as graded manifolds of type $\Delta$. 

}

\medskip
\noindent{\it Proof.} Let us define a sheaf isomorphism $\Phi: \mathcal O \to \mathcal O'$ using $(\varphi_{\delta})$. It is sufficient to define $\Phi$ on all sections $f\in \mathcal O|_{U}$, where $U$ is a sufficiently small open set. 
 We put
\begin{enumerate}
	\item $\Phi (f) = \varphi_{\delta} (f)$, if $f\in \mathcal O_{\delta}$ and $\delta\in \Delta$.
	
	\item $\Phi (f) = \varphi_{\delta_1} (f_1)\cdots \varphi_{\delta_k} (f_k)$, if $f= f_1\cdots f_p$ and $f_i\in \mathcal O_{\delta_i}|_{U}$, $\delta_i\in \Delta$.  
\end{enumerate}
Let us prove that $\Phi$ is well-defined. It may happen that $f\in \mathcal O_{\delta}|_{U}$ and $\delta\in \Delta$, and that $f= f_1\cdots f_p$, where $f_i\in \mathcal O_{\delta_i}|_{U}$ and $\delta_1+ \cdots + \delta_k=\delta$, where $\delta_i\in \Delta$. In this case both definitions coincide since Diagrams (\ref{diagr}) are commutative.$\Box$

\section{A functor $\mathbb F$ from $\Delta$\textsf{Man} to  $\Delta'$\textsf{VB} }

In this section we construct a functor $\mathbb F:\Delta$\textsf{Man} $\to \Delta'$\textsf{VB}, where $\Delta'$ is a weight system that will be defined later.
Recall that we denoted by $\Delta$\textsf{Man} the category of graded manifolds of type $\Delta$ and by $\Delta'$\textsf{VB} the category of $r'$-fold vector bundles of type $\Delta'$, where $\Delta'$ is a multiplicity free weight system and $r'$ is the rank of $\Delta'$. If $\mathcal N$ is a graded manifold, we denote by $\mathcal O_{\mathcal N}$ its structure sheaf.

\subsection{Preliminaries}

  Let $\Delta$ be a weight system of rank $r$, $\alpha_i$ be basic weights, see (\ref{eq Delta as a subset}), and $\mathcal N$ be a graded manifold of type $\Delta$. By definition of a weight system, any weight $\delta\in \Delta$ has the form 
\begin{equation}\label{eq delta = sum a_i alpha_i}
\delta= \sum_{i=1}^r a_i(\delta)\alpha_i,
\end{equation}
where $a_i(\delta)$ are non-negative integers. Denote by 
\begin{equation}\label{eq definition of n_i}
n_i:= \max_{\delta\in \Delta}\{a_i(\delta)\} , \quad  i=1,\ldots r.
\end{equation}
In other words $n_i$ is the {\it maximal multiplicity} of the basic weight $\alpha_i$ in the weight system $\Delta$. Note that $n_i>0$. To construct the functor $\mathbb F$, we will need the following set of {\it additional formal weights}:
\begin{equation}\label{eq additional weights beta def}
\{\beta_{ji} \,\, |\,\, j=2,\ldots, n_i,\, \,i=1,\ldots, r \}.
\end{equation}

For our construction we will use sequentially the weights
\begin{equation}\label{eq sequence of beta's}
\beta_{21}, \ldots, \beta_{n_11},\,\, \beta_{22}, \ldots, \beta_{n_22}, \ldots,\,\, \beta_{2r}, \ldots, \beta_{n_rr}. 
\end{equation}
Let us take the first weight from Sequence (\ref{eq sequence of beta's}). For simplicity we assume that the first weight is $\beta_{21}$. As above we denote by $T\mathcal N$ the tangent space of $\mathcal N$ and we denote by $\d_{\beta_{21}}: \mathcal O_{T\mathcal N} \to \mathcal O_{T\mathcal N}$ the corresponding de Rham differential. We assume that the map $\d_{\beta_{21}}$ has the weight $\beta_{21} - \alpha_1$ and we indicate our assumption by the subscript $\beta_{21}$ in $\d_{\beta_{21}}$. In other words, we have 
$$
\d_{\beta_{21}}\big((\mathcal O_{\mathcal N})_{\delta}\big) \subset (\mathcal O_{T\mathcal N})_{\delta+ \beta_{21} - \alpha_1}.
$$ 
(Compare with Example $4$, Section $3$. We have seen there that such weight agreements is well-defined.) Here $(\mathcal O_{\mathcal N})_{\delta}$ is the subsheaf in $\mathcal O_{\mathcal N}$ of weight $\delta$. Note that $\delta$ is not necessary from $\Delta$ in this case. Using our assumption about the weight of $\d_{\beta_{21}}$, we see that the weight system $\Delta_{T\mathcal N}$ of $T\mathcal N$ is given by the following formula:
\begin{equation}\label{eq weight system of TN}
\Delta_{T\mathcal N}= \Delta \cup \{\delta + \beta_{21} - \alpha_1\,\, |\,\, \delta \in \Delta\}.
\end{equation}
Compare this with Formula (\ref{eq weight of the tangent space}). As in Example $4$, we see that $\Delta_{T\mathcal N}$ is not a non-negatively graded manifold anymore. For instance $\Delta_{T\mathcal N}$ contains the weight $\beta_{21}-\alpha_1$
that has a negative coefficient.

Further, let us take the next weight from Sequence (\ref{eq sequence of beta's}). Say the next weight is $\beta_{31}$. Denote by 
$$
\d_{\beta_{31}}: \mathcal O_{TT\mathcal N} \to \mathcal O_{TT\mathcal N}
$$
 the de Rham differential on the tangent space $TT\mathcal N$ of $T\mathcal N$. We assume that $\d_{\beta_{31}}$ has the weight $\beta_{31} - \alpha_1$. We denote the tangent prolongation of $\d_{\beta_{21}}$ on $\mathcal O_{TT\mathcal N}$ also by $\d_{\beta_{21}}$. (The map $\d_{\beta_{21}}$ is a vector field on the graded manifold $T\mathcal N$. Hence the action of $\d_{\beta_{21}}$ is defined on all tensors on $T\mathcal N$ by the Lie derivative.)  By definition of the tangent prolongation the vector fields $\d_{\beta_{21}}$ and $\d_{\beta_{31}}$ on $\mathcal O_{TT\mathcal N}$ commute. Summing up on $TT\mathcal N$ we have two odd homological commuting vector fields $\d_{\beta_{21}}$ and $\d_{\beta_{31}}$. In other words, we have 
 $$
 \d_{\beta_{21}}\circ \d_{\beta_{21}} = 0,\quad \d_{\beta_{31}}\circ \d_{\beta_{31}} = 0, \quad [\d_{\beta_{21}}, \d_{\beta_{31}}] = 0.
  $$  

We continue this process. Altogether we iterate this procedure $n$ times, where
\begin{equation}\label{eq definition of n}
n:=\sum_{i=1}^r n_i-r,
\end{equation}
 using sequentially the de Rham differentials 
\begin{equation}\label{eq order of de Rham diff}
\begin{split}
\d_{\beta_{21}}, \ldots, \d_{\beta_{n_11}},\,\, \d_{\beta_{22}}, \ldots, \d_{\beta_{n_22}}, \ldots,\,\, \d_{\beta_{2r}}, \ldots, \d_{\beta_{n_rr}}. 
\end{split}
\end{equation}
We assume that the de Rham differential $\d_{\beta_{ji}}$ has the  weight $\beta_{ji}-\alpha_i$.
 The result of this procedure is the following  iterated tangent bundle:
$$
\widetilde{\mathcal N}:= \underbrace { T\cdots T}_{n}(\mathcal N)
$$
with $n$ odd operators $\d_{\beta_{ji}}$ such that
\begin{equation}\label{eq commuting de Rham operators on N'}
 [\d_{\beta_{ji}}, \d_{\beta_{j'i'}}] = 0\,\,\, \text{for all}\,\,\, (ji) \,\,\text{and}\,\,\, (j'i').
\end{equation}

Further let $\mathcal R$ be a $\mathbb Z^r$-graded manifold that is not necessary non-negatively graded. In this paper we consider only the case, when $\mathcal R$ is an iterated tangent bundle of a graded manifold of type $\Delta$. Denote by $\mathcal J^{-}_{\mathcal R}$ the ideal in $\mathcal O_{\mathcal R}$ that is generated by all elements with weights that have at least one negative coefficient. To simplify notations usually we will write $\mathcal J^{-}$ instead of $\mathcal J^{-}_{\mathcal R}$.

 \subsection{Construction of $\mathbb F$}

We are ready to define the functor $\mathbb F$. Let $\mathcal N$ be a graded manifold as above. Our goal now is to construct an $r'$-fold vector bundle $\mathbf D_{\mathcal N}$, where $r'$ is defined below. 
Consider the sheaf $ \mathcal O_{\widetilde{\mathcal N}} / \mathcal J^{-}$, where $\widetilde{\mathcal N}$ is defined in the previous section.  Clearly this is the structure sheaf of a certain non-negatively graded manifold. We denote by $\tilde\Delta$ the weight system of this graded manifold and by  $\Delta'\subset \tilde\Delta$
its maximal multiplicity free subsystem. That is $\Delta'$ is the weight subsystem  in $\tilde\Delta$ that contains all weights from $\tilde\Delta$ that have coefficients $0$ or $1$ before the basic weights $\alpha_i$ and $\beta_{ji}$. It is easy to see that $\Delta'$ satisfies conditions of Lemma \ref{lem Delta' subset Delta denote a smf}. We denote by  $\mathbf D_{\mathcal N}$ the corresponding graded manifold of type $\Delta'$, see the construction in  Lemma \ref{lem Delta' subset Delta denote a smf}.

Applying the rule (\ref{eq weight system of TN}), we see that $\alpha_i, \beta_{ji}\in \tilde\Delta$ and we note  that these weights are multiplicity free. Hence, $\alpha_i, \beta_{ji}\in \Delta'$ and the rank of  $\Delta'$ is equal to $r':=n+r$. Further, by definition the weight system $\Delta'$ is multiplicity free, hence, the graded manifold {\it $\mathbf D_{\mathcal N}$ is an  $r'$-fold vector bundle}, see Definition \ref{de n-fold vector bundle general}. We put 
\begin{equation}\label{eq def of functor F}
\mathbb F(\mathcal N) := \mathbf D_{\mathcal N}.
\end{equation}
Note that $\Delta'$ depends only on $\Delta$, but not on a particular choice of $\mathcal N$.

Further, let us take a morphism $\Phi:\mathcal N\to \mathcal N_1$ of two graded  manifolds $\mathcal N$ and $\mathcal N_1$ of type $\Delta$. By definition, $\Phi$ preserves all weights. We have the corresponding map in the iterated tangent bundles
$$
\big(\underbrace { T\cdots T}_{n}\Phi\big):\underbrace { T\cdots T}_{n}(\mathcal N) \to \underbrace { T\cdots T}_{n}(\mathcal N_1),
$$
that  preserves all weights. Therefore, the map 
$$
(\Phi')^*: \mathcal O_{\widetilde{\mathcal N}_1}/ \mathcal J^{-}_{\widetilde{\mathcal N}_1}\to \mathcal O_{\widetilde{\mathcal N}} /\mathcal J^{-}_{\widetilde{\mathcal N}}
$$
is well-defined. Here we use the following notations:
$$
\widetilde{\mathcal N}:= \underbrace { T\cdots T}_{n}(\mathcal N),\quad \quad\widetilde{\mathcal N}_1:= \underbrace { T\cdots T}_{n}(\mathcal N_1),
$$
and we denote by $\mathcal J^{-}_{\widetilde{\mathcal N}}$ and $\mathcal J^{-}_{\widetilde{\mathcal N}_1}$ the ideals in $\mathcal O_{\widetilde{\mathcal N}}$ and $\mathcal O_{\widetilde{\mathcal N}_1}$, respectively, that are defined in the previous section, i.e. these ideals are generated by all elements with weights that have at least one negative coefficient. Since $(\Phi')^*$ preserves all weights, we get
$$
(\Phi')^* \big( \mathcal O_{\mathbf D_{\mathcal N_1}} \big) \subset \mathcal O_{\mathbf D_{\mathcal N}}.
$$
Hence the map 
$$
\mathbb F (\Phi): \mathbf D_{\mathcal N}\to \mathbf D_{\mathcal N_1} 
$$
is defined. 
 Clearly, the correspondence $\Phi \mapsto \mathbb F (\Phi)$ sends a composition of morphisms to a composition of morphisms. Hence we obtain the following theorem. 

\medskip

\t\label{teor F is a functor} {\sl 
The correspondence $\mathbb F$ is a functor from the category of graded manifolds of type $\Delta$ to the category of $r'$-fold vector bundles of type $\Delta'$. 
	
}

\medskip

\subsection{Explicit description of $\Delta' = \Delta'(\Delta)$ }

Let us describe $\Delta' = \Delta'(\Delta)$ explicitly. We take $\delta \in \Delta$ with coefficients $a_i(\delta)\in \mathbb Z$ as in (\ref{eq delta = sum a_i alpha_i}), and we put
$$
\beta_{I_ii} := \sum_{s\in I_i} \beta_{si},\quad \text{where  $I_i\subset\{2,\ldots, n_i\}$.}
$$

\medskip

\prop\label{prop obtain Delta_D from Delta_N}
{\sl Let $\Delta$ be a weight system and $\Delta' = \Delta'(\Delta)$ be the weight system contracted in Section $5.2$. 
	We have
	\begin{equation}\label{eq Delta_D_N expression new}
	\Delta' = \bigcup_{\delta\in \Delta} \Delta'_{\delta},
	\end{equation}
	where $\Delta'_{\delta}$ is given by the following formula:
	\begin{equation}\label{eq Delta_delta new}
	\begin{split}
	\Delta'_{\delta}:= \Big\{\delta+ \sum_{i=1}^r  (\beta_{I_i}  - |I_i| \alpha_i)\,\,| \,\, I_i\subset\{2,\ldots, n_i\} \,\, \text{such that}\,\,\\
	|I_i|= a_i(\delta) \,\, \text{or} \,\, |I_i|= a_i(\delta)-1	\Big\}.
	\end{split}
	\end{equation}

}

\medskip

\noindent{\it Proof.} Consider the weight system $\Delta_{T\mathcal N}$ that is given by Formula (\ref{eq weight system of TN}). If we iterate this process we see that 
the weight system of $\widetilde{\mathcal N}$ is given by
$$
	\Delta_{\widetilde{\mathcal N}} = \bigcup_{\delta\in \Delta} \widetilde{\Delta}_{\delta},\quad\text{where}\quad \widetilde{\Delta}_{\delta}:= \Big\{\delta+ \sum_{i=1}^r  (\beta_{I_i}  - |I_i| \alpha_i)\,\,| \,\,
	|I_i|= 0,1,\cdots, n_i	\Big\}.
$$
 If we remove from $\widetilde\Delta_{\delta}$ all weights with at least one negative coefficient and all weights with non-trivial multiplicities, we get (\ref{eq Delta_delta new}). This finishes the proof.$\Box$

\medskip

\noindent{\bf Example.} Consider the weight system $\Delta_{\mathcal M_2}$. Let us compute $\Delta' =\Delta'(\Delta_{\mathcal M_2})$. According (\ref{eq additional weights beta def}) in this case we need one additional weight $\beta_{21}$, since $n_1=2$. Using (\ref{eq Delta_delta new}), we have
\begin{align*}
\Delta'_{0}= \Big\{ 0 +  (\beta_{I_1}  - |I_1| \alpha_1)\,\,| \,\, I_1\subset\{2\} \,\, \text{such that}\,\,
|I_1|= 0 \,\, \text{or} \,\, |I_1|= -1	\Big\} = \{ 0\},\\
\Delta'_{\alpha_1}= \Big\{\alpha_1 + (\beta_{I_1}  - |I_1| \alpha_1)\,\,| \,\, I_1\subset\{2\} \,\, \text{such that}\,\,
|I_1|= 1 \,\, \text{or} \,\, |I_1|= 0	\Big\}\\
 = \{ \alpha_1, \beta_{21}\},
\\
\Delta'_{2\alpha_1}= \Big\{2\alpha_1 +  (\beta_{I_1}  - |I_1| \alpha_1)\,\,| \,\, I_1\subset\{2\} \,\, \text{such that}\,\,
|I_1|= 2 \,\, \text{or} \,\, |I_1|= 1	\Big\}=\\
  = \{ \alpha_1+\beta_{21}\}. 
\end{align*}
Summing up, $\Delta'= \{ 0, \alpha_1, \beta_{21}, \alpha_1+\beta_{21}\}$. Therefore $\Delta' = \Delta_{\mathbf D_2}$ is the weight system of a double vector bundle. 

To clarify the behavior of the parities in this case, we need to consider two cases.  Firstly let us take $\bar \alpha_1=\bar 0$. Recall that the de Rham differential $\d_{\beta_{21}}$ is odd and has weight $\alpha_1 -\beta_{21}$. Therefore the parity of $\beta_{21}$ is equal to $\bar 1$. Secondly, we take $\bar \alpha_1=\bar 1$. In this case the parity of $\beta_{21}$ is equal to $\bar 0$. In general $\beta_{ji}$ has always the opposite parity to $\alpha_i$.

Consider another example. Let us now compute $\Delta'$ in the case that 
$$
\Delta = \{0,\alpha_1,\alpha_2, \alpha_1+\alpha_2, 2\alpha_1+\alpha_2 \},
$$
see Section $3.5$, case $B_2$.
In this case $n_1=2$ and $n_2=1$, hence we need one additional weight $\beta_{21}$. Using (\ref{eq Delta_delta new}), we have
\begin{align*}
&\Delta'_{0} = \{ 0\},\quad \Delta'_{\alpha_1} = \{ \alpha_1, \beta_{21}\}, \quad \Delta'_{\alpha_2} = \{ \alpha_2\},\\ 
&\Delta'_{\alpha_1+\alpha_2} = \{ \alpha_1+\alpha_2,  \alpha_2+\beta_{21},  \},\quad
\Delta'_{2\alpha_1+\alpha_2} = \{ \alpha_1+\alpha_2+\beta_{21} \}.
\end{align*}
Here $\bar \beta_{21} = \bar \alpha_1 +\bar 1$.

\section{Additional structures on $\mathbf D_{\mathcal N}= \mathbb F(\mathcal N)$}

\subsection{Odd commuting vector fields on $\mathbf D_{\mathcal N}$}

Recall that we denoted by $\mathcal N$ a graded manifold of type $\Delta$ and by $\widetilde{\mathcal N}$ the $n$-times iterated tangent bundle of $\mathcal N$.  It is a graded manifold of type $\Delta_{\widetilde{\mathcal N}}$. Further, $n=r'-r$, where $r$ is the rank of $\Delta$ and $r'$ is the rank of $\Delta'$.  On  $\widetilde{\mathcal N}$  there are $n$ odd commuting homological vector fields $\d_{\beta_{ji}}$ of weights $\beta_{ji}-\alpha_i$, see Section $5.1$.  Our goal now is to show that these vector fields induce odd commuting homological vector fields on $\mathbf D_{\mathcal N}$.

\medskip
\prop\label{prop operators are defined on D_N} {\sl The de Rham differentials (\ref{eq order of de Rham diff}) defined on $\widetilde{\mathcal N}$ induce $n$ odd commuting homological vector fields on $\mathbf D_{\mathcal N}$:
	\begin{equation}\label{eq order of induced operators}
	\D_{\beta_{21}}, \ldots, \D_{\beta_{n_11}},\, \D_{\beta_{22}}, \ldots, \D_{\beta_{n_22}}, \ldots, \D_{\beta_{2r}}, \ldots, \D_{\beta_{n_rr}}.
	\end{equation}

}

\medskip

\noindent{\it Proof.}  By our weight agreement any vector field $\d_{\beta_{ji}}$ preserves the ideal $\mathcal J^-$. Hence $\d_{\beta_{ji}}$ determines the vector field $\D_{\beta_{ji}}$ acting on the sheaf $\mathcal O_{\widetilde{\mathcal N}}/\mathcal J^-$.   Furthermore, by definition we have $\mathcal O_{\mathbf D_{\mathcal N}} \subset  \mathcal O_{\widetilde{\mathcal N}}/\mathcal J^-$. We need to show that $\D_{\beta_{ji}}\big( (\mathcal O_{\mathbf D_{\mathcal N}})_{\delta}\big) \subset \mathcal O_{\mathbf D_{\mathcal N}}$ for any $\delta\in \Delta'$. 
Consider the following inclusion:
$$
\D_{\beta_{ji}}\big( (\mathcal O_{\mathbf D_{\mathcal N}})_{\delta}\big)=\D_{\beta_{ji}}\big( (\mathcal O_{\widetilde{\mathcal N}}/\mathcal J^-)_{\delta}\big) \subset (\mathcal O_{\widetilde{\mathcal N}}/\mathcal J^-)_{\delta + \beta_{ji}-\alpha_i}. 
$$
Since $\delta$ is multiplicity free, the coefficient $a_i(\delta)$ before $\alpha_i$ is equal to $0$ or $1$, see (\ref{eq delta = sum a_i alpha_i}) for notations. In case $a_i(\delta)=0$, the weight $\delta + \beta_{ji}-\alpha_i$ has a negative coefficient, hence 
$$
(\mathcal O_{\widetilde{\mathcal N}}/\mathcal J^-)_{\delta + \beta_{ji}-\alpha_i}=\{0\}.
$$ 
In case $a_i(\delta)=1$, the weight $\delta + \beta_{ji}-\alpha_i$ has no negative coefficients. Since $\beta_{ji}$ has no multiplicities in the weight system $\Delta_{\widetilde{\mathcal N}}$, the sheaf $(\mathcal O_{\widetilde{\mathcal N}}/\mathcal J^-)_{\delta + \beta_{ji}-\alpha_i}$ is a product of subsheaves in $\mathcal O_{\widetilde{\mathcal N}}/\mathcal J^-$ with multiplicity free weights. Hence it is a subsheaf in  $\mathcal O_{\mathbf D_{\mathcal N}}$. The proof is complete.$\Box$

\medskip

Some properties of the vector fields $\D_{\beta_{ji}}$ are described in the next propositions. 

\medskip

\prop\label{prop D_ji are linear} {\sl
The vector fields $\D_{\beta_{ji}}$, see (\ref{eq order of induced operators}), are $ (\mathcal O_{\mathbf D_{\mathcal N}})_0$-linear. 
}

\medskip

\noindent{\it Proof.} Let us take $f\in (\mathcal O_{\widetilde{\mathcal N}})_0$ and a vector field $\d_{\beta_{ji}}$. Then the weight of $\d_{\beta_{ji}}(f)$ is equal to $\beta_{ji}-\alpha_i$. It is a weight with a negative coefficient, therefore, 
$\D_{\beta_{ji}}(f)=0.
$
The result follows from the Leibniz rule.$\Box$

\medskip

Let $\mathcal N$ and $\mathcal N^1$ be two graded manifolds of type $\Delta$. Denote by $\D_{\beta_{ji}}$ and $\D^1_{\beta_{ji}}$ two vector fields given  on $\mathbf D_{\mathcal N}$ and on $\mathbf D_{\mathcal N^1}$, respectively, and defined  as in Proposition \ref{prop operators are defined on D_N}.

\medskip

\prop\label{prop image of a morphism commutes with operators} {\sl 
	Let $\psi: \mathcal N \to \mathcal N_1$ be a morphism of graded manifolds of type $\Delta$ and $\mathbb F(\psi): \mathbf D_{\mathcal N} \to \mathbf D_{\mathcal N_1}$ be the corresponding morphism of $r'$-fold vector bundles of type $\Delta'$. Then
	$$
	\mathbb F(\psi)^* \circ \D^1_{\beta_{ji}} = \D_{\beta_{ji}} \circ \mathbb F(\psi)^*.
	$$ 

}

\medskip

\noindent{\it Proof.} This follows from the definition of $\mathbb F(\psi)$ and the fact that all morphisms and the induced morphisms between tangent spaces commute with de Rham differentials.$\Box$

\medskip

\subsection{Description of $\mathbf D_{\mathcal N}$ in local coordinates}

Let us take $\delta \in \Delta$ and $\delta'=\sum_ia_i\alpha_i + \sum_{ji}b_{ji}\beta_{ji} \in \Delta'_{\delta}$, see (\ref{eq Delta_delta new}) for the definition of $\Delta'_{\delta}$. Then there exists the unique up to sign operator $\D_{\delta \to \delta'}$ that is equal to a composition of some $\D_{\beta_{ji}}$ or equal to the identity such that $\D_{\delta \to \delta'}(\delta) = \delta'.$ The operator $\D_{\delta \to \delta'}$ is explicitly given by  
$
\D_{\delta \to \delta'} =\pm \D_{\beta_{j_1i_1}} \circ \cdots \circ \D_{\beta_{j_ki_k}}$, where this composition is taken over all $\beta_{j_si_s}$ such that $b_{j_si_s}\ne 0$ in the expression for $\delta'$. 
Let us choose a local chart $\mathcal U$ on $\mathcal N$ with local coordinates $(\xi_i^{\delta})_{\delta \in \Delta}$. Here the superscript $\delta$ indicates the weight of coordinates. By our construction of $\mathbf D_{\mathcal N}$ we obtain the following proposition.

\medskip
\prop\label{prop coordinates on D_N}{\sl The ringed space $(\mathcal U_0, \mathcal O_{\mathbf D_{\mathcal N}}|_{\mathcal U_0})$ is a local chart on $\mathbf D_{\mathcal N}$ with the following local coordinates:
	\begin{equation}\label{eq local coordinates on D_N}
	\bigcup_{\delta\in \Delta}\{ \D_{\delta \to \delta'}(\xi_i^{\delta})\,\, | \,\, \delta' \in \Delta'_{\delta}\}.
	\end{equation}
	
}

\medskip

The coordinates (\ref{eq local coordinates on D_N}) satisfy the following property. Let  $\delta$, $\delta'$ and $\D_{\delta \to \delta'}$ be as above.
Then $\eta_{i}^{\delta'}:=\D_{\delta \to \delta'}(\xi_i^{\delta})$ satisfies the equation
$
\D_{\beta_{j_si_s}}(\eta_{i}^{\delta'})=0
$
for all $\D_{\beta_{j_si_s}}$ such that $b_{j_si_s}\ne 0$ in the expression for $\delta'$.  In other words, $\eta_{i}^{\delta'}\in \Ker \D_{\beta_{j_si_s}}$ for all such $\beta_{j_si_s}$. This means that the local coordinates of weight $\delta'$ of our chart 
$(\mathcal U_0, \mathcal O_{\mathbf D_{\mathcal N}}|_{\mathcal U_0})$ satisfy the condition $\eta_{i}^{\delta'}\in \Ker \D_{\beta_{j_si_s}}$ for all $\beta_{j_si_s}$ such that $b_{j_si_s}\ne 0$ in the expression for $\delta'$. 
Consider the sheaf $(\mathcal O_{\mathbf D_{\mathcal N}})_{\delta'}|_{\mathcal U_0}$. This sheaf is generated over $(\mathcal O_{\mathbf D_{\mathcal N}})_{0}|_{\mathcal U_0}$ by the coordinates $\eta_{i}^{\delta'}$ and by
  the sheaf 
$$
\bigoplus_{\delta'_1+\delta'_2= \delta'} (\mathcal O_{\mathbf D_{\mathcal N}})_{\delta'_1} (\mathcal O_{\mathbf D_{\mathcal N}})_{\delta'_2}|_{\mathcal U_0},
$$	
where $\delta'_i\ne \delta'$. The global version of this observation is stated in the following Proposition.

\medskip
\prop\label{prop the sheaf O_delta is generated by ker and smaller}
{\sl Let us take $\delta' = \sum\limits_i a_i \alpha_i + \sum\limits_{ji} b_{ji}\beta_{ji}\in \Delta'$.  The sheaf $(\mathcal O_{\mathbf D_{\mathcal N}})_{\delta'}$ possesses the following decomposition:
	$$
	(\mathcal O_{\mathbf D_{\mathcal N}})_{\delta'}=	\bigoplus_{\delta'_1+\delta'_2= \delta'} (\mathcal O_{\mathbf D_{\mathcal N}})_{\delta'_1} (\mathcal O_{\mathbf D_{\mathcal N}})_{\delta'_2} + (\mathcal O_{\mathbf D_{\mathcal N}})_{\delta'}\!\!\bigcap_{b_{j_si_s} \ne 0}  \Ker \D_{\beta_{j_si_s}},
	$$
	where $\delta'_i\ne \delta'$.$\Box$	
}

\medskip

\subsection{Properties of the structure sheaf of $\mathbf D_{\mathcal N}$}

Recall that $\mathcal N$ is a graded manifold of type $\Delta$, $\widetilde{\mathcal N}$ is the $n$-times iterated tangent bundle of $\mathcal N$, it is a graded manifold of type $\Delta_{\widetilde{\mathcal N}}$, and $\mathbf D_{\mathcal N}=\mathbb F(\mathcal N)$ is a graded manifold of type $\Delta'$, see Proposition \ref{prop obtain Delta_D from Delta_N} for the definition of $\Delta'$. If $\gamma$ is a certain weight in the weight lattice generated by $\alpha_i$, $\beta_{ji}$, then by definition we put
\begin{equation}\label{eq d(delta) :=}
	\d_{\beta_{ji}} (\gamma) = 
	\D_{\beta_{ji}} (\gamma) := \gamma + \beta_{ji}- \alpha_i.
\end{equation}

Let us take a subset $\Lambda= \{\gamma_1,\ldots, \gamma_s\}$ in the set (\ref{eq additional weights beta def}). Denote by $\bar\Lambda=(\gamma_1,\ldots, \gamma_s)$ the same set $\Lambda$, but with a certain order, and by $\D^{\bar\Lambda}$ the following composition: 
\begin{equation}\label{eq D_delta = circ circ}
\D^{\bar\Lambda}: \mathcal O_{\mathcal N}\hookrightarrow  \mathcal O_{\widetilde{\mathcal N}} \to \mathcal O_{\widetilde{\mathcal N}}/ \mathcal J^-,\quad  \D^{\bar\Lambda}:= \d_{\gamma_1} \circ \cdots  \circ  \d_{\gamma_s}  \,\, \mod\,\, \mathcal J^-.
\end{equation}
Note that the underlying spaces $\mathcal N_0$ and $(\mathbf D_{\mathcal N})_0$ of graded manifolds $\mathcal N$ and $\mathbf D_{\mathcal N}$, respectively,  coincide. Hence we can identify their structure sheaves $(\mathcal O_{\mathcal N})_0 = (\mathcal O_{\mathbf D_{\mathcal N}})_0$. As in the proof of Proposition \ref{prop D_ji are linear} we can see that the map of sheaves $\D^{\bar\Lambda}$ is $(\mathcal O_{\mathbf D_{\mathcal N}})_0$-linear. Hence, $\D^{\bar\Lambda}$ is a morphism of sheaves of $(\mathcal O_{\mathbf D_{\mathcal N}})_0$-modules.

Let us take a weight $\delta = \sum\limits_{i=1}^r a_i\alpha_i$, where $a_i\geq 0$, in the weight lattice generated by $\alpha_i$. (Note that $\delta$ is not necessary from $\Delta$.) Then the weight $\D^{\bar\Lambda}(\delta)$ is defined by (\ref{eq d(delta) :=}).
Let $\D^{\bar\Lambda}(\delta)$ be multiplicity free and does not have negative coefficients. Then we have the following morphism:
\begin{equation}\label{eq Psi_delta=}
\begin{split}
\D^{\bar\Lambda}: (\mathcal O_{\mathcal N})_{\delta} \to (\mathcal O_{\mathbf D_{\mathcal N}})_{\D^{\bar\Lambda}(\delta)}.
\end{split}
\end{equation}
Note that in (\ref{eq Psi_delta=}) the weight $\D^{\bar\Lambda}(\delta)$ is not necessary from $\Delta'$. However, for any multiplicity free weight $\theta$ we have $(\mathcal O_{\mathbf D_{\mathcal N}})_{\theta} =  (\mathcal O_{\widetilde{\mathcal N}}/ \mathcal J^-)_{\theta}$.

\medskip
\noindent{\bf Example.} Consider a weight system $\Delta= \{ 0, \alpha_1,3\alpha_1 \}$. In this case Set (\ref{eq additional weights beta def})  is equal to $\{\beta_{21}, \beta_{31}\}$ and $\Delta'= \{0, \alpha_1, \beta_{21}, \beta_{31}, \alpha_1+ \beta_{21}+ \beta_{31}\}$. Let us take $\delta=2\alpha_1\notin \Delta$. Let $\D^{\bar\Lambda}(\delta)$ be multiplicity free. Then for $\bar\Lambda$ we have the following three possibilities up to order of $\beta$'s.
\begin{itemize}
	\item $\bar\Lambda_1 = (\beta_{21}, \beta_{31})$ and $\D^{\bar\Lambda_1}(2\alpha_1) = 2\alpha_1 +(\beta_{31} - \alpha_1) + (\beta_{21} - \alpha_1) = \beta_{21}+ \beta_{31}$; 
	\item $\bar\Lambda_2 = (\beta_{21})$ and $\D^{\bar\Lambda_2}(2\alpha_1) = 2\alpha_1 + (\beta_{21} - \alpha_1) = \alpha_1 + \beta_{21} $; 
	\item $\bar\Lambda_3 = (\beta_{31})$ and $\D^{\bar\Lambda_3}(2\alpha_1) = 2\alpha_1 + (\beta_{31} - \alpha_1) = \alpha_1 + \beta_{31}$.
\end{itemize}
For instance we see that $\D^{\bar\Lambda}(2\alpha_1)\notin \Delta'$.

For any graded manifold $\mathcal N$ of type $\Delta$ and the corresponding $\mathbf D_{\mathcal N}=\mathbb F(\mathcal N)$, we can consider the following maps
\begin{align*}
&\D^{\bar\Lambda_1}: (\mathcal O_{\mathcal N})_{2\alpha_1} \to (\mathcal O_{\mathbf D_{\mathcal N}})_{\beta_{21}+ \beta_{31}};\quad 
\D^{\bar\Lambda_2}: (\mathcal O_{\mathcal N})_{2\alpha_1} \to (\mathcal O_{\mathbf D_{\mathcal N}})_{\alpha_1 + \beta_{21}};\\
&\D^{\bar\Lambda_3}: (\mathcal O_{\mathcal N})_{2\alpha_1} \to (\mathcal O_{\mathbf D_{\mathcal N}})_{\alpha_1 + \beta_{31}}.
\end{align*}

 Proposition \ref{prop D_delta = d circ d circ is injective} and \ref{prop D_delta is an iso onto kernel} below establish some properties of these maps. 
More precisely, we will show that the maps $\D^{\bar\Lambda_i}$, $i=1-3$, are injective and in the case $i=2,3$  we will find the image of $\D^{\bar\Lambda_i}$. For instance consider the map $\D^{\bar\Lambda_1}: (\mathcal O_{\mathcal N})_{2\alpha_1} \to (\mathcal O_{\mathbf D_{\mathcal N}})_{\beta_{21}+ \beta_{31}}$. Let as take a graded domain $\mathcal U$ on $\mathcal N$ with coordinates $(x_i,\xi^{\alpha_1}_j,\xi^{3\alpha_1}_k)$ of weights $0$, $\alpha_1$ and $3\alpha_1$, respectively. Assume that $\bar \alpha_1= \bar 1$.  The sheaf $(\mathcal O_{\mathcal N})_{2\alpha_1}$ is generated locally over $(\mathcal O_{\mathcal N})_{0}$  by the monomials $\xi^{\alpha_1}_i \cdot \xi^{\alpha_1}_j$, where $i\ne j$.  We have
\begin{align*}
&\D^{\bar\Lambda_1}(\xi^{\alpha_1}_i\cdot \xi^{\alpha_1}_j) =  \d_{\beta_{21}} \circ  \d_{\beta_{31}}(\xi^{\alpha_1}_i \cdot \xi^{\alpha_1}_j) \,\, \mod\,\, \mathcal J^- =  \d_{\beta_{21}} (\d_{\beta_{31}}(\xi^{\alpha_1}_i)\cdot \xi^{\alpha_1}_j -\\
&\xi^{\alpha_1}_i\cdot \d_{\beta_{31}}(\xi_2)) \,\, \mod\,\, \mathcal J^- =  \d_{\beta_{31}}(\xi^{\alpha_1}_i)\cdot \d_{\beta_{21}}(\xi^{\alpha_1}_j) -
\d_{\beta_{21}} (\xi^{\alpha_1}_i)\cdot \d_{\beta_{31}}(\xi^{\alpha_1}_j) \,\, \mod\,\, \mathcal J^-. 
\end{align*}
The element $\d_{\beta_{31}}(\xi^{\alpha_1}_i)\cdot \d_{\beta_{21}}(\xi^{\alpha_1}_j) -
\d_{\beta_{21}} (\xi^{\alpha_1}_i)\cdot \d_{\beta_{31}}(\xi^{\alpha_1}_j) \,\, \mod\,\, \mathcal J^-$ is not trivial in $(\mathcal O_{\mathbf D_{\mathcal N}})_{\beta_{21}+ \beta_{31}}$, hence the restriction $\D^{\bar\Lambda_1}|(\mathcal O_{\mathcal N})_{2\alpha_1}$ is injective.

\medskip

 We will need the following proposition. 

\medskip

\prop\label{prop D_delta = d circ d circ is injective} {\sl Let us take $\delta = \sum\limits_{i=1}^r a_i\alpha_i$, where $a_i\geq 0$, and assume that the weight $\D^{\bar\Lambda}(\delta)$ is multiplicity free and does not have negative coefficients. Then the morphism (\ref{eq Psi_delta=}) 
is injective.
}

\medskip

\noindent{\it Proof.} The map $\D^{\bar\Lambda}$ is a composition of de Rham differentials  $mod\, \mathcal J^-$. The idea of the proof is to use the following fact: the kernel of the de Rham differential for graded manifolds (as for usual manifolds) restricted to functions coincides with the vector space of constant functions. Then we use the fact that the sheaf $(\mathcal O_{\mathcal N})_{\delta} $ does not contain constant functions for $\delta\ne 0$.   A detailed proof can be found in Appendix.$\Box$

\medskip

The next proposition describes the image of $\D^{\bar\Lambda}$ in some particular cases.

\medskip
\prop\label{prop D_delta is an iso onto kernel} {\sl Let us take $\delta = \sum\limits_{i=1}^r a_i\alpha_i$, where $a_i\geq 0$, and assume that  $\D^{\bar\Lambda}(\delta) = \sum\limits_{i=1}^r a_i\alpha_i + \sum\limits_{ji} b_{ji}\beta_{ji}$ is multiplicity free, does not have negative coefficients and that $\D^{\bar\Lambda}(\delta)$ satisfies the following property: if $b_{st}\ne 0$, then $a_t\ne 0$. Then we have
$$
\D^{\bar\Lambda} \big((\mathcal O_{\mathcal N})_{\delta}\big) =  \Big(	(\mathcal O_{\mathbf D_{\mathcal N}})_{\D^{\bar\Lambda}(\delta)} \Big)\bigcap_{k=1}^s \Ker \D_{\gamma_{k}}
$$
 and the map 
$$
\D^{\bar\Lambda}: (\mathcal O_{\mathcal N})_{\delta} \to \Big(	(\mathcal O_{\mathbf D_{\mathcal N}})_{\D^{\bar\Lambda}(\delta)} \Big)\bigcap_{k=1}^s \Ker \D_{\gamma_{k}} 
$$
is an isomorphism. 
}

\medskip

\noindent{\it Proof.} The idea of the proof is to use the Poincar\'{e} Lemma for graded manifolds: any closed differential form is locally exact. Details can be found in Appendix.$\Box$

\medskip

Consider again the map $\D^{\bar\Lambda_1}$ as above. In this case $\D^{\bar\Lambda_1}(2\alpha_1) = \beta_{21}+ \beta_{31}$. 
 We see that $\D^{\bar\Lambda_1}(\delta)$ does not satisfy the property: if $b_{st}\ne 0$, then $a_t\ne 0$. In this case Proposition \ref{prop D_delta is an iso onto kernel} is wrong since $\Im (\D^{\bar\Lambda_1})$ does not contain for example $\d_{\beta_{21}} (\xi^{\alpha_1}_i)\cdot \d_{\beta_{31}}(\xi^{\alpha_1}_j) \,\, \mod\,\, \mathcal J^-\in \cap_{k=2}^3 \Ker \D_{\beta_{k1}}$.

We will need the following corollary:

\medskip

\noindent{\bf Corollary.} {\sl Let us take $\delta = \sum\limits_{i=1}^r a_i\alpha_i + \sum\limits_{ji} b_{ji}\beta_{ji}\in \Delta'$ such that $a_{i_0}, b_{j_oi_0}\ne 0$ for some indexes $i_0$ and $(j_0i_0)$, and  
	$$
	f\in (\mathcal O_{\mathbf D_{\mathcal N}})_{\delta}\bigcap_{b_{st} \ne 0}  \Ker \D_{\beta_{st}}.
	$$
	Then there exists $F\in (\mathcal O_{\widetilde{\mathcal N}}/\mathcal J^-)_{\delta - \beta_{j_0i_0} +\alpha_i}$
	such that 
	$$
	(\d_{\beta_{j_0i_0}} \mod \mathcal J^-)(F) =f \quad \text{and}\quad (\d_{\beta_{st}}\mod \mathcal J^-)(F)=0 
	$$
	for any $(st)\ne (j_0i_0)$ such that $b_{st} \ne 0$.$\Box$ 
	
}

\medskip

\noindent{Proof} follows from the proof of Proposition \ref{prop D_delta is an iso onto kernel}, see Appendix.$\Box$ 

\medskip

Further properties of the commuting vector fields $\D_{\beta_{ji}}$ are described in the following proposition.  

\medskip

\prop\label{prop D_ji are non-degenerate} {\sl
	Let us take $\delta \in \Delta'$. If $\D_{\beta_{ji}}(\delta)\in \Delta'$, then 
	\begin{equation}\label{eq non-degenerate operators new}
		\D_{\beta_{ji}}: (\mathcal O_{\mathbf D_{\mathcal N}} )_{\delta} \to (\mathcal O_{\mathbf D_{\mathcal N}} )_{\D_{\beta_{ji}}(\delta)}
	\end{equation}
	is an isomorphism of $(\mathcal O_{\mathbf D_{\mathcal N}} )_{0}$-locally free sheaves. In particular, all maps 
	$$
	\D_{\beta_{ji}}: (\mathcal O_{\mathbf D_{\mathcal N}} )_{\alpha_i} \to (\mathcal O_{\mathbf D_{\mathcal N}} )_{\beta_{ji}}
	$$
	are isomorphisms of $(\mathcal O_{\mathbf D_{\mathcal N}} )_{0}$-locally free sheaves. 
}

\medskip

\noindent{\it Proof.} Recall that $\D_{\beta_{ji}}$ is $(\mathcal O_{\mathbf D_{\mathcal N}} )_{0}$-linear by Proposition \ref{prop D_ji are linear}. 
Consider a chart $\mathcal U$ on $\mathcal M$. Clearly this chart determines a chart on $\mathbf D_{\mathcal N}$. We choose coordinates   $(x_p)$, $(\xi_q)$ and $(\eta_t)$ such that  $x_p$ are local coordinates of weight $0$, $\xi_q$ are coordinates with weights in the form $\alpha_{i} + \ldots$, and $\eta_t$ are other local coordinates. 

Any $f\in (\mathcal O_{\mathbf D_{\mathcal N}} )_{\delta}$ has the following form $f=\sum\limits_{kI} f_{kI} \xi_k \eta^{I}$, where $I$ as a multi-index and $f_{kI}$ are functions of weight $0$. Since $\D_{\beta_{ji}}(\delta) \in \Delta'$ is multiplicity free, we see that $\delta$ does not depend on $\beta_{ji}$. 
 The map (\ref{eq non-degenerate operators new}) in coordinates is given by the following formula:
$$
\D_{\beta_{ji}}(f) = \D_{\beta_{ji}}(\sum_{kI} f_{kI} \xi_k \eta^{I}) = \sum_{kI} f_{kI} \D_{\beta_{ji}}(\xi_k) \eta^{I}.
$$
We see that $\D_{\beta_{ji}}(\xi_k)$ and $\eta_{q}$ form a subset of independent local coordinates in $\mathbf D_{\mathcal N}$, since $\D_{\beta_{ji}}(\xi_k)$ and $\eta_{q}$ have different weights. Note that any function in $\mathcal O_{\mathbf D_{\mathcal N}} $ of weight $\D_{\beta_{ji}}(\delta)$ has the form $\sum\limits_{kI} f_{kI} \D_{\beta_{ji}}(\xi_k) \eta^{I}$. Therefore, the inverse map $\sum\limits_{kI} f_{kI} \D_{\beta_{ji}}(\xi_k) \eta^{I} \mapsto \sum\limits_{kI} f_{kI} \xi_k \eta^{I}$ of the map (\ref{eq non-degenerate operators new}) is well-defined.$\Box$

\medskip

The vector fields satisfying (\ref{eq non-degenerate operators new}) we will call {\it non-degenerate}.

\subsection{Combinatorical properties of odd commuting vector fields $\D_{\beta_{ji}}$}

Some properties of  odd commuting vector fields $\D_{\beta_{ji}}$ can be described using the combinatorics of the weight system $\Delta'$. Let us take $\delta,\delta'\in \Delta'$ and two vector fields $\D_{\beta_{ji}}$ and $\D_{\beta_{si}}$ such that $
\D_{\beta_{ji}}(\delta) = \D_{\beta_{si}}(\delta') \in \Delta'$.
Explicitly this means that $\delta = \alpha_i+\beta_{si}+\theta$ and $\delta' = \alpha_i+\beta_{ji}+\theta$ for a certain weight $\theta$. 
Then
$$
\D_{\beta_{ji}}: (\mathcal O_{\mathbf D_{\mathcal N}})_{\delta} \to (\mathcal O_{\mathbf D_{\mathcal N}})_{\D_{\beta_{ji}}(\delta)} \quad \text{and}\quad \D_{\beta_{si}}: (\mathcal O_{\mathbf D_{\mathcal N}})_{\delta'} \to (\mathcal O_{\mathbf D_{\mathcal N}})_{\D_{\beta_{ji}}(\delta)}
$$ 
are isomorphisms, see Proposition \ref{prop D_ji are non-degenerate}. Hence the following isomorphism of sheaves is defined
$$
\D_{\beta_{si}}^{-1} \circ \D_{\beta_{ji}} :  (\mathcal O_{\mathbf D_{\mathcal N}})_{\delta} \to  (\mathcal O_{\mathbf D_{\mathcal N}})_{\delta'}.
$$
Explicitly on weights we have 
$$
\D_{\beta_{si}}^{-1} \circ \D_{\beta_{ji}}(\delta) = \D_{\beta_{si}}^{-1} \circ \D_{\beta_{ji}}(\alpha_i+\beta_{si}+\theta ) = \alpha_i+\beta_{ji}+\theta = \delta'.  
$$ 
Assume that $\delta = \sum\limits_{i=1}^r a_i\alpha_i + \sum\limits_{ji} b_{ji}\beta_{ji}$ and $\delta'= 	 \sum\limits_{i=1}^r a'_i\alpha_i + \sum\limits_{ji} b'_{ji}\beta_{ji}.$ 
For $\delta\in \Delta'$ we put 
\begin{equation}\label{eq sheaf S_delta}
\mathcal S_{\delta}:= \big((\mathcal O_{\mathbf D_{\mathcal N}})_{\delta}\bigcap\limits_{b_{pt} \ne 0}  \Ker \D_{\beta_{pt}}\big).
\end{equation}
Similarly we define the sheaf $\mathcal S_{\delta'}$

We will need the following proposition.

\medskip

\prop\label{prop properties of change of a root 1}
{\sl Assume that $b_{ji}=0$ and $a_{i}, b_{si}\ne 0$ for indexes $i$, $(ji)$ and $(si)$. 
Then we have
	$$
	(\D_{\beta_{si}}^{-1} \circ \D_{\beta_{ji}}) \big(\mathcal S_{\delta}\big) = \mathcal S_{\delta'}.
	$$

}

\medskip

\noindent{\it Proof.} It is enough to show only the following inclusion
$$
(\D_{\beta_{si}}^{-1} \circ \D_{\beta_{ji}}) \big(\mathcal S_{\delta}\big) \subset \mathcal S_{\delta'}.
$$
Let us take $f\in\mathcal S_{\delta}$. In Corollary of Proposition \ref{prop D_delta is an iso onto kernel}, we have seen that since 
$\D_{\beta_{si}}(f) =0$, there exists $F\in (\mathcal O_{\widetilde{\mathcal N}}/\mathcal J^-)_{\delta - \beta_{si} +\alpha_{i}}$ 
such that $f= \d_{\beta_{si}}(F) \mod \mathcal J^-$. We have
\begin{align*}
(\D_{\beta_{si}}^{-1} \circ \D_{\beta_{ji}})(f) &= (\D_{\beta_{si}}^{-1} \circ \D_{\beta_{ji}} \circ (\d_{\beta_{si}} \mod \mathcal J^-))(F) =\\
& -  (\D_{\beta_{si}}^{-1} \circ \D_{\beta_{si}} \circ (\d_{\beta_{ji}} \mod\, \mathcal J^-))(F)=\\
&-  (\d_{\beta_{ji}} \mod\, \mathcal J^-)(F)\in (\mathcal O_{\mathbf D_{\mathcal N}})_{\delta'} \cap \Ker \D_{\beta_{ji}}.
\end{align*}
Further, again by Corollary of Proposition \ref{prop D_delta is an iso onto kernel}, we have $\d_{\beta_{pt}}(F)\, \mod \,\mathcal J^-=0$, where $b_{pt}\ne 0$.  Hence
\begin{align*}
\D_{\beta_{pt}}\circ (\D_{\beta_{si}}^{-1} \circ \D_{\beta_{ji}})(f) &= - \D_{\beta_{pt}}\circ (\d_{\beta_{ji}} \mod\, \mathcal J^-)(F)=0.
\end{align*}
The proof is complete.$\Box$

\medskip

Let us take $\delta, \delta_1, \delta_2\in \Delta'$ in the following form:
$$
\delta= \alpha_i + \beta_{ji} +\theta,\quad  \delta_1= \alpha_i + \beta_{j_1i}+\theta,\quad \delta_2= \alpha_i + \beta_{j_2i}+\theta,
$$
where $j\ne j_1$, $j\ne j_2$ and $j_1\ne j_2$. Note that since $\Delta'$ is multiplicity free, $\theta$ does not depend on $\alpha_i$, $\beta_{ji}$, $\beta_{j_1i}$ and  $\beta_{j_2i}$.

\medskip

\prop\label{prop properties of change of a root cocycle condition} {\it Let $\delta, \delta_1,\delta_2$ be as above. Then we have: 
\begin{equation}\label{eq cocycle condition formula, operators}
\begin{split}
(\D_{\beta_{ji}}^{-1} \circ \D_{\beta_{j_2i}})|_{(\mathcal O_{\mathbf D_{\mathcal N}})_{\delta} \cap \Ker \D_{\beta_{ji}}}  = - &(\D_{\beta_{j_1i}}^{-1} \circ \D_{\beta_{j_2i}}) \circ \\	
&(\D_{\beta_{ji}}^{-1} \circ \D_{\beta_{j_1i}})|_{(\mathcal O_{\mathbf D_{\mathcal N}})_{\delta} \cap \Ker \D_{\beta_{ji}}}.
\end{split}
\end{equation}

}

\medskip

We will call (\ref{eq cocycle condition formula, operators}) the {\bf cocycle like condition} or just {\bf cocycle condition} for our vector fields. This name was inspired by classical cocycle conditions. However in our case we have an additional sign and the order of entries does not agree with the classical case. 

\medskip

\noindent{\it Proof.} Let us take $f\in (\mathcal O_{\mathbf D_{\mathcal N}})_{\delta} \cap \Ker \D_{\beta_{ji}}$. Again by Corollary of Proposition \ref{prop D_delta is an iso onto kernel}, we can find $F$ such that $f= (\d_{\beta_{ji}} \mod\, \mathcal J^-)(F)$. We have
\begin{align*}
(\D_{\beta_{ji}}^{-1} \circ \D_{\beta_{j_2i}}) (\d_{\beta_{ji}} \mod \,\mathcal J^-)(F) = -  (\d_{\beta_{j_2i}} \mod\, \mathcal J^-)(F).
\end{align*}
On the other hand,
\begin{align*}
(\D_{\beta_{j_1i}}^{-1} \circ \D_{\beta_{j_2i}}) \circ 	(\D_{\beta_{ji}}^{-1} \circ \D_{\beta_{j_1i}}) (\d_{\beta_{ji}} \mod\, \mathcal J^-)(F) &=\\ 
-(\D_{\beta_{j_1i}}^{-1} \circ \D_{\beta_{j_2i}}) (\d_{\beta_{j_1i}} \mod \mathcal J^-)(F)&=
(\d_{\beta_{j_2i}} \mod\, \mathcal J^-)(F).
\end{align*}
The proof is complete.$\Box$

\medskip

\noindent{\bf Remark.} Let $\delta, \delta_1,\delta_2$ be as is Proposition \ref{prop properties of change of a root cocycle condition}. Let us show that (\ref{eq cocycle condition formula, operators}) does not hold for any $f\in (\mathcal O_{\mathbf D_{\mathcal N}})_{\delta}$. In other words the assumption $f\in (\mathcal O_{\mathbf D_{\mathcal N}})_{\delta} \cap \Ker \D_{\beta_{ji}}$ is essential.  Let us take two variables $\xi_1,\xi_2$ of weight $\alpha_i$. Then $f=\xi_1\cdot \D_{\beta_{ji}}(\xi_2)$ has the weight $\alpha_i + \beta_{ji}$.
Further, $\D_{\beta_{ji}}(f) = \D_{\beta_{ji}}(\xi_1)\cdot \D_{\beta_{ji}}(\xi_2)\ne 0$. Hence, $f\notin \Ker \D_{\beta_{ji}}$.
 Applying the left hand side of (\ref{eq cocycle condition formula, operators}), we get:
$$
(\D_{\beta_{ji}}^{-1} \circ \D_{\beta_{j_2i}})(\xi_1\cdot \D_{\beta_{ji}}(\xi_2)) = \D_{\beta_{ji}}^{-1} (\D_{\beta_{j_2i}}(\xi_1)\cdot \D_{\beta_{ji}}(\xi_2)) = \pm \D_{\beta_{j_2i}}(\xi_1)\cdot \xi_2.
$$
Further,
\begin{align*}
(\D_{\beta_{j_1i}}^{-1} \circ \D_{\beta_{j_2i}}) \circ (\D_{\beta_{ji}}^{-1} \circ &\D_{\beta_{j_1i}}) ( \xi_1\cdot \D_{\beta_{ji}}(\xi_2)) = \\
&\pm (\D_{\beta_{j_1i}}^{-1} \circ \D_{\beta_{j_2i}}) ( \D_{\beta_{j_1i}} (\xi_1)\cdot \xi_2 )=
\pm \xi_1\cdot \D_{\beta_{j_2i}}(\xi_2).
\end{align*}
We see that the results are different.

\section{Equivalence of categories}

\subsection{The category of $r'$-fold vector bundles with $n$ odd commuting non-degenerate vector fields} 

In this section we introduce the category $\Delta'$\textsf{VBVect}. This is a category of $r'$-fold vector bundles of type $\Delta'$ with $n$ odd commuting non-degenerate vector fields. More precisely, let $\Delta'$ be a  weight system  with the following set of basic weights: 
$$
\{ \alpha_i, \,\,\beta_{ji}\,\,|\,\, i=1,\ldots, r,\, j=2,\ldots, n_{i}\},
$$
where $n_i \geq 2$ and $i=1,\ldots, r$ are some non-negative integers. (See (\ref{eq Delta as a subset}) for the definition of basic weights.) 
We put $n:= \sum\limits_{i=1}^rn_i-r$ and $r'= n+r$. Note that $r'$ is the rank of $\Delta'$.  Let $\mathbf D$ be an $r'$-fold vector bundle of type $\Delta'$ 
with $n$ odd  vector fields $\D_{\beta_{ji}}$ of weights $\beta_{ji}-\alpha_i$. Assume that these {\bf vector fields} have the following {\bf properties}:
\begin{enumerate}
	\item The vector fields $\D_{\beta_{ji}}$ are {\bf $(\mathcal O_{\mathbf D})_0$-linear}.
	
	\item The vector fields $\D_{\beta_{ji}}$ {\bf super-commute}:
	$$
	[\D_{\beta_{ji}}, \D_{\beta_{j'i'}}] =0
	$$ 
	for all $(ji)$ and $(j'i')$. In particular, any $\D_{\beta_{ji}}$ satisfy the condition $\D^2_{\beta_{ji}}=0$.  
	
	\item The operators $\D_{\beta_{ji}}$ are {\bf non-degenerate} in the following sense. Let us take $\delta\in \Delta'$. As above we put $\D_{\beta_{ji}}(\delta): = \delta + \beta_{ji}-\alpha_i.$
	We call an odd vector field $\D_{\beta_{ji}}$ of weight $\beta_{ji}-\alpha_i$ {\bf non-degenerate}, if it satisfies conditions of Proposition \ref{prop D_ji are non-degenerate} for any $\delta$. More precisely, if $\D_{\beta_{ji}}(\delta)\in \Delta'$ for a certain $\delta \in \Delta'$, then the following map 
	$$
	\D_{\beta_{ji}}: (\mathcal O_{\mathbf D} )_{\delta} \to (\mathcal O_{\mathbf D} )_{\D_{\beta_{ji}}(\delta)}
	$$
	is an isomorphism of sheaves of $(\mathcal O_{\mathbf D})_0$-modules.

	\item Let us take $\delta = \sum\limits_i a_i\alpha_i + \sum\limits_{ji} b_{ji}\beta_{ji}\in \Delta'$. We assume that the sheaf $(\mathcal O_{\mathbf D})_{\delta}$ possesses the following {\bf decomposition}:
	$$
	(\mathcal O_{\mathbf D})_{\delta}= (\mathcal O_{\mathbf D})_{\delta}\bigcap_{b_{st} \ne 0}  \Ker \D_{\beta_{st}}+	\bigoplus_{\delta_1+\delta_2= \delta} (\mathcal O_{\mathbf D})_{\delta_1} (\mathcal O_{\mathbf D})_{\delta_2},
	$$
	where $\delta_1, \delta_2\ne 0$.
	
	\item Let $\delta$, $\D_{\beta_{ji}}$, $\D_{\beta_{j_2i}}$ and $\D_{\beta_{j_1i}}$ be as in Proposition \ref{prop properties of change of a root cocycle condition}. The vector fields $\D_{\beta_{ji}}$, $\D_{\beta_{j_2i}}$ and $\D_{\beta_{j_1i}}$ satisfy the following {\bf cocycle condition}:
	\begin{align*}
		(\D_{\beta_{ji}}^{-1} \circ \D_{\beta_{j_2i}})|_{(\mathcal O_{\mathbf D})_{\delta} \cap \Ker \D_{\beta_{ji}}}  = - (\D_{\beta_{j_1i}}^{-1} \circ \D_{\beta_{j_2i}}) \circ 
		(\D_{\beta_{ji}}^{-1} \circ \D_{\beta_{j_1i}})|_{(\mathcal O_{\mathbf D})_{\delta} \cap \Ker \D_{\beta_{ji}}}.
	\end{align*}

	\item Let $\delta, \delta'\in \Delta'$, $\D_{\beta_{ji_0}}$ and $\D_{\beta_{j_0i_0}}$ be as in Proposition \ref{prop properties of change of a root 1}.  Our vector fields 	{\bf preserve the kernels} in the following sense:
	$$
	(\D_{\beta_{j_0i_0}}^{-1} \circ \D_{\beta_{ji_0}}) \big((\mathcal O_{\mathbf D})_{\delta}\bigcap_{b_{st} \ne 0}  \Ker \D_{\beta_{st}}\big) = (\mathcal O_{\mathbf D})_{\delta'}\bigcap_{b'_{st}\ne 0}  \Ker \D_{\beta_{st}}.
	$$
	In other words this means that the operator $\D_{\beta_{j_0i_0}}^{-1} \circ \D_{\beta_{ji_0}}$ preserves the decomposition from item $4$.

\end{enumerate}

The category of $r'$-fold vector bundles of type $\Delta'$ with $n$ odd  vector fields of weight $\beta_{ji}-\alpha_i$ satisfying Properties $1-6$ we  denote by $\Delta'_{(\bar \alpha_1, \ldots, \bar \alpha_r)}$\textsf{VBVect} or just by $\Delta'$\textsf{VBVect}.  A morphism in this category is a morphism in the category of $r'$-fold vector bundles of type $\Delta'$ that commutes with all vector fields.

It follows from Propositions \ref{prop the sheaf O_delta is generated by ker and smaller}, \ref{prop D_ji are linear}, \ref{prop image of a morphism commutes with operators}, \ref{prop D_ji are non-degenerate}, \ref{prop properties of change of a root 1} and \ref{prop properties of change of a root cocycle condition} that the image of the functor $\mathbb F$ is contained in $\Delta'$\textsf{VBVect}. 
In the next sections we will prove that $\mathbb F$ defines an equivalence of categories. To do this we will use the following definition.

\medskip

\de\label{de equivalence of categories} Two categories $\mathcal C$ and $\mathcal C'$ are called {\it equivalent} if there is a functor $F:\mathcal C \to \mathcal C'$ such that:
\begin{itemize}
	\item $F$ is full and faithful,  this is $Hom_{\mathcal C}(c_1,c_2)$ is in bijection with \newline $Hom_{\mathcal C'}(Fc_1, Fc_2)$.
	\item $F$ is essentially surjective, this is for any $a\in \mathcal C'$ there exists $b\in \mathcal C$ such that $a$ is isomorphic to $F(b)$. 
\end{itemize}

\medskip

\noindent{\bf Example.} Let us illustrate our construction of the functor $\mathbb F$ and Properties $1-6$ on an example.  Consider a {\bf $\mathbb Z$-graded manifold $\mathcal M$ of degree $3$}, that is a graded manifold of type $\Delta_{\mathcal M_3}$. In this case $r=1$, $n_1=3$ and $n=n_1-r= 2$. Therefore we need to take  twice iterated tangent bundle $TT(\mathcal M) = T[\beta_{31}-\alpha_1](T[\beta_{21}-\alpha_1](\mathcal M))$ and to use two additional weights $\beta_{21}$ and $\beta_{31}$. As it was noticed above our construction works for both parity agreements: $\bar\alpha_1 =\bar 0$ or $\bar\alpha_1 =\bar 1$. Recall that $\bar \beta_{21} = \bar\alpha_1+\bar 1$ and $\bar \beta_{31} = \bar\alpha_1+\bar 1$. 

Let us explicitly describe $\mathbb F(\mathcal M)= \mathbf D_{\mathcal M}$. Consider a local chart on $\mathcal M$ with coordinates $(x,\xi^{\alpha_1},\xi^{2\alpha_1}, \xi^{3\alpha_1})$. We omit here all subscripts. As above a superscript indicates the weight of a coordinate. Recall that we denoted by $\d_{\beta_{21}}$ and $\d_{\beta_{31}}$ the first and the second de Rham differentials in $TT(\mathcal M)$. Therefore the standard local coordinates on $TT(\mathcal M)$ have the following form
\begin{equation}\label{eq loc coord TTM}
\begin{split}
(x,\xi^{\alpha_1},\xi^{2\alpha_1}, \xi^{3\alpha_1}, \d_{\beta_{21}}\!x, \d_{\beta_{21}}\!\xi^{\alpha_1}, \d_{\beta_{21}}\!\xi^{2\alpha_1}, \d_{\beta_{21}}\! \xi^{3\alpha_1},
\d_{\beta_{31}}\! x, \d_{\beta_{31}}\!\xi^{\alpha_1},\\ \d_{\beta_{31}}\!\xi^{2\alpha_1}, \d_{\beta_{31}}\!\xi^{3\alpha_1}, \d_{\beta_{31}}\!\d_{\beta_{21}}\!x, \d_{\beta_{31}}\!\d_{\beta_{21}}\!\xi^{\alpha_1}, \d_{\beta_{31}}\!\d_{\beta_{21}}\!\xi^{2\alpha_1}, \d_{\beta_{31}}\!\d_{\beta_{21}}\! \xi^{3\alpha_1}).
\end{split}
\end{equation}
 The first and the second de Rham differentials $\d_{\beta_{21}}$ and $\d_{\beta_{31}}$ have weights $\beta_{21} -\alpha_1$ and $\beta_{31} -\alpha_1$, respectively. The coordinates (\ref{eq loc coord TTM}) have the following weights, respectively.
\begin{equation}\label{eq weights loc coord TTM}
\begin{split}
(0, \alpha_1, 2\alpha_1, 3\alpha_1, \beta_{21}-\alpha_1, \beta_{21}, \alpha_1+\beta_{21}, 2\alpha_1+\beta_{21}, \beta_{31} -\alpha_1, \beta_{31},\\ \beta_{31} +\alpha_1, \beta_{31} +2\alpha_1, \beta_{21}+  \beta_{31}-2\alpha_1,  \beta_{31} - \alpha_1 + \beta_{21},\\ \beta_{31} +\beta_{21}, \beta_{31}+\alpha_1+\beta_{21}).
\end{split}
\end{equation}
According Section $5.2$ to obtain $\mathbf D_{\mathcal M}$ we need to factorize the structure sheaf $\mathcal O_{TT(\mathcal M)}$ of $TT(\mathcal M)$ by the ideal $\mathcal J^-$, see Section $5.1$. Recall that this ideal is generated by all local coordinates of weights with at least one negative coefficient. As it was noticed above the sheaf $\mathcal O_{TT(\mathcal M)}/ \mathcal J^-$ is a structure sheaf of a graded manifold, say $\mathcal M'$, of type
\begin{align*}
\tilde\Delta=\{ 0, \alpha_1, 2\alpha_1, 3\alpha_1,  \beta_{21}, \alpha_1+\beta_{21}, 2\alpha_1+\beta_{21}, \beta_{31}, \beta_{31} +\alpha_1,\\ \beta_{31} +2\alpha_1, \beta_{31} +\beta_{21}, \beta_{31}+\alpha_1+\beta_{21}\}
\end{align*}
To obtain $\tilde\Delta$ we removed all weights with at least one negative coefficient in (\ref{eq weights loc coord TTM}). The corresponding to (\ref{eq loc coord TTM}) local coordinates on $\mathcal M'$ have the following form 
\begin{align*}
(x,\xi^{\alpha_1},\xi^{2\alpha_1}, \xi^{3\alpha_1},  \d_{\beta_{21}}\!\xi^{\alpha_1}, \d_{\beta_{21}}\!\xi^{2\alpha_1}, \d_{\beta_{21}}\! \xi^{3\alpha_1},
 \d_{\beta_{31}}\!\xi^{\alpha_1}, \d_{\beta_{31}}\!\xi^{2\alpha_1},\\ \d_{\beta_{31}}\!\xi^{3\alpha_1},   \d_{\beta_{31}}\!\d_{\beta_{21}}\!\xi^{2\alpha_1}, \d_{\beta_{31}}\!\d_{\beta_{21}}\! \xi^{3\alpha_1}).
\end{align*}
Note that more precisely we should write $x + \mathcal J^-$, $\xi^{\alpha_1}+ \mathcal J^-$ and so on. We omit $\mathcal J^-$ for notational simplicity.

Further, by definition $\Delta'$ is the maximal multiplicity free subsystem in $\tilde\Delta$. Explicitly we have
\begin{align*}
\Delta' = \{ 0, \alpha_1, \beta_{21}, \alpha_1+\beta_{21},  \beta_{31}, \beta_{31} +\alpha_1, \beta_{31} +\beta_{21}, \beta_{31}+\alpha_1+\beta_{21}\}.
\end{align*}
The graded manifold $\mathbf D_{\mathcal M}$ is a graded manifold of type $\Delta'$. Locally its structure sheaf $\mathcal O_{\mathbf D_{\mathcal M}}$ is generated by the coordinates
\begin{align*}
(x,\xi^{\alpha_1}, \d_{\beta_{21}}\!\xi^{\alpha_1}, \d_{\beta_{21}}\!\xi^{2\alpha_1}, 
\d_{\beta_{31}}\!\xi^{\alpha_1}, \d_{\beta_{31}}\!\xi^{2\alpha_1},   \d_{\beta_{31}}\!\d_{\beta_{21}}\!\xi^{2\alpha_1}, \d_{\beta_{31}}\!\d_{\beta_{21}}\! \xi^{3\alpha_1}).
\end{align*}
To obtain the transition functions between two local charts we need to write transition functions for $TT(\mathcal M)$ and factorize by $ \mathcal J^-$. This completes the construction of $\mathbf D_{\mathcal M}$.

Note that $\mathbf D_{\mathcal M}$ is a graded manifold of type $\Delta'$, where $\Delta'$ is multiplicity free of rank $3$. Therefore $\mathbf D_{\mathcal M}$ is a $3$-fold vector bundle of type $\Delta'$. Moreover we have two operators $\D_{\beta_{21}}$ and  $\D_{\beta_{31}}$ on $\mathbf D_{\mathcal M}$ that by Propositions \ref{prop the sheaf O_delta is generated by ker and smaller}, \ref{prop D_ji are linear}, \ref{prop image of a morphism commutes with operators}, \ref{prop D_ji are non-degenerate}, \ref{prop properties of change of a root 1} and \ref{prop properties of change of a root cocycle condition}, satisfy Properies $1-6$. Let us describe these operators explicitly in our coordinates.

First of all $\D_{\beta_{21}}$ and  $\D_{\beta_{31}}$ are induced by the first and by the second de Rham differentials $\d_{\beta_{21}}$ and $\d_{\beta_{31}}$, respectively. To define for example $\D_{\beta_{21}}$ in our coordinates we need to apply  $\d_{\beta_{21}}$ to a coordinate and to factorize by $\mathcal J^-$. For example we have
\begin{align*}
\D_{\beta_{21}} (x) = \d_{\beta_{21}}\!x + \mathcal J^- = \mathcal J^-,\quad
\D_{\beta_{21}} (\d_{\beta_{21}}\!\xi^{\alpha_1} )  = \d_{\beta_{21}}  \d_{\beta_{21}}\!\xi^{\alpha_1} + \mathcal J^- = \mathcal J^-,\\
\D_{\beta_{21}} (\d_{\beta_{31}}\!\xi^{2\alpha_1} ) = \d_{\beta_{21}} \d_{\beta_{31}}\!\xi^{2\alpha_1} + \mathcal J^- = - \d_{\beta_{31}} \d_{\beta_{21}}\!\xi^{2\alpha_1} + \mathcal J^-.
\end{align*}
We see that $\D_{\beta_{21}}$ is $(\mathcal O_{\mathbf D_{\mathcal M}})_0$-linear, since it sends functions of degree $0$, for example $x$, to $0$. The operators $\D_{\beta_{21}}$ and  $\D_{\beta_{31}}$ super-commute since $\d_{\beta_{21}}$ and $\d_{\beta_{31}}$ are super-commutative. Hence we have Properties $1$ and $2$.

Further for example the sheaves $(\mathcal O_{\mathbf D_{\mathcal M}})_{\alpha_1+\beta_{31}}$ and  
$(\mathcal O_{\mathbf D_{\mathcal M}})_{\beta_{21}+\beta_{31}}$ are locally generated over $(\mathcal O_{\mathbf D_{\mathcal M}})_0$ by the monomials $\{\xi^{\alpha_1} \cdot \d_{\beta_{31}}\!\xi^{\alpha_1},\, \d_{\beta_{31}}\!\xi^{2\alpha_1}\}$ and by the monomials $\{\d_{\beta_{21}}\!\xi^{\alpha_1} \cdot \d_{\beta_{31}}\!\xi^{\alpha_1},\,\, \d_{\beta_{21}}\!\d_{\beta_{31}}\!\xi^{2\alpha_1}\}$, respectively. We have 
$$
\D_{\beta_{21}} (\xi^{\alpha_1} \cdot \d_{\beta_{31}}\!\xi^{\alpha_1}) = \d_{\beta_{21}}\!\xi^{\alpha_1} \cdot \d_{\beta_{31}}\!\xi^{\alpha_1},\quad 
\D_{\beta_{21}} (\d_{\beta_{31}}\!\xi^{2\alpha_1}) = \d_{\beta_{21}}\!\d_{\beta_{31}}\!\xi^{2\alpha_1}. 
$$
Therefore, $\D_{\beta_{21}}: (\mathcal O_{\mathbf D_{\mathcal M}})_{\alpha_1+\beta_{31}} \to  (\mathcal O_{\mathbf D_{\mathcal M}})_{\beta_{21}+\beta_{31}}$ is an isomorphism. This observation leads to Property $3$.

Consider again the local generators $\{\xi^{\alpha_1} \cdot \d_{\beta_{31}}\!\xi^{\alpha_1},\, \d_{\beta_{31}}\!\xi^{2\alpha_1}\}$ of $(\mathcal O_{\mathbf D_{\mathcal M}})_{\alpha_1+\beta_{31}}$. We see that the first monomial $\xi^{\alpha_1} \cdot \d_{\beta_{31}}\!\xi^{\alpha_1}$ 
is decomposable, while for the second we have $\d_{\beta_{31}}\!\xi^{2\alpha_1}\in \Ker \D_{\beta_{31}}$. So we get Property $4$ for this sheaf. Properties $5$ and $6$ are technical, we can check them directly in our coordinates.

\medskip
\noindent {\bf Remark.} The main result of our paper is that the graded manifold $\mathbf D_{\mathcal M}$ contains all information about the original graded manifold $\mathcal M$. Moreover if a $3$-fold vector bundle $\mathbf D_{\mathcal M}$ with operators $\D_{\beta_{21}}$ and  $\D_{\beta_{31}}$, satisfying Properties $1-6$, is given, we can recover the graded manifold $\mathcal M$. In this case the operators $\D_{\beta_{21}}$ and  $\D_{\beta_{31}}$ are images after a factorization of 
the first and the second de Rham differentials $\d_{\beta_{21}}$ and $\d_{\beta_{31}}$, respectively. 
\medskip

\subsection{\bf $\mathbb{Z}$-graded manifolds of degree $2$ and double vector bundles with an odd homological vector field}

 In this section we establish a correspondence between $\mathbb{Z}$-graded manifolds of degree $2$ and double vector bundles with an odd non-degenerate homological vector field.  Recall that graded manifolds of type $\{0,\alpha,2\alpha \}$ are usually called in the literature  {\it $\mathbb{Z}$-graded manifolds of degree $2$}. Below we give two constructions. First of all we assign a double vector bundle to a $\mathbb{Z}$-graded manifold of degree $2$ and then we reconstruct a graded manifold corresponding to a double vector bundle with an odd non-degenerate vector field.

 Constructions $1$ and $2$ in what follows establish an equivalence between the category of $\mathbb{Z}$-graded manifolds of degree $2$ and the category of double vector with an odd homological vector field.
 Another result of this type about the equivalence of categories of $\mathbb{Z}$-graded manifolds of degree $2$ and of double vector bundles with an involution (or a metric) was obtained in \cite{Fernando} (in \cite{JL}).

\bigskip

\noindent{\bf Construction 1.} Consider a $\mathbb{Z}$-graded manifold $\mathcal M_2$ of degree $2$ or, in other words, a  graded manifold $\mathcal M_2$ of type $\Delta = \{0,\alpha,2\alpha \}$. In this case $\mathbb F(\mathcal M_2) = :\mathbf D_{\mathcal M_2}$ is a double vector bundle with basic weights $\alpha:=\alpha_1$ and $\beta:=\beta_{21}$. The weight system of  $\mathbf D_{\mathcal M_2}$ has the following form:
$$
\Delta' =\{0,\,\alpha,\,\, \beta, \,\, \alpha+\beta\}.
$$
On $\mathbf D_{\mathcal M_2}$ we have an odd linear homological vector field $\D_{\beta} :=\d_{\beta} \mod \mathcal J^-$
such that $\D_{\beta}$ is {\bf non-degenerate} and has {\bf weight} $\beta- \alpha$. In this case the non-degeneracy of $\D_{\beta}$ means that the following map 
$$
\D_{\beta} : (\mathcal O_{\mathbf D_{\mathcal M_2}})_{\alpha} \longrightarrow 
(\mathcal O_{\mathbf D_{\mathcal M_2}})_{\beta}
$$
is an isomorphism of sheaves of $(\mathcal O_{\mathbf D_{\mathcal M_2}})_{0}$-modules.

\bigskip

\noindent{\bf Construction 2.}  Let us show that  any  double vector bundle $\mathbf D$ or a graded manifold of type $\Delta' = \{0,\,\alpha,\,\, \beta, \,\, \alpha+\beta\}$ with an odd non-degenerate linear homological vector field $\D$ of weight $\beta - \alpha$ is isomorphic to a double vector bundle in the form $\mathbf D_{\mathcal M_2}:=\mathbb F(\mathcal M_2)$, where $\mathcal M_2$ is a certain graded manifold of type $\Delta= \{0,\,\alpha,\,\, 2\alpha\}$. We also will show that this isomorphism commutes with operators $\D_{\beta}$ and $\D$, which are defined on $\mathbf D_{\mathcal M_2}$ and $\mathbf D$, respectively.

\medskip

\noindent{\bf Step 1, exact sequence.} Consider the subsheaf $(\mathcal O_{\mathbf D})_{\alpha +\beta}$ in $\mathcal O_{\mathbf D}$, where $\mathcal O_{\mathbf D}$ is the structure sheaf of $\mathbf D$. 
We have the following  exact sequence of sheaves of $(\mathcal O_{\mathbf D})_0$-modules:
$$
0\to (\mathcal O_{\mathbf D})_{\alpha} (\mathcal O_{\mathbf D})_{\beta}
\longrightarrow (\mathcal O_{\mathbf D})_{\alpha +\beta} \longrightarrow \mathcal E\to 0.
$$
 Here $\mathcal E$ is a certain locally free sheaf of $(\mathcal O_{\mathbf D})_0$-modules. A standard argument shows that the following sequence is also exact
$$
0\to \Ker \D \cap \big((\mathcal O_{\mathbf D})_{\alpha} (\mathcal O_{\mathbf D})_{\beta}\big)
\longrightarrow \Ker \D \cap(\mathcal O_{\mathbf D})_{\alpha +\beta} \longrightarrow \Ker \D' \to 0,
$$
where 
$$
\D':\mathcal E\to  \D((\mathcal O_{\mathbf D})_{\alpha +\beta})/ \D((\mathcal O_{\mathbf D})_{\alpha} (\mathcal O_{\mathbf D})_{\beta})
$$
is the map induced by $\D$. Since $\D$ is non-degenerate and $\D ((\mathcal O_{\mathbf D})_{\beta}) =\D^2 ((\mathcal O_{\mathbf D})_{\alpha}) = 0$, the following map 
$$
\D: (\mathcal O_{\mathbf D})_{\alpha} (\mathcal O_{\mathbf D})_{\beta} \to (\mathcal O_{\mathbf D})_{\beta} (\mathcal O_{\mathbf D})_{\beta}
$$
is surjective. Since $\D$ has weight $\beta-\alpha$ and since $2\beta\notin \Delta'$, we have
$$
\D\big( (\mathcal O_{\mathbf D})_{\alpha +\beta} \big) \subset (\mathcal O_{\mathbf D})_{2\beta} \quad \text{and} \quad (\mathcal O_{\mathbf D})_{2\beta}= (\mathcal O_{\mathbf D})_{\beta} (\mathcal O_{\mathbf D})_{\beta}.
$$
Hence, 
$$
\D((\mathcal O_{\mathbf D})_{\alpha +\beta}) = (\mathcal O_{\mathbf D})_{\beta} (\mathcal O_{\mathbf D})_{\beta},
$$
the map $\D'$ is trivial and $\Ker \D' =\mathcal E$. 
Therefore we have the following exact sequence:
\begin{equation}\label{eq exact sequence case 2}
0\to \Ker\D \cap \big( (\mathcal O_{\mathbf D})_{\alpha} (\mathcal O_{\mathbf D})_{\beta} \big)
\longrightarrow \Ker\D \cap (\mathcal O_{\mathbf D})_{\alpha +\beta} \longrightarrow \mathcal E\to 0.
\end{equation}

\medskip

\noindent{\bf Step 2, construction of $\mathcal M_2$.} The idea is to show that the following data:
\begin{align*}
	\mathcal O_{0} =  (\mathcal O_{\mathbf D})_{0}, \quad 
	\mathcal O_{\alpha} =  (\mathcal O_{\mathbf D})_{\alpha}, \quad \mathcal O_{2\alpha} =  (\mathcal O_{\mathbf D})_{\alpha +\beta}\cap \Ker \D
\end{align*}
defines a $\mathbb{Z}$-graded manifold $\mathcal M_2$ of degree $2$. To simplify notations we denoted here by $\mathcal O$ the structure sheaf of $\mathcal M_2$.  First of all consider the graded manifold $\mathcal M_1$ of type $\{0,\alpha\}$ with the structure sheaf $S^*_{\mathcal O_{0}}(\mathcal O_{\alpha})$ and the sheaf $\mathcal O_{T\mathcal M_1}/ \mathcal J^-$ with the vector field $\d \mod \mathcal J^-:= \d_{\beta} \mod \mathcal J^-$ defined as above, i.e induced by the de Rham differential on $\mathcal O_{T\mathcal M_1}$. We set $\mathcal O_{\beta} := (\d \mod \mathcal J^-)(\mathcal O_{\alpha}).
$
By Proposition \ref{prop D_delta is an iso onto kernel} or by a direct computation, we have
\begin{equation}\label{eq iso on kernel case 2}
(\d \mod \mathcal J^-)(\mathcal O_{\alpha}\cdot \mathcal O_{\alpha}) = (\mathcal O_{\alpha}\cdot \mathcal O_{\beta}) \cap \Ker (\d \mod \mathcal J^-).
\end{equation}
Now we can define an isomorphism of sheaves of ringed spaces 
$$
\Theta: \mathcal O_{T\mathcal M_1}/ \mathcal J^-= S^*(\mathcal O_{\alpha}\oplus \mathcal O_{\beta}) \longrightarrow
S^*((\mathcal O_{\mathbf D})_{\alpha}\oplus (\mathcal O_{\mathbf D})_{\beta}) 
$$
 in the following way:
$$
\Theta|_{\mathcal O_{\alpha}} := \id, \quad \Theta|_{\mathcal O_{\beta}}: = \D\circ (\d \mod \mathcal J^-)^{-1}. 
$$
Clearly, $\Theta$ preserves all weights and $\Theta \circ (\d \mod \mathcal J^-) = \D\circ\, \Theta.
$ Therefore,
\begin{equation}\label{eq iso between kernels case 2}
\Theta \big((\mathcal O_{\alpha}\cdot \mathcal O_{\beta}) \cap \Ker (\d \mod \mathcal J^-)\big) = \big( (\mathcal O_{\mathbf D})_{\alpha} (\mathcal O_{\mathbf D})_{\beta} \big) \cap \Ker\D.
\end{equation}
Combining (\ref{eq exact sequence case 2}), (\ref{eq iso on kernel case 2}) and (\ref{eq iso between kernels case 2}), we get the following exact sequence:
\begin{equation}\label{eq exact sequence for M_2 case 2}
0\to \mathcal O_{\alpha}\cdot \mathcal O_{\alpha}
\stackrel{\chi}{\longrightarrow} (\mathcal O_{\mathbf D})_{\alpha +\beta}\cap \Ker \D \longrightarrow \mathcal E\to 0.
\end{equation}
Here $\chi = \Theta\circ (\d \mod \mathcal J^-)$ is an injective map.  Further, we put 
$$
\mathcal O_{2\alpha} :=  (\mathcal O_{\mathbf D})_{\alpha +\beta}\cap \Ker \D.
$$
By Construction $2$, Section $4.2$, the exact sequence (\ref{eq exact sequence for M_2 case 2}) determines a $\mathbb{Z}$-graded manifold of degree $2$,  which we denote by $\mathcal M_2$.

\medskip

\noindent{\bf Step 3, construction of an isomorphism $\mathbb F(\mathcal M_2)\simeq \mathbf D$.} 
To define an isomorphism $\mathbb F(\mathcal M_2)\simeq \mathbf D$ we use Proposition \ref{prop isomorphism of graded mnf}. For simplicity of notations we denote the structure sheaf of $\mathbb F(\mathcal M_2)$ by $\mathcal O'$. By definition and by properties of the functor $\mathbb F$ we have
$$
\mathcal O'_{\alpha} = \mathcal O_{\alpha}, \quad \mathcal O'_{\beta} = \mathcal O_{\beta}, \quad \mathcal O'_{\alpha+\beta} = \mathcal O'_{\alpha}\cdot \mathcal O'_{\beta} +  (\d \mod \mathcal J^-) (\mathcal O_{2\alpha}) .
$$
The last equality follows from Proposition \ref{prop the sheaf O_delta is generated by ker and smaller}. We will use this decomposition to define an isomorphism $\mathbb F(\mathcal M_2)\to \mathbf D$.

Consider the following commutative diagram:
$$
\begin{CD}
(\d \mod \mathcal J^-) (\mathcal O_{2\alpha}) @<{(\d \mod \mathcal J^-)}<< \mathcal O_{2\alpha} @= (\mathcal O_{\mathbf D})_{\alpha +\beta}\cap \Ker \D\\
@AAA @AAA @AAA\\
(\d \mod \mathcal J^-) (\mathcal O_{\alpha}\cdot \mathcal O_{\alpha}) @<{(\d \mod \mathcal J^-)}<< \mathcal O_{\alpha}\cdot \mathcal O_{\alpha} @>{\chi}>>  \big( (\mathcal O_{\mathbf D})_{\alpha} (\mathcal O_{\mathbf D})_{\beta} \big)\cap \Ker \D\\
\end{CD}.
$$
Note that the right square is commutative by definition. Here all horizontal maps are isomorphisms and all vertical maps are inclusions. By Proposition \ref{prop isomorphism of graded mnf}, we need to define isomorphisms $\varphi_{\delta}$, where $\delta \in \Delta'$. We put
\begin{align*}
\varphi_{\alpha}|_{\mathcal O'_{\alpha}} = \Theta|_{\mathcal O_{\alpha}}, \quad \varphi_{\beta}|_{\mathcal O'_{\alpha}} = \Theta|_{\mathcal O_{\beta}},\quad
 \varphi_{\alpha+\beta}|_{\mathcal O'_{\alpha}\mathcal O'_{\beta}} = \varphi_{\alpha} \varphi_{\beta}|_{\mathcal O'_{\alpha}\mathcal O'_{\beta}},
 \\
 \varphi_{\alpha+\beta}|_{(\d \mod \mathcal J^-) (\mathcal O_{2\alpha})} = (\d \mod \mathcal J^-)^{-1}.
\end{align*}
In the last line we use the identification $\mathcal O_{2\alpha} =  (\mathcal O_{\mathbf D})_{\alpha +\beta}\cap \Ker \D$.
Since the diagram above is commutative, the conditions of Proposition \ref{prop isomorphism of graded mnf} hold and the injective map is defined. Since the sequence (\ref{eq exact sequence case 2}) is exact, $(\mathcal O_{\mathbf D})_{\alpha +\beta} = (\mathcal O_{\mathbf D})_{\alpha +\beta}\cap \Ker \D + (\mathcal O_{\mathbf D})_{\alpha}(\mathcal O_{\mathbf D})_{\beta}$.  Hence the map $\varphi$ defined by $(\varphi_{\delta})$ is an isomorphism.  Clearly $\varphi$ commutes with $\D_{\beta}=(\d \mod \mathcal J^-)$ and  $\D$.

\medskip

\noindent{\bf Remark.} We have seen that the decomposition $(\mathcal O_{\mathbf D})_{\alpha +\beta} = (\mathcal O_{\mathbf D})_{\alpha +\beta}\cap \Ker \D + (\mathcal O_{\mathbf D})_{\alpha}(\mathcal O_{\mathbf D})_{\beta}$ follows from exactness of  (\ref{eq exact sequence case 2}) in the case of $\mathbb{Z}$-graded manifolds of degree $2$. Hence in this case it is enough to assume that the vector field $\D$ is linear and non-degenerate. (Note that the vector field $\D$ is homological due to the weight agreement.)

\subsection{Two additional functors}

Let us fix a weight system $\Delta$ satisfying Definition \ref{de weight system} and the weight system $\Delta'=\Delta'(\Delta)$ as in Proposition \ref{prop obtain Delta_D from Delta_N}. Recall that we denoted by $r$ and $r'$ the ranks of $\Delta$ and $\Delta'$, respectively.  There is a projection $\mathsf{G}: \Delta'\to \Delta$ that is
defined in the following way.
Let us take 
$$
\delta' = \sum_{k\in K} \alpha_{i_k} + \sum_{(s,t)\in S\times T} \beta_{j_si_t}\in \Delta' 
$$
for certain sets $K$, $S$ and $T$. 
We set 
$$
\mathsf{G}(\delta):= \sum\limits_{k\in K} \alpha_{i_k} + \sum\limits_{(s,t)\in S\times T} \alpha_{i_t}.
$$ 
In other words, we replace any $\beta_{j_si_t}$ by $\alpha_{i_t}$. 

\medskip

\prop\label{[prop properties of weight system and projections]}
{\sl Let us take $\delta \in \Delta$.  We have
	$$
	\mathsf{G}^{-1} (\delta ) = \Delta'_{\delta},
	$$
	where $\Delta'_{\delta}$ is given by (\ref{eq Delta_delta new})
}

\medskip

\noindent{\it Proof}\, follows from definitions.$\Box$

\medskip

Denote by $\Delta'_{<\beta_{ji}}$ and $\Delta'_{=\beta_{ji}}$
the weight subsystems in $\Delta'$ generated by 
the sets 
\begin{align*}
\mathbf{A}_{<\beta_{ji}} = \{ \alpha_s,\, \beta_{st}\,\,|\,\, s=1,\ldots,r,\, \, \text{$t<i$ or $t=i$ and $s<j$}   \};\\
\mathbf{A}_{=\beta_{ji}}=\{ \alpha_s,\, \beta_{st}\,\,|\,\, s=1,\ldots,r,\, \, \text{$t<i$ or $t=i$ and $s\leq j$}   \},
\end{align*}
respectively. We put $\Delta_{<\beta_{ji}} := \mathsf{G}(\Delta'_{<\beta_{ji}})$ and  $\Delta_{=\beta_{ji}} := \mathsf{G}(\Delta'_{=\beta_{ji}})$.

In Section $5.2$ we constructed the functor $\mathbb F$ from the category of graded manifolds of type $\Delta$ to  the category of $r'$-fold vector bundles of type $\Delta'$. Now we need to construct in a similar way two additional functors  $\mathbb F_{<\beta_{ji}}$  and $\mathbb F_{=\beta_{ji}}$. The functor $\mathbb F_{<\beta_{ji}}$ is a functor from the category of graded manifolds of type $\Delta_{<\beta_{ji}}$ to  the category of $r'$-fold vector bundles of type $\Delta'_{<\beta_{ji}}$ and the functor $\mathbb F_{=\beta_{ji}}$ is a functor from the category of graded manifolds of type $\Delta_{=\beta_{ji}}$ to  the category of $r'$-fold vector bundles of type $\Delta'_{=\beta_{ji}}$, respectively. Note that always we deal with $r'$-fold vector bundles. 
Recall that to construct the functor $\mathbb F$ we used the additional formal weights $(\beta_{ji})$, where $j=2,\ldots, n_i$ and $i=i,\ldots, r$, see (\ref{eq additional weights beta def}). Similarly we define the functor $\mathbb F_{<\beta_{ji}}$  using the additional weights $\beta_{st}\in \mathbf{A}_{<\beta_{ji}}$ and the functor $\mathbb F_{=\beta_{ji}}$ using the additional weights $\beta_{st}\in \mathbf{A}_{=\beta_{ji}}$.

More precisely, let us describe for example the functor $\mathbb F_{<\beta_{ji}}$ in more details. We set $n':= | \mathbf{A}_{<\beta_{ji}}|$, i.e. $n'$ is the number of elements in $\mathbf{A}_{<\beta_{ji}}$. Further, we take a graded manifold $\mathcal N_{<\beta_{ji}}$ of type $\Delta_{<\beta_{ji}}$. Then we define $\mathbb F_{<\beta_{ji}}(\mathcal N_{<\beta_{ji}})$ in the following way. We take $n'$-iterated tangent bundle $\mathcal N'_{<\beta_{ji}}$ of $\mathcal N_{<\beta_{ji}}$ using sequentially additional weights from $\mathbf{A}_{<\beta_{ji}}$. Further, we consider the graded manifold with the structure sheaf $\mathcal O_{\mathcal N'_{<\beta_{ji}}}/\mathcal J^-$. Assume that it has type $\tilde{\Delta}_{<\beta_{ji}}$. We choose the maximal multiplicity free subset $\mathbb F_{<\beta_{ji}}( \Delta_{<\beta_{ji}})$ in $\tilde{\Delta}_{<\beta_{ji}}$ and denote by $\mathbb F_{<\beta_{ji}}(\mathcal N_{<\beta_{ji}})$ the corresponding to $\Delta'_{<\beta_{ji}}$ graded manifold.  The definition of the functor $\mathbb F_{=\beta_{ji}}$ is similar: we should replace the set $\mathbf{A}_{<\beta_{ji}}$ by $\mathbf{A}_{=\beta_{ji}}$.

Note that the constructed functors $\mathbb F_{<\beta_{ji}}$ and $\mathbb F_{=\beta_{ji}}$ are from the categories of graded manifolds of type $\Delta_{<\beta_{ji}}$ and of type $\Delta_{=\beta_{ji}}$ to the categories of graded manifolds of type $\mathbb F_{<\beta_{ji}}( \Delta_{<\beta_{ji}})$ and of type $\mathbb F_{=\beta_{ji}}( \Delta_{=\beta_{ji}})$, respectively. The weight system $\mathbb F_{<\beta_{ji}}( \Delta_{<\beta_{ji}})$ is defined by formulas (\ref{eq Delta_D_N expression new}) and (\ref{eq Delta_delta new}), where $\beta_{st}\in \mathbf{A}_{<\beta_{ji}}$. In the same way we define $\mathbb F_{=\beta_{ji}}( \Delta_{=\beta_{ji}})$. However in fact we have the following equalities.

\medskip

\prop\label{prop properties of weight system <beta}
{\sl We have
	\begin{align*}
		\mathbb F_{<\beta_{ji}}( \Delta_{<\beta_{ji}}) = \Delta'_{<\beta_{ji}},\quad \mathbb F_{=\beta_{ji}}( \Delta_{=\beta_{ji}}) = \Delta'_{=\beta_{ji}}.
	\end{align*}
	
}

\medskip

\noindent{\it Proof.} It is enough to prove only the second statement. Let us take 
$$
\delta' = \sum_{k} \alpha_{i_k} + \sum_{st} \beta_{j_si_t} \in \mathbb F_{=\beta_{ji}}( \Delta_{=\beta_{ji}}).
$$
Then $\alpha_{i_k},\beta_{j_si_t}\in \Delta'_{=\beta_{ji}}$. Hence, $\mathbb F_{=\beta_{ji}}( \Delta_{=\beta_{ji}}) \subset \Delta'_{=\beta_{ji}}$.
On other hand assume that 
$$
\delta' = \sum_{k} \alpha_{i_k} + \sum_{st} \beta_{j_si_t} \in \Delta'_{=\beta_{ji}}.
$$
Then $\delta'\in \Delta'_{\mathsf{G}(\delta')}$. Since $\delta'$ depends only on $\beta_{st}\in \mathbf{A}_{=\beta_{ji}}$, we see that $\delta'\in \mathbb F_{=\beta_{ji}}( \Delta_{=\beta_{ji}})$.$\Box$

\medskip

\subsection{$\mathbb F$ is an equivalence of categories}

In the previous section we constructed a functor $\mathbb F$ from $\Delta$\textsf{Man} to  $\Delta'$\textsf{VBVect}. 
Let us prove that  $\mathbb F$ is {\bf essentially surjective}. This is the most difficult part of our paper.

 Consider the set of additional weights $(\beta_{ji})$, see (\ref{eq additional weights beta def}), with the lexicographical order:  $\beta_{ji}< \beta_{j'i'}$ if $i<i'$ or $i=i'$ and $j<j'$. We will prove the essential surjectivity of $\mathbb F$ by induction on this order. 
 
  Let $\mathbf D$ be an object in $\Delta'$\textsf{VBVect}. Consider the weight systems $\Delta'_{<\beta_{ji}}$, $ \Delta_{<\beta_{ji}} := \mathsf{G}(\Delta'_{<\beta_{ji}})$, $\Delta'_{=\beta_{ji}}$, $\Delta_{=\beta_{ji}} := \mathsf{G}(\Delta'_{=\beta_{ji}})$ and the functors $\mathbb F_{<\beta_{ji}}$, $\mathbb F_{=\beta_{ji}}$ constructed in the previous section.  Clearly the weight system $\Delta'_{<\beta_{ji}}$ satisfies the conditions of Lemma \ref{lem Delta' subset Delta denote a smf}. We denote by $\mathbf D_{<\beta_{ji}}$ the  graded manifold of type $\Delta'_{<\beta_{ji}}\subset \Delta'$, see Construction 1, Section $4.1$.  Further, assume by induction that there exists a graded manifold $\mathcal N_{<\beta_{ji}}$ of type $\Delta_{<\beta_{ji}}$ such that 
 $$
 \mathbb F_{<\beta_{ji}}(\mathcal N_{<\beta_{ji}}) \simeq \mathbf D_{<\beta_{ji}},
 $$ 
and this isomorphism that we denote by $\varphi'$ commutes with all vector fields $\D_{\beta_{st}}$. Now our goal is to show that there exists a graded manifold $\mathcal N_{=\beta_{ji}}$ such that
$$
\mathbb F_{=\beta_{ji}}(\mathcal N_{=\beta_{ji}}) \simeq \mathbf D_{=\beta_{ji}},
$$ 
 where $\mathbf D_{=\beta_{ji}}$ is the  graded manifold of type $\Delta'_{=\beta_{ji}}\subset \Delta'$, see Construction 1, Section $4.1$.

Note that any graded manifold of type $\Delta_{<\beta_{ji}}$ is also a graded manifold of type $\Delta_{=\beta_{ji}}$. Hence we can apply functor $\mathbb F_{=\beta_{ji}}$ to $\mathcal N_{<\beta_{ji}}$ and get the graded manifold  
 $\mathbb F_{=\beta_{ji}}(\mathcal N_{<\beta_{ji}})$ of type $\Delta'_{=\beta_{ji}}$. However,  $\mathbb F_{=\beta_{ji}}(\mathcal N_{<\beta_{ji}})$ is also a graded manifold of type $\mathbb F_{=\beta_{ji}}(\Delta_{<\beta_{ji}})$, where 
 \begin{equation}\label{eq decomposition of Delta}
 \mathbb F_{=\beta_{ji}}(\Delta_{<\beta_{ji}}) = \Delta'_{<\beta_{ji}}\cup \Delta'_1\cup \Delta'_2 \subset \Delta'.
\end{equation}
 Here $\Delta'_1$ and $\Delta'_2$ are subsets in $\Delta'$ that are defined by the following formulas:
 \begin{equation}\label{eq def of Delta_1 and Delta_2}
 \begin{split}
 &\Delta'_1 = \{\D_{\beta_{ji}} (\delta)\,\,\, |\,\,\, \delta\in \Delta'_{<\beta_{ji}} \, :\, \D_{\beta_{ji}} (\delta)\in \Delta' \}; \\
& \Delta'_2 = \{\delta\in \Delta' \,\,\, |\,\,\,  \exists\beta_{j_0i} <\beta_{ji}\, :\, \D_{\beta_{j_0i}} (\delta)\in \Delta'_1\}.
 \end{split}
 \end{equation}
Explicitly, $\delta_1\in \Delta'_1$ and $\delta_2\in \Delta'_2$ have the following form:
\begin{equation}\label{eq def of Delta_i, explicite}
\delta_1= \beta_{ji} + \theta_1 \quad \text{and} \quad \delta_2= \alpha_i+  \beta_{ji} + \theta_2,
\end{equation}
where $\theta_i$ are independent on $\beta_{ji}$ and $\alpha_i$. In addition we assume that for $\delta_2$ there exists $\beta_{j_0i}<\beta_{ji}$ such that $\D_{\beta_{j_0i}} (\delta)\in \Delta'_1$. 
Again the weight system $\mathbb F_{=\beta_{ji}}(\Delta_{<\beta_{ji}})\subset \Delta'$ satisfies the conditions of Lemma \ref{lem Delta' subset Delta denote a smf}.
Denote by $\mathbf D_{\mathbb F_{=\beta_{ji}}(\Delta_{<\beta_{ji}})}$ the graded manifold of type $\mathbb F_{=\beta_{ji}}(\Delta_{<\beta_{ji}})$, see Construction 1, Section $4.1$.

 For simplicity of notations we denote the structure sheaves of graded manifolds $\mathbb F_{=\beta_{ji}}(\mathcal N_{<\beta_{ji}})$ and $\mathbf D_{\mathbb F_{=\beta_{ji}}(\Delta_{<\beta_{ji}})}$ by $\mathcal O$ and $\mathcal O'$, respectively. We also denote operators $\D_{\beta_{st}}$ in $\mathcal O$ and $\mathcal O'$ by the same letter.

\medskip

\prop\label{prop raznesenie} {\sl
Let $\mathcal N_{<\beta_{ji}}$ and $\mathbf D_{\mathbb F_{=\beta_{ji}}(\Delta_{<\beta_{ji}})}$ be the graded manifolds of type $\Delta_{<\beta_{ji}}$ and $\mathbb F_{=\beta_{ji}}(\Delta_{<\beta_{ji}})$ as above.
Then there exists an isomorphism
$$
\varphi:\mathbb F_{=\beta_{ji}}(\mathcal N_{<\beta_{ji}}) \to \mathbf D_{\mathbb F_{=\beta_{ji}}(\Delta_{<\beta_{ji}})}
$$	
of graded manifolds of type $\mathbb F_{=\beta_{ji}}(\Delta_{<\beta_{ji}})$.
} 
 
\medskip

\noindent{\it Proof.} The idea of the proof is to extend the isomorphism $\varphi'$ using the vector fields $\D_{\beta_{st}}$. By Proposition \ref{prop isomorphism of graded mnf} our goal is to construct compatible bundle maps $\varphi_{\delta}$ for all $\delta \in \mathbb F_{=\beta_{ji}}(\Delta_{<\beta_{ji}})$. We use the decomposition (\ref{eq decomposition of Delta}). By induction we have an isomorphism $\varphi': \mathbb F_{<\beta_{ji}}(\mathcal N_{<\beta_{ji}}) \to\mathbf D_{<\beta_{ji}}$. For any $\delta \in \Delta'_{<\beta_{ji}}$ we put $\varphi_{\delta}:= \varphi'_{\delta}$. 
Further, let us take $\D_{\beta_{ji}} (\delta) \in \Delta'_1$. We put
\begin{equation}\label{eq def of phi, group A}
\varphi_{\D_{\beta_{ji}}(\delta)}:= \D_{\beta_{ji}} \circ\varphi_{\delta}\circ \D^{-1}_{\beta_{ji}}.
\end{equation}
Note that (\ref{eq def of phi, group A}) is well-defined since $ \delta\in \Delta'_{<\beta_{ji}}$. If $\delta\in \Delta'_2$ we put
\begin{equation}\label{eq def of phi, group B}
\varphi_{\delta}:= \D^{-1}_{\beta_{j_0i}} \circ\varphi_{\D_{\beta_{j_0i}}(\delta)}\circ \D_{\beta_{j_0i}},
\end{equation}
where $\beta_{j_0i}$ is as in the definition of $\Delta'_2$. 
Note that $\varphi_{\D_{\beta_{j_0i}}(\delta)}$ is defined by (\ref{eq def of phi, group A}), since $\D_{\beta_{j_0i}}(\delta)\in  \Delta'_1$. 
Combining (\ref{eq def of phi, group A}) and (\ref{eq def of phi, group B}), we get
$$
\varphi_{\delta}:= \D^{-1}_{\beta_{j_0i}} \circ\D_{\beta_{ji}} \circ\varphi_{\lambda}\circ \D^{-1}_{\beta_{ji}}\circ \D_{\beta_{j_0i}}.
$$
Explicitly, if $\delta = \alpha_i+  \beta_{ji} + \theta_2\in \Delta'_2$, see (\ref{eq def of Delta_i, explicite}), then $\lambda= \alpha_i+  \beta_{j_0i} + \theta_2\in \Delta'_{<\beta_{ji}}$.

The formula (\ref{eq def of phi, group A}) is well-defined, while  (\ref{eq def of phi, group B}) dependents on the choice of $\beta_{j_0i}$. Let us show that in fact this is not the case. Assume that 
\begin{align*}
&\varphi_{\delta}= \D^{-1}_{\beta_{j_0i}} \circ\D_{\beta_{ji}} \circ\varphi_{\lambda}\circ \D^{-1}_{\beta_{ji}}\circ \D_{\beta_{j_0i}},
&\varphi_{\delta}= \D^{-1}_{\beta_{j_1i}} \circ\D_{\beta_{ji}} \circ\varphi_{\lambda'}\circ \D^{-1}_{\beta_{ji}}\circ \D_{\beta_{j_1i}}.
\end{align*}
We need to show that
\begin{equation}\label{eq phi_delta is well defined}
\varphi_{\lambda}= \D^{-1}_{\beta_{ji}}\circ \D_{\beta_{j_0i}} \circ \D^{-1}_{\beta_{j_1i}} \circ\D_{\beta_{ji}} \circ\varphi_{\lambda'}\circ \D^{-1}_{\beta_{ji}}\circ \D_{\beta_{j_1i}}\circ \D_{\beta_{j_0i}}^{-1} \circ \D_{\beta_{ji}},
\end{equation}
where $\lambda= \alpha_i+  \beta_{j_0i} + \theta_2$ and $\lambda'= \alpha_i+  \beta_{j_1i} + \theta_2$ such that $\lambda,\lambda'\in \Delta'_{<\beta_{ji}}$. 
First of all consider $\varphi_{\lambda}|_{\Ker \D_{\beta_{j_0i}}}$. Then we can apply the cocycle condition for our vector fields:
$$
(\D^{-1}_{\beta_{ji}}\circ \D_{\beta_{j_1i}}\circ \D_{\beta_{j_0i}}^{-1} \circ \D_{\beta_{ji}})|_{\Ker \D_{\beta_{j_0i}}} = - \D_{\beta_{j_0i}}^{-1} \circ \D_{\beta_{j_1i}}|_{\Ker \D_{\beta_{j_0i}}}.
$$
Further, we use the relation:
$$
(\D_{\beta_{j_0i}}^{-1} \circ \D_{\beta_{j_1i}}) (\Ker \D_{\beta_{j_0i}}) = \Ker \D_{\beta_{j_1i}}.
$$
 Therefore we can apply the cocycle condition again
 $$
 (\D^{-1}_{\beta_{ji}}\circ \D_{\beta_{j_0i}} \circ \D^{-1}_{\beta_{j_1i}} \circ\D_{\beta_{ji}})|_{\Ker \D_{\beta_{j_1i}}} = - \D_{\beta_{j_1i}}^{-1} \circ \D_{\beta_{j_0i}}|_{\Ker \D_{\beta_{j_1i}}}.
 $$
 Now we can rewrite (\ref{eq phi_delta is well defined}) in the following form
 \begin{align*}
 \varphi_{\lambda}|_{\Ker \D_{\beta_{j_0i}}}= (\D_{\beta_{j_1i}}^{-1} \circ \D_{\beta_{j_0i}} \circ\varphi_{\lambda'}\circ \D_{\beta_{j_0i}}^{-1} \circ \D_{\beta_{j_1i}})|_{\Ker \D_{\beta_{j_0i}}}
 \end{align*}
 or
 \begin{align*}
(\D_{\beta_{j_0i}}^{-1} \circ \D_{\beta_{j_1i}} \circ \varphi_{\lambda})|_{\Ker \D_{\beta_{j_0i}}}= (\varphi_{\lambda'}\circ \D_{\beta_{j_0i}}^{-1} \circ \D_{\beta_{j_1i}})|_{\Ker \D_{\beta_{j_0i}}}.
 \end{align*}
 This equation holds because $\varphi_{\lambda}$ and $\varphi_{\lambda'}$ commute with vector fields $\D_{\beta_{j_0i}}$ and $\D_{\beta_{j_1i}}$ by induction.

 To show (\ref{eq phi_delta is well defined}) our next step is to use the decomposition 
 $$
 \mathcal O_{\lambda} =  \mathcal O_{\lambda}\cap \Ker \D_{\beta_{j_0i}} + \bigoplus_{\lambda_1+\lambda_2=\lambda} \mathcal O_{\lambda_1} \mathcal O_{\lambda_2},
 $$
 where $\lambda_i\ne  0$. (Note that the existence of such kind of decompositions follows from the definition of the category $\Delta'$\textsf{VBVect}.) Consider the following two cases:

 \medskip
 \noindent{\bf 1.} Assume that $\lambda_1= \alpha_i+ \gamma_1$ and $\lambda_2= \beta_{j_0i}+ \gamma_2$ such that $\lambda_1+\lambda_2=\lambda$, and $f_i$ is a function of weight $\lambda_i$, $i=1,2$.  
We have
\begin{align*}
(\D^{-1}_{\beta_{ji}}\circ \D_{\beta_{j_1i}}\circ \D_{\beta_{j_0i}}^{-1} \circ \D_{\beta_{ji}})(f_1f_2) = (\D^{-1}_{\beta_{ji}}\circ \D_{\beta_{j_1i}}\circ \D_{\beta_{j_0i}}^{-1})( \D_{\beta_{ji}}(f_1) \cdot f_2) =\\
(-1)^{(f_1+1)} \D^{-1}_{\beta_{ji}}\circ \D_{\beta_{j_1i}} ( \D_{\beta_{ji}}(f_1) \cdot \D_{\beta_{j_0i}}^{-1}(f_2))
 = \\
 \D^{-1}_{\beta_{ji}} ( \D_{\beta_{ji}}(f_1) \cdot \D_{\beta_{j_1i}} \circ\D_{\beta_{j_0i}}^{-1}(f_2))=
 ( f_1 \cdot \D_{\beta_{j_1i}} \circ\D_{\beta_{j_0i}}^{-1}(f_2)).
\end{align*} 
Further, since $\varphi_{\lambda'} = \varphi'_{\lambda'}$, we have
\begin{align*}
\varphi_{\lambda'}  ( f_1 \cdot \D_{\beta_{j_1i}} \circ\D_{\beta_{j_0i}}^{-1}(f_2)) =  \varphi_{\lambda'_1}(f_1) \cdot \varphi_{\lambda'_2}(\D_{\beta_{j_1i}} \circ\D_{\beta_{j_0i}}^{-1}(f_2)), 
\end{align*}
 where $\lambda'_1 = \lambda_1= \alpha_i+ \gamma_1$ and $\lambda'_1 = \beta_{j_1i}+ \gamma_2$.
 Similarly we get
 \begin{align*}
( \D^{-1}_{\beta_{ji}}\circ \D_{\beta_{j_0i}} \circ \D^{-1}_{\beta_{j_1i}} \circ\D_{\beta_{ji}}) (\varphi_{\lambda'_1}(f_1) \cdot \varphi_{\lambda'_2}(\D_{\beta_{j_1i}} \circ\D_{\beta_{j_0i}}^{-1}(f_2))) 
 =\\
 \varphi_{\lambda'_1}(f_1) \cdot (\D_{\beta_{j_1i}} \circ\D_{\beta_{j_0i}}^{-1}) \circ \varphi_{\lambda'_2}(\D_{\beta_{j_1i}} \circ\D_{\beta_{j_0i}}^{-1}(f_2))=\\
  \varphi_{\lambda_1}(f_1) \cdot  \varphi_{\lambda_2}(f_2) = \varphi_{\lambda} (f_1 \cdot f_2 ).
  \end{align*}
 
  \medskip
  
  \noindent{\bf 2.} Assume that $\lambda_1= \alpha_i+ \beta_{j_0i}+ \gamma_1$ and $\lambda_2= \gamma_2$, where $\gamma_i$ are not depending on $\alpha_i$ and $\beta_{j_0i}$. Consider the restriction of  (\ref{eq phi_delta is well defined}) on $\mathcal O_{\lambda_1} \cdot \mathcal O_{\lambda_2}$, We get
  \begin{align*}
  &\D^{-1}_{\beta_{ji}}\circ \D_{\beta_{j_0i}} \circ \D^{-1}_{\beta_{j_1i}} \circ\D_{\beta_{ji}} \circ\varphi_{\lambda'}\circ \D^{-1}_{\beta_{ji}}\circ \D_{\beta_{j_1i}}\circ \D_{\beta_{j_0i}}^{-1} \circ \D_{\beta_{ji}}|_{\mathcal O_{\lambda_1} \cdot \mathcal O_{\lambda_2}}
  =\\
   &(\D^{-1}_{\beta_{ji}}\circ \D_{\beta_{j_0i}} \circ \D^{-1}_{\beta_{j_1i}} \circ\D_{\beta_{ji}} \circ\varphi_{\lambda'_1}\circ \D^{-1}_{\beta_{ji}}\circ \D_{\beta_{j_1i}}\circ \D_{\beta_{j_0i}}^{-1} \circ \D_{\beta_{ji}})|_{\mathcal O_{\lambda_1}} \cdot \varphi_{\lambda'_2}|_{ \mathcal O_{\lambda_2}},
  \end{align*}
  where $\lambda'_1 = \alpha_i+ \beta_{j_1i}+ \gamma_1$ and $\lambda'_2=\lambda_2$. Hence to prove (\ref{eq phi_delta is well defined}) we need to prove
  $$
 \varphi_{\lambda_1}|_{\mathcal O_{\lambda_1}} = (\D^{-1}_{\beta_{ji}}\circ \D_{\beta_{j_0i}} \circ \D^{-1}_{\beta_{j_1i}} \circ\D_{\beta_{ji}} \circ\varphi_{\lambda'_1}\circ \D^{-1}_{\beta_{ji}}\circ \D_{\beta_{j_1i}}\circ \D_{\beta_{j_0i}}^{-1} \circ \D_{\beta_{ji}})|_{\mathcal O_{\lambda_1}} 
  $$
  Since $\lambda_2\ne 0$, the last equality follows by  induction on length of $\gamma_1$, where the length of $\gamma_1$ is equal to the number of summands in $\gamma_1$. Note that in case $\lambda_1= \alpha_i+ \beta_{j_0i}$ (in other words the length of $\gamma_1=0$, the basis of our induction)  
  the result follows from the decomposition 
  $$
   \mathcal O_{\alpha_i+ \beta_{j_0i}} =  (\mathcal O_{\alpha_i+ \beta_{j_0i}})\cap \Ker \D_{\beta_{j_0i}} +  \mathcal O_{\alpha_i} \mathcal O_{\beta_{j_0i}}
  $$
  and the case $1$. Hence we proved that $\varphi_{\delta}$ is well-defined.  
  
It remains to prove the compatibility condition of Proposition \ref{prop isomorphism of graded mnf}
$$
\varphi_{\delta}|_{\mathcal O_{\delta_1} \cdot \mathcal O_{\delta_2}} = (\varphi_{\delta_1}\cdot \varphi_{\delta_2})|_{\mathcal O_{\delta_1} \cdot \mathcal O_{\delta_2}} 
$$
for $\delta = \delta_1+\delta_2$, where $\delta_i\ne 0$.
Consider first the case $\delta\in \Delta'_1$. Without loss of generality we may assume that $\delta_1$ depends on $\beta_{ji}$.
For $f_i$ of weight $\delta_i$, we have
\begin{align*}
\varphi_{\delta}(f_1\cdot f_2) = (\D_{\beta_{ji}} \circ\varphi'_{\D^{-1}_{\beta_{ji}}(\delta)}\circ \D^{-1}_{\beta_{ji}}) (f_1\cdot f_2) = 
(\D_{\beta_{ji}} \circ\varphi'_{\D^{-1}_{\beta_{ji}}(\delta)}) (\D^{-1}_{\beta_{ji}}(f_1)\cdot f_2)
 =
  \\
\D_{\beta_{ji}} (\varphi'_{\D^{-1}_{\beta_{ji}}(\delta_1)}(\D^{-1}_{\beta_{ji}}(f_1))\cdot \varphi'_{\delta_2}(f_2)) = (\D_{\beta_{ji}} \circ\varphi'_{\D^{-1}_{\beta_{ji}}(\delta_1)}\circ \D^{-1}_{\beta_{ji}})(f_1)\cdot \varphi'_{\delta_2} (f_2)=\\
\varphi_{\delta_1}(f_1)\cdot \varphi_{\delta_2}(f_2).
\end{align*}
We used here the fact that $\varphi'$ is an isomorphism by induction.

Further, assume that $\delta=\delta_1+\delta_2\in \Delta'_2$. Again without loss of generality we may assume that $\delta_1$ depends on $\alpha_{i}$. Similarly we have for $f_i$ of weight $\delta_i$
\begin{align*}
\varphi_{\delta}(f_1\cdot f_2) = &(\D^{-1}_{\beta_{j_0i}} \circ\varphi_{\D_{\beta_{j_0i}}(\delta)}\circ \D_{\beta_{j_0i}}) (f_1\cdot f_2)
=
\\
& (\D^{-1}_{\beta_{j_0i}} \circ\varphi_{\D_{\beta_{j_0i}}(\delta_1)}\circ \D_{\beta_{j_0i}})(f_1)\cdot \varphi_{\delta_2} (f_2) = \varphi_{\delta_1}(f_1)\cdot \varphi_{\delta_2}(f_2).
\end{align*}
By Proposition \ref{prop isomorphism of graded mnf}, an isomorphism $\varphi$ is defined by the collection $(\varphi_{\delta})$, where  $\delta \in \mathbb F_{=\beta_{ji}}(\Delta_{<\beta_{ji}})$. The proof is complete.$\Box$ 
    \medskip

  Note that in the proof of Proposition  \ref{prop raznesenie} we used the folowing decomposition  
  	$$
  \mathcal O_{\delta}= \mathcal O_{\delta}\bigcap_{b_{st} \ne 0}  \Ker \D_{\beta_{st}}+	\bigoplus_{\delta_1+\delta_2= \delta} \mathcal O_{\delta_1} \mathcal O_{\delta_2},
  $$
  where $\delta_i\ne 0$. The idea was the following. First we show a certain equality on kernels $\mathcal O_{\delta}\bigcap\limits_{b_{st} \ne 0}  \Ker \D_{\beta_{st}}$ and then  using induction on $\bigoplus\limits_{\delta_1+\delta_2= \delta} \mathcal O_{\delta_1} \mathcal O_{\delta_2}$. Further we will use this idea several time. We will call this argument the {\bf decomposition and induction argument}. 
  
  Our goal now is to prove that the constructed isomorphism $\varphi$ is a morphism in the category $\Delta'$\textsf{VBVect}.  
  
  \medskip
 
  \prop\label{prop phi icommute with operators} {\sl
  	The isomorphism $	\varphi$ from Proposition \ref{prop raznesenie}
  commutes with operators $\D_{\beta_{st}}$. 
  
  } 
  
  \medskip
  
\noindent{\it Proof.} We need to show that
\begin{equation}\label{eq phi commuts with operators}
\D_{\beta_{st}}\circ \varphi_{\delta} = \varphi_{\D_{\beta_{st}}(\delta)} \circ \D_{\beta_{st}}
\end{equation}
for any operator $\D_{\beta_{st}}$ and for any $\delta$ from our weight lattice. Since $\D_{\beta_{st}}$ is a vector field, it is enough to show (\ref{eq phi commuts with operators}) for $\delta \in \mathbb F_{=\beta_{ji}}(\Delta_{<\beta_{ji}})$.  We use the decomposition (\ref{eq decomposition of Delta}). 
 Consider the following cases.

\medskip

\noindent{\bf 1.} {\bf Case $\delta \in \Delta'_{<\beta_{ji}}$ and $(st)\ne (ji)$.} In this case  (\ref{eq phi commuts with operators}) holds by the assumption that $\varphi'$ commutes with all operators.

\medskip
\noindent{\bf 2.} {\bf Case $\delta \in \Delta'_{<\beta_{ji}}$
  and $(st)= (ji)$.} If $\delta$ is independent on $\alpha_i$, the equality (\ref{eq phi commuts with operators}) holds trivially. If $\delta = \alpha_i + \theta$, where $\theta$ is independent on $\beta_{ji}$, then $\D_{\beta_{ji}}(\delta)\in \Delta'_1$ and (\ref{eq phi commuts with operators}) holds by definition of $\varphi_{\D_{\beta_{ji}}(\delta)}$.

\medskip
\noindent{\bf 3.} {\bf Case $\delta \in \Delta'_1$ and $(st)\ne (ji)$.} 
If $t=i$, the equality (\ref{eq phi commuts with operators}) holds trivially, since $\delta$ does not depend on $\alpha_i$. Assume that $t\ne i$. Consider the decomposition of the weight $\D_{\beta_{st}}(\delta) = \gamma_1+ \gamma_2$, where $\gamma_1=\beta_{ji}+\cdots\in \Delta'_1$ and $\gamma_2=\D_{\beta_{st}}(\delta) - \gamma_1$. Note that $\gamma_2$ is independent on $\beta_{ji}$. Then the following bundle isomorphism is defined
$$
\D_{\beta_{ji}}^{-1} : \mathcal O_{\gamma_1+\gamma_2} \to \mathcal O_{\D_{\beta_{ji}}^{-1}(\gamma_1)+\gamma_2}, \quad \D_{\beta_{ji}}^{-1}(f_1\cdot f_2):= \D_{\beta_{ji}}^{-1}(f_1)\cdot f_2,
$$
where $f_i\in \mathcal O_{\gamma_i}$. Therefore, we have
$$
\D_{\beta_{st}} \circ \D_{\beta_{ji}}^{-1}|_{ \mathcal O_{\delta}} = \D^{-1}_{\beta_{ji}} \circ \D_{\beta_{st}}|_{ \mathcal O_{\delta}}.
$$
Further, 
\begin{align*}
\varphi_{\D_{\beta_{st}}(\delta)} \circ \D_{\beta_{st}}= (\varphi_{\gamma_1}\cdot \varphi_{\gamma_2}) \circ \D_{\beta_{st}}= ((\D_{\beta_{ji}} \circ\varphi_{\D^{-1}_{\beta_{ji}}(\gamma_1)}\circ \D^{-1}_{\beta_{ji}})\cdot \varphi_{\gamma_2}) \circ \D_{\beta_{st}}.
\end{align*}
On the other hand we have
\begin{align*}
\D_{\beta_{st}}\circ \varphi_{\delta} = \D_{\beta_{st}}\circ (\D_{\beta_{ji}} \circ\varphi_{\D^{-1}_{\beta_{ji}}(\delta)}\circ \D^{-1}_{\beta_{ji}}) = (\D_{\beta_{ji}} \circ\varphi_{(\D_{\beta_{st}} \circ\D^{-1}_{\beta_{ji}})(\delta)}\circ \D^{-1}_{\beta_{ji}})\circ \D_{\beta_{st}}=\\
\big((\D_{\beta_{ji}} \circ\varphi_{\D^{-1}_{\beta_{ji}}(\gamma_1)}\circ \D^{-1}_{\beta_{ji}})\cdot \varphi_{\gamma_2}\big)\circ \D_{\beta_{st}}.
\end{align*}
 Hence,  (\ref{eq phi commuts with operators}) is proven for this case. 

\medskip
\noindent{\bf 4.} {\bf Case $\delta \in \Delta'_2$ and $(st)\ne (ji)$.} Assume that $t=i$ and $\D_{\beta_{si}}(\delta)$ is a weight. In this case $\D_{\beta_{si}}(\delta)\in \Delta'_1$. We have
\begin{align*}
\D_{\beta_{st}}\circ \varphi_{\delta} =  \D_{\beta_{st}}\circ (\D^{-1}_{\beta_{si}} \circ\varphi_{\D_{\beta_{si}}(\delta)}\circ \D_{\beta_{si}}) = \varphi_{\D_{\beta_{si}}(\delta)}\circ \D_{\beta_{si}}.
\end{align*}
Hence, (\ref{eq phi commuts with operators}) holds.

Assume that $t=i$ and $\D_{\beta_{si}}(\delta)$ is not a weight. Then since $\D_{\beta_{si}}(\delta)\notin \mathbb F_{=\beta_{ji}}(\Delta_{<\beta_{ji}})$, it follows that $\delta = \beta_{si} + \ldots$, i.e. $\delta$ depends on $\beta_{si}$ non-trivially.  Hence, $\D_{\beta_{si}}(\delta)=\gamma_1 + \gamma_2$, where $\gamma_1\in \Delta'_1$.  In this case (\ref{eq phi commuts with operators}) is equivalent to the following equality:
\begin{equation}\label{eq commutativity of operators, case 3}
\begin{split}
\D_{\beta_{si}}\circ (\D^{-1}_{\beta_{j_0i}} \circ \D_{\beta_{ji}} \circ\varphi_{\D^{-1}_{\beta_{ji}} \circ \D_{\beta_{j_0i}}(\delta)}\circ \D^{-1}_{\beta_{ji}} \circ \D_{\beta_{j_0i}}) & = \\
\big( (\D_{\beta_{ji}} \circ\varphi_{\D^{-1}_{\beta_{ji}}(\gamma_1)}&\circ \D^{-1}_{\beta_{ji}}) \circ \varphi_{\gamma_2} \big)\circ  \D_{\beta_{si}},
\end{split}
\end{equation}
where $j_0\ne s,j$. Now we use the decomposition and induction argument. By Proposition \ref{prop the sheaf O_delta is generated by ker and smaller}, we have the following decomposition:
	$$
\mathcal O_{\delta}= \mathcal O_{\delta}\bigcap_{b_{j_si_s}(\delta) \ne 0}  \Ker \D_{\beta_{j_si_s}}+	\bigoplus_{\delta_1+\delta_2= \delta} \mathcal O_{\delta_1} \mathcal O_{\delta_2},
$$
where $\delta_i\ne 0$ and $b_{j_si_s}(\delta)\in \mathbb C$ is the multiplicity of $\beta_{j_si_s}$ in $\delta$. Since $\delta$ depends on $\beta_{si}$ non-trivially and since the maps  $\D^{-1}_{\beta_{ji}} \circ \D_{\beta_{j_0i}}$ and
 $\D^{-1}_{\beta_{j_0i}} \circ \D_{\beta_{ji}}$ preserve this decomposition, (\ref{eq commutativity of operators, case 3}) follows from the previous cases and by induction.

Consider now the case $\delta \in \Delta'_2$ and $t\ne i$.
If $\delta$ does not depends on $\alpha_t$, then (\ref{eq phi commuts with operators}) holds trivially. Assume that $\delta = \alpha_t +\cdots$ and $\delta$ does not depend on $\beta_{st}$. Then $\D_{\beta_{st}}(\delta)$ is a weight.  In this case we have:
\begin{align*}
\D_{\beta_{st}} \circ (\D^{-1}_{\beta_{j_0i}} \circ \varphi_{\D_{\beta_{j_0i}} (\delta)} 
\circ \D_{\beta_{j_0i}}) =  \D^{-1}_{\beta_{j_0i}} \circ \D_{\beta_{st}} \circ \varphi_{\D_{\beta_{j_0i}} (\delta)} 
\circ \D_{\beta_{j_0i}} = \\
(\D^{-1}_{\beta_{j_0i}}  \circ \varphi_{\D_{\beta_{st}}\circ\D_{\beta_{j_0i}} (\delta)} \circ \D_{\beta_{j_0i}}) \circ \D_{\beta_{st}}.
\end{align*}
Hence  (\ref{eq phi commuts with operators}) holds. 
Further, assume that $\delta = \alpha_t + \beta_{st} + \cdots$. 
 In this case (\ref{eq phi commuts with operators}) is equivalent to the following equality:
\begin{align*}
\D_{\beta_{st}}\circ (\D^{-1}_{\beta_{j_0i}} \circ \D_{\beta_{ji}} \circ\varphi_{\D^{-1}_{\beta_{ji}} \circ \D_{\beta_{j_0i}}(\delta)}\circ \D^{-1}_{\beta_{ji}} \circ \D_{\beta_{j_0i}}) & = \\
\big((\D^{-1}_{\beta_{j_0i}} \circ \D_{\beta_{ji}} \circ\varphi_{\D^{-1}_{\beta_{ji}} \circ \D_{\beta_{j_0i}}\circ\D_{\beta_{st}}(\delta)}\circ &\D^{-1}_{\beta_{ji}} \circ \D_{\beta_{j_0i}})\big)\circ  \D_{\beta_{st}},
\end{align*}
where $j_0\ne j$. This holds by the decomposition and induction argument.

\medskip
\noindent{\bf 5.} {\bf Case $\delta\in \Delta'_1\cup \Delta'_2$ and  $(st)= (ji)$.} If $\delta \in \Delta'_1$, then (\ref{eq phi commuts with operators}) holds trivially. Further, assume that $\delta \in \Delta'_2$. In this case the result follows from the decomposition and induction argument. 
The proof is complete. $\Box$

\medskip

Recall that  we denoted by $\mathbf D_{=\beta_{ji}}$ the  graded manifold of type $\Delta'_{=\beta_{ji}}\subset \Delta'$, see Construction 1, Section $4.1$. Now we have the following situation. By induction we assumed that there exists an isomorphism $\varphi': \mathbb F_{<\beta_{ji}}(\mathcal N_{<\beta_{ji}}) \to\mathbf D_{<\beta_{ji}}$ of graded manifolds of type $\Delta'_{<\beta_{ji}}$ that commutes with all operators. By Propositions \ref{prop raznesenie} and \ref{prop phi icommute with operators} there exists an isomorphism $
\varphi:\mathbb F_{=\beta_{ji}}(\mathcal N_{<\beta_{ji}}) \to \mathbf D_{\mathbb F_{=\beta_{ji}}(\Delta_{<\beta_{ji}})}
$ of graded manifolds of type $\mathbb F_{=\beta_{ji}}(\Delta_{<\beta_{ji}})$ compatible with $\varphi'$ that also commutes with all operators. Our goal now is to prove the following proposition.

\medskip

\prop \label{prop prisoedinenie weights} {\sl
There exists a graded manifold $\mathcal N_{=\beta_{ji}}$ of type $\Delta_{=\beta_{ji}} $  such that
	$$
	\mathbb F_{=\beta_{ji}}(\mathcal N_{=\beta_{ji}}) \simeq \mathbf D_{=\beta_{ji}}
	$$	
and this isomorphism say $\psi$ commutes with all operators.
}

\medskip

\noindent{\it Proof.} First of all note that by Proposition \ref{prop properties of weight system <beta}, $\psi$ is an isomorphism of graded manifolds of type $\Delta'_{=\beta_{ji}}$. Further, clearly we have
 $$
 \Delta'_{=\beta_{ji}} = \mathbb F_{=\beta_{ji}}(\Delta_{<\beta_{ji}}) \cup S,
 $$
where $S$ is the subset in $\Delta'$ that contains all weights in the form $\alpha_i+ \sum\limits_{q=2}^j\beta_{qi} +\gamma$. We prove this proposition by induction on the length of $\gamma$. (Recall that the length $|\theta|$ of a multiplicity free weight $\theta$ is the number of summands in $\theta$.)

If the set $S$ is empty by Propositions \ref{prop raznesenie} and \ref{prop phi icommute with operators} we are done. Assume that $S\ne \emptyset$. Denote by $S_p$ the subset in $S$ such that $|\gamma|=p$. So we have $S= \cup_{p\geq 0} S_p$. Let us take $\delta'\in S_p$, where $p\geq 0$, satisfying the following property: if $\delta'$ depends on $\beta_{st}$, then $\delta'$ depends also on $\alpha_t$. (Note that by  definition of $\Delta'$ we always can find such $\delta'$.)
In other words, if we put $\delta' = \sum\limits_s a'_s\alpha_s + \sum\limits_{pq} b'_{pq}\beta_{pq}$, then from $b'_{st}=1$ it follows that $a'_t=1$. 
 
 \medskip
\noindent{\bf Step 1, construction of graded manifold.}
 Consider the following exact sequence
 $$
 0\to \bigoplus_{\delta'_1+\delta'_2= \delta'} (\mathcal O_{\mathbf D})_{\delta'_1} (\mathcal O_{\mathbf D})_{\delta'_2}
 \longrightarrow (\mathcal O_{\mathbf D})_{\delta'} \longrightarrow \mathcal E_{\delta'}\to 0,
 $$
 where $\delta'_i\ne 0$, $i=1,2$, and $\mathcal E_{\delta'}$ is a certain locally free sheaf. As in (\ref{eq sheaf S_delta}),  for any $\delta'\in \Delta'$, we put 
 $$
 \mathcal S_{\delta'}:= (\mathcal O_{\mathbf D})_{\delta'}\bigcap_{b'_{st} \ne 0}  \Ker \D_{\beta_{st}}.
 $$
  By definition of the category $\Delta'$\textsf{VBVect}, see Section $7.1$, we have the following decomposition
$$
(\mathcal O_{\mathbf D})_{\delta'}= \mathcal S_{\delta'} +	\bigoplus_{\delta'_1+\delta'_2= \delta'} (\mathcal O_{\mathbf D})_{\delta'_1} (\mathcal O_{\mathbf D})_{\delta'_2},
$$
where $\delta_i\ne 0$. Hence the following sequence is also exact
\begin{align*}
0\to \bigoplus_{\delta'_1+\delta'_2= \delta'} (\mathcal O_{\mathbf D})_{\delta'_1} (\mathcal O_{\mathbf D})_{\delta'_2}\bigcap_{b_{st} \ne 0}  \Ker \D_{\beta_{st}}
\stackrel{\iota}{\longrightarrow} 
 \mathcal S_{\delta'} \longrightarrow \mathcal E_{\delta'}\to 0.
\end{align*}
By induction we assume that  there is a graded manifold $\mathcal M$ of type $\Delta_{\mathcal M} = \mathsf{G} (\Delta'_{\mathcal M})$, where
$$
\Delta'_{\mathcal M}= \mathbb F_{=\beta_{ji}}(\Delta_{<\beta_{ji}}) \cup \bigcup_{q < p} S_q,
$$
and an isomorphism $\varphi_{\mathcal M}: \mathbb F_{=\beta_{ji}}(\mathcal M) \to \mathbf D_{\Delta'_{\mathcal M}}$ that is compatible with $\varphi$ and commutes with all operators. 
Here $\mathbf D_{\Delta'_{\mathcal M}}$ is again defined by  Construction 1, Section $4.1$.

We put $\delta:= \mathsf{G}(\delta')$. By Proposition \ref{prop D_delta is an iso onto kernel}, there is an isomorphism 
$$
\varphi_{\mathcal M} \circ \D^{\bar\Lambda} : (\mathcal O_{\mathcal M})_{\delta}\to   \bigoplus_{\delta'_1+\delta'_2= \delta'} (\mathcal O_{\mathbf D})_{\delta'_1} (\mathcal O_{\mathbf D})_{\delta'_2}\bigcap_{b'_{st} \ne 0}  \Ker \D_{\beta_{st}}.
$$
Here we assume that $ \D^{\bar\Lambda}(\delta) = \delta'$ and that $\bar\Lambda$ has lexicographical order. 

Note that $(\mathcal O_{\mathcal M})_{\delta} = \bigoplus\limits_{\delta_1+\delta_2= \delta} (\mathcal O_{\mathcal M})_{\delta_1} (\mathcal O_{\mathcal M})_{\delta_2}$, where $\delta_i\ne 0$, since $\delta\notin \Delta_{\mathcal M}$. 
Hence the following sequence is exact
\begin{align*}
0\to \bigoplus\limits_{\delta_1+\delta_2= \delta} (\mathcal O_{\mathcal M})_{\delta_1} (\mathcal O_{\mathcal M})_{\delta_2}
\stackrel{\iota\circ \varphi_{\mathcal M}\circ \D^{\bar\Lambda}}{\longrightarrow} 
\mathcal S_{\delta'} \longrightarrow \mathcal E_{\delta'}\to 0.
\end{align*}
Now we use Construction $2$, Section $4.2$ to build a graded manifold $\widetilde{\mathcal M}$ of type 
$\Delta_{\mathcal M}\cup \mathsf{G}(S_p)$. In more details, we put 
$$
(\mathcal O_{\widetilde{\mathcal M}})_{\delta}:= \mathcal S_{\delta'}.
$$
If there is  $\lambda'\in S_p$, where $\lambda\ne \delta'$, satisfying the following property: if $\lambda$ depends on $\beta_{st}$, then $\delta'$ depends also on $\alpha_t$, then we repeat this construction and define $(\mathcal O_{\widetilde{\mathcal M}})_{\lambda}$, where $\lambda= \mathsf{G} (\lambda')$. So we defined a graded manifold $\widetilde{\mathcal M}$ of type 
$\Delta_{\mathcal M}\cup \mathsf{G}(S_p)$.

 \medskip
\noindent{\bf Step 2, construction of an isomorphism.}
Our goal now is to construct an isomorphism $\widetilde{\psi} : \mathbb F_{=\beta_{ji}}(\widetilde{\mathcal M}) \to  \widetilde{\mathbf D}$ , where $\widetilde{\mathbf D}$ is a graded manifold of type $\Delta'_{\mathcal M}\cup S_p$ that is defined  by  Construction 1, Section $4.1$. We put $\widetilde{\psi}_{\theta} = (\varphi_{\mathcal M})_{\theta}$ for any $\theta \in \Delta'_{\mathcal M}$. To define $\widetilde{\psi}_{\delta'}$, where $\delta'\in S_p$ satisfy the following property: if $\delta'$ depends on $\beta_{st}$, then $\delta'$ depends also on $\alpha_t$,  we use the decomposition:
$$
\big(\mathcal O_{\mathbb F_{=\beta_{ji}}(\widetilde{\mathcal M})}\big)_{\delta'} = \D^{\bar\Lambda}\big((\mathcal O_{\widetilde{\mathcal M}})_{\delta}\big) +   \bigoplus_{\delta'_1+\delta'_2= \delta'}\big(\mathcal O_{\mathbb F_{=\beta_{ji}}(\widetilde{\mathcal M})}\big)_{\delta'_1} \cdot \big(\mathcal O_{\mathbb F_{=\beta_{ji}}(\widetilde{\mathcal M})}\big)_{\delta'_2},
$$
where $\delta'_i\ne 0$. Note that $\widetilde{\psi}_{\delta'}$ is already defined on the second summand. Further we put 
\begin{equation}\label{eq def of psi}
\begin{split}
&\widetilde{\psi}_{\delta'}|_{\D^{\bar\Lambda}\big((\mathcal O_{\widetilde{\mathcal M}})_{\delta}\big)} : \D^{\bar\Lambda}\big((\mathcal O_{\widetilde{\mathcal M}})_{\delta}\big) \to \mathcal S_{\delta'}, \\
&\widetilde{\psi}_{\delta'}|_{\D^{\bar\Lambda}\big((\mathcal O_{\widetilde{\mathcal M}})_{\delta}\big)}= (\D^{\bar\Lambda})^{-1}.
\end{split}
\end{equation}
Let us show that $\widetilde{\psi}_{\delta'}$ is well-defined. Assume that 
$$
f\in \D^{\bar\Lambda}\big((\mathcal O_{\widetilde{\mathcal M}})_{\delta}\big) \cap   \bigoplus_{\delta'_1+\delta'_2= \delta'}\big(\mathcal O_{\mathbb F_{=\beta_{ji}}(\widetilde{\mathcal M})}\big)_{\delta'_1} \cdot \big(\mathcal O_{\mathbb F_{=\beta_{ji}}(\widetilde{\mathcal M})}\big)_{\delta'_2}.
$$
Then we have
\begin{align*}
\widetilde{\psi}_{\delta'}(f) = (\D^{\bar\Lambda})^{-1} (f) = \iota\circ \varphi_{\mathcal M}\circ \D^{\bar\Lambda}((\D^{\bar\Lambda})^{-1} (f) ) = \iota\circ \varphi_{\mathcal M} (f).
\end{align*}

Now our goal is to define $\widetilde{\psi}_{\theta'}$ for other $\theta'\in  S_p$. Let us take any $\theta'\in S_p$. Then there exists operators $\D_{\gamma_1},\ldots, \D_{\gamma_t}$ and a weight $\delta'$ as above such that $\theta'= \D_{\gamma_1}\circ \cdots \circ \D_{\gamma_t}(\delta')$. In this case we put 
$$
\widetilde{\psi}_{\theta'}:= (\D_{\gamma_1}\circ \cdots \circ \D_{\gamma_t}) \circ \widetilde{\psi}_{\delta'}\circ (\D_{\gamma_1}\circ \cdots \circ \D_{\gamma_t})^{-1}.
$$  
Note that the composition $\D_{\gamma_1}\circ \cdots \circ \D_{\gamma_t}$ is unique up to sign. Hence $\widetilde{\psi}_{\theta'}$ is well-defined. It can be easily shown that $(\widetilde{\psi}_{\theta'})$ satisfies the conditions of Proposition \ref{prop isomorphism of graded mnf}. Therefore the following morphism of graded manifolds $\widetilde{\psi}=(\widetilde{\psi}_{\theta'})$, where $\theta'\in \Delta'_{\mathcal M}\cup S_p$, is defined.

 \medskip
 
\noindent{\bf Step 3, isomorphism commutes with the operators.} It remains to show that 
\begin{equation}\label{eq iso psi commutes with operators}
 \D_{\beta_{st}} \circ \widetilde{\psi}_{\theta'} = \widetilde{\psi}_{\D_{\beta_{st}}(\theta')} \circ \D_{\beta_{st}}
\end{equation}
for $(s,t) \leq (j,i)$. If $\theta'\in \Delta'_{\mathcal M}$, then (\ref{eq iso psi commutes with operators}) holds by induction. Assume that $\theta'\in S_p$ and consider the following cases. 

\medskip

\noindent{\bf 1.} Assume that $\D_{\beta_{st}}(\theta')\in \Delta'$. Then (\ref{eq iso psi commutes with operators}) follows by definition. 

\medskip

\noindent{\bf 2.} Assume that $\D_{\beta_{st}}(\theta')\notin \Delta'$. We also may assume that $\theta'$ depends on $\alpha_t$, since otherwise (\ref{eq iso psi commutes with operators}) holds trivially. In this case $\theta'$ depends on $\beta_{st}$ since otherwise $\D_{\beta_{st}}(\theta')\in \Delta'$. 
 
 {\bf 2.1.} Assume that $\theta'$ satisfy the following property: if $\theta'$ depends on $\beta_{st}$, then $\theta'$ depends also on $\alpha_t$. By (\ref{eq def of psi}), the equality  (\ref{eq iso psi commutes with operators}) holds by the decomposition and induction argument. 
 \smallskip
 
 {\bf 2.2.} Assume that $\theta'= (\D_{\gamma_1}\circ \cdots \circ \D_{\gamma_p})(\theta'_1)$, where $\theta'_1$ is from $(2.1)$. Then we have
 \begin{align*}
 & \D_{\beta_{st}} \circ \widetilde{\psi}_{\theta'} =  \D_{\beta_{st}} \circ (\D_{\gamma_1}\circ \cdots \circ \D_{\gamma_p}) \circ \widetilde{\psi}_{\delta'}\circ (\D_{\gamma_1}\circ \cdots \circ \D_{\gamma_p})^{-1}= \\
  &(\D_{\gamma_1}\circ \cdots \circ \D_{\gamma_p}) \circ \widetilde{\psi}_{\D_{\beta_{st}}(\delta')}\circ (\D_{\gamma_1}\circ \cdots \circ \D_{\gamma_p})^{-1} \circ \D_{\beta_{st}} = \widetilde{\psi}_{\D_{\beta_{st}}(\theta')} \circ \D_{\beta_{st}}. 
 \end{align*}
 The proof is complete.$\Box$

\medskip

\prop \label{prop reconstruction of a graded manifold} {\sl
Let $\mathbf D$ be an $r'$-fold vector bundle of type $\Delta'$ with a family of $(r'-r)$ odd commuting non-degenerate operators $\D_{\beta_{ij}}$ of weights $\beta_{ji}- \alpha_i$, where $i=1,\ldots, r$ and $j=1\ldots, n_i$, satisfying properties $1-6$, Section $7.1$. Then there exists a graded manifold $\mathcal N$ of type $\Delta = \mathsf{G}(\Delta')$ such that $\mathbb F(\mathcal N) \simeq \mathbf D$.

}

\medskip

\noindent{\it Proof.} The proof follows by induction from Propositions \ref{prop raznesenie}, \ref{prop phi icommute with operators}  and \ref{prop prisoedinenie weights}.$\Box$

\medskip

It is remaining to show that  $\mathbb F$ is full and faithful, see Definition \ref{de equivalence of categories}. 

\medskip
\prop\label{prop F is full and faithful} {\sl The functor $\mathbb F$ is full and faithful.}

\medskip
\noindent{\it Proof.} Let us take two objects in the category  $\Delta'$\textsf{VBVect}, i.e. two $r'$-fold vector bundles $\mathbf D_1$ and $\mathbf D_2$ of type $\Delta'$ and a morphism $\Psi:\mathbf D_1\to \mathbf D_2$ that commutes with all vector fields $\D_{\beta_{ji}}$. We have seen in Proposition \ref{prop reconstruction of a graded manifold} that there exist graded manifolds $\mathcal N_i$ of type $\Delta$ such that $\mathbf D_i\simeq \mathbb F(\mathcal N_i)$, where  $i=1,2$. 

Further, let us take two chats $\mathcal U_1$ and $\mathcal U_2$ on  $\mathbf D_1$ and $\mathbf D_2$, respectively, such that we can consider the restriction $\Psi: \mathcal U_1 \to \mathcal U_2$. Denote by $\mathcal V_i$ the corresponding to $\mathcal U_i$  chart on $\mathcal N_i$.  Let us take $\delta$ and $\delta'$ are as in Proposition \ref{prop D_delta is an iso onto kernel}. By Proposition \ref{prop coordinates on D_N}, we can chose coordinates $\zeta^{\delta'}$ and $\xi^{\delta}$ in $\mathcal U_2$ and $\mathcal V_2$, respectively, and there exists unique up to even permutation the operator  $\D^{\bar{\Lambda}}$ such that $\D^{\bar{\Lambda}}(\xi^{\delta}) = \zeta^{\delta'}$. Consider $f=\Psi^*(\zeta^{\delta'})$. By Proposition \ref{prop D_delta is an iso onto kernel} there exists unique function $F\in \mathcal O_{\mathcal V_1}$ such that $\D^{\bar{\Lambda}}(F)=f$. Now we can define the morphism $\Phi|_{\mathcal V_1}: \mathcal V_1\to \mathcal V_2$ by $\Phi(\xi^{\delta}):= F$. (Note that for any $\delta\in \Delta$ there exists $\delta'\in \Delta'$ as in Proposition \ref{prop D_delta is an iso onto kernel}.) Since function $F$ if unique, the morphisms $\Phi|_{\mathcal V_p}$ coincide in all intersections $\mathcal V_s\cap \mathcal V_t$ and defone the global morphism $\Phi$. Clearly, $\mathbb F(\Phi)= \Psi$. The proof is complete.$\Box$ 

\medskip

Now we can formulate our main result.

\bigskip 
\t\label{teor main result} {\sl The categories $\Delta$\textsf{Man} and $\Delta'$\textsf{VBVect} are equivalent. }

\medskip
\noindent{\it Proof.} The proof follows from Propositions \ref{prop reconstruction of a graded manifold} and \ref{prop F is full and faithful}.$\Box$

\medskip

\section{Appendix}

In this section we will prove Propositions \ref{prop D_delta = d circ d circ is injective} and \ref{prop D_delta is an iso onto kernel}.

\medskip
\noindent{\bf Proof of Proposition \ref{prop D_delta = d circ d circ is injective}.} Let us fix an operator $\d_{\gamma_p}$, where $p\in\{1,\ldots, s\}$, from Sequence (\ref{eq D_delta = circ circ}). Clearly, we can rewrite $\D^{\bar\Lambda}$ in the following form
\begin{align*}
\D^{\bar\Lambda}= (\d_{\gamma_1}\,\, \mod\,\, \mathcal J^-) \circ \cdots  \circ (\d_{\gamma_p}\,\, \mod\,\, \mathcal J^-) \circ \cdots \circ (\d_{\gamma_s} \,\, \mod\,\, \mathcal J^-).
\end{align*}
We put 
$$
\D^{\bar\Lambda} = \D_1\circ (\d_{\gamma_p}\,\, \mod\,\, \mathcal J^-) \circ \D_2,
$$
where the notations $\D_1$ and $\D_2$ have obvious meaning. 
Assume that $\D_2(f) \ne 0$ in $\mathcal O_{\widetilde{\mathcal N}}/ \mathcal J^-$, where $f\in (\mathcal O_{\mathcal N})_{\delta}$. Our goal is to show that 
$$
(\d_{\gamma_p}\,\, \mod\,\, \mathcal J^-) (\D_2(f))\ne 0.
$$

Assume that $\gamma_p=\beta_{ji}$. We work locally in a chart $\mathcal U$ on the non-negatively graded manifold $\mathcal M:=(\mathcal N_0, \mathcal O_{\widetilde{\mathcal N}}/ \mathcal J^-)$. Note that we can divide all local homogeneous coordinates in $\mathcal U$  into three groups: coordinates with weight $0$; coordinates with weights of the form $c \alpha_i + \cdots$,  where $c> 0$; and all other coordinates. Hence we can write $\D_2(f)$ in the following form:
$$
\D_2(f) = \sum_{IJ} f_{IJ} \xi^I\eta^J.
$$ 
Here $I$ and $J$ are multi-indexes, $f_{IJ}\in (\mathcal O_{\mathcal N})_0$, $\xi^I$ are monomials in local homogeneous coordinates from the second group, and $\eta^J$ are monomials in local homogeneous coordinates from the third group. Since $\D^{\bar\Lambda}(\delta)$ does not have negative coefficients, we observe that the weight $(\d_{\gamma_p} \circ \D_2) (\delta)$ also does not have negative coefficients. Hence the weight $\D_2 (\delta)$ is a weight of the form $c \alpha_i + \cdots$,  where $c> 0$.  It follows that $\D_2(f)$ depends on coordinates from the second group non-trivially. Let us apply $(\d_{\beta_{ji}}\mod \mathcal J^-)$ to the function $\D_2(f)$. Since $\d_{\beta_{ji}} (f_{IJ})\in \mathcal J^- $ and $\d_{\beta_{ji}} (\eta^J)\in \mathcal J^-,$ 
we get
$$
(\d_{\beta_{ji}}\mod \mathcal J^-) \Big(\sum_{IJ} f_{IJ} \xi^I\eta^J\Big) =  \sum_{IJ} f_{IJ} \d_{\beta_{ji}}( \xi^I)\eta^J \mod \mathcal J^-.
$$
Assume that
\begin{equation}\label{eq f_IJ xi eta =0}
\sum_{IJ} f_{IJ} \d_{\beta_{ji}}( \xi^I)\eta^J \mod \mathcal J^- = 0.
\end{equation} 
Since $\D^{\bar\Lambda}(\delta)$ is multiplicity free and the weight of $\d_{\beta_{ji}}( \xi^I)$ contains the summand $\beta_{ji}$, we see that in the expression for weights of  $\eta^J$ we do not have the summand $\beta_{ji}$. Therefore,  the equation (\ref{eq f_IJ xi eta =0}) is equivalent to the vanishing of all coefficients before $\eta^J$:
$$
\sum_{I} f_{IJ} \d_{\beta_{ji}}( \xi^I)=0\quad \text{for any} \,\,\, J.
$$
Note that here we do not need to assume that the equality holds $\mod \mathcal J^-$.

Further, $f_{IJ}$ is a function of weight $0$. In other words it is a usual smooth (or holomorphic) function that is defined in $\mathcal U_0$.  Let us take a point $x\in \mathcal U_0$ and evaluate the function $f_{IJ}$ at $x$. We get
$$
\sum_{I} f_{IJ}(x) \d_{\beta_{ji}}( \xi^I)= \d_{\beta_{ji}}\big(\sum_{I}   f_{IJ}(x)\xi^I\big) =0\quad \text{for any} \,\,\, x\,\, \text{and}\,\, J.
$$
Since $\d_{\beta_{ji}}$ is the de Rham differential and we can consider $f_{IJ}(x)\xi^I$ as an exterior form of degree $0$ (or just a function) for this operator, we conclude from this equation that $f_{IJ}(x)\xi^I$ is a constant function. Therefore $f_{IJ}\xi^I$ does not depend on $\xi^I$ for any $x$. This contradicts to the fact that the weight of $\D_{2}(f)$ depends on $\alpha_i$ non-trivially and $(\d_{\gamma_p}\,\, \mod\,\, \mathcal J^-) (\D_2(f))\ne 0$ is proven. Since this holds for any $p$, the result follows.$\Box$

\medskip

\medskip
\noindent{\bf Proof of Proposition \ref{prop D_delta is an iso onto kernel}.} Since $\d_{\beta_{ji}}\circ \d_{\beta_{ji}}=0$, see (\ref{eq commuting de Rham operators on N'}), and therefore $\d_{\beta_{ji}}\circ \d_{\beta_{ji}} \,\,\mod \mathcal J^- =0$,  we have
\begin{equation}\label{eq image of Psi_delta}
\D^{\bar\Lambda} \big((\mathcal O_{\mathcal N})_{\delta}\big) \subset  \Big(	(\mathcal O_{\mathbf D_{\mathcal N}})_{\D^{\bar\Lambda}(\delta)} \Big)\bigcap_{k=1}^s \Ker \D_{\gamma_{k}}.
\end{equation}
Our goal is to prove that the inclusion (\ref{eq image of Psi_delta}) is in fact the equality. As in the proof of Proposition \ref{prop D_delta = d circ d circ is injective}, let us write the sequence (\ref{eq D_delta = circ circ}) in the following form:
$$
\D_1\circ (\d_{\gamma_{p}}\,\, \mod\,\, \mathcal J^-) \circ \D_2,
$$ 
where $\d_{\gamma_{p}} = \d_{\beta_{ji}}$,
and let us take any 
$$
f'\in \Big(	(\mathcal O_{\mathbf D_{\mathcal N}})_{\D^{\bar\Lambda}(\delta)} \Big)\bigcap_{k=1}^s \Ker \D_{\gamma_{k}}.
$$
Assume by induction that we found an element 
$$
f\in (\mathcal O_{\widetilde{\mathcal N}}/\mathcal J^-)_{(\d_{\gamma_{p}}\circ \D_2)(\delta)}
$$
such that $\D_1(f)=f'$ and $(\d_{\gamma_{q}} \mod \mathcal J^-)(f) = 0$
for  $q=p,p+1,\ldots ,s$. We have to show that there exists 
$
F\in (\mathcal O_{\widetilde{\mathcal N}}/\mathcal J^-)_{\D_2(\delta)}
$
such that 
$$
(\d_{\beta_{ji}} \mod \mathcal J^-)(F) =f \quad \text{and}\quad (\d_{\gamma_{q}}\mod \mathcal J^-)(F)=0 
$$
for any  $q=p+1,p+2,\cdots, s$.

Consider a chart $\mathcal U$ on $\mathcal M$ as in the proof of Proposition \ref{prop D_delta = d circ d circ is injective} 
with coordinates $(x_i)$, $(\xi_j)$ and $(\eta_t)$ from the groups $1$, $2$ and $3$, respectively. We can write $f$ locally in the following form:
$$
f= \sum_{I,J,u} f_{IJu} \xi^I \d_{\beta_{ji}}(\xi_u) \eta^J.
$$
Here $I$ and $J$ are multi-indexes, $f_{IJu}$ are functions of weight $0$; $\xi^I$ and $\eta^J$ are monomials in local homogeneous coordinates from the second and third  groups, respectively;  and $\d_{\beta_{ji}}(\xi_u)$ are local coordinates from the third group which weights contain $\beta_{ji}.$ 
By our assumption, we have 
\begin{equation}\label{eq d_beta_ji (f=0)}
(\d_{\beta_{ji}} \mod \mathcal J^-)(f) = \sum_{I,J,u} f_{IJu} \d_{\beta_{ji}}(\xi^I) \d_{\beta_{ji}}(\xi_u) \eta^J =0 \,\, \mod\,\, \mathcal J^-.
\end{equation}
Note that  to obtain the first equality in (\ref{eq d_beta_ji (f=0)}), we use the following facts 
$$
\d_{\beta_{ji}}(f_{IJu})\in \mathcal J^-,\quad \d_{\beta_{ji}}(\eta^J )\in \mathcal J^-\quad \text{and} \quad \d_{\beta_{ji}}\circ \d_{\beta_{ji}}=0.
$$
Since weights of monomials $\eta^J$ do not contain $\beta_{ji}$, the equality (\ref{eq d_beta_ji (f=0)}) is equivalent to 
$$
\d_{\beta_{ji}} \big(\sum_{I,u} f_{IJu}(x_0) \xi^I \d_{\beta_{ji}}(\xi_u)\big)  =0 \quad \text{for any} \,\, J\,\,\text{and any}\,\,x_0\in \mathcal U_0.
$$
We see that $ \sum\limits_{I,u} f_{IJu}(x_0) \xi^I \d_{\beta_{ji}}(\xi_u)$ is a closed $1$-form in the superdomain with coordinates $(\xi_j)$, with respect to the de Rham differential $\d_{\beta_{ji}}$.  By the Poincar\'{e} Lemma for graded manifolds, for any $x_0$ and $J$ there exists 
$$
F_{J}(x_0) = \sum_{K} F_{KJ}(x_0) \xi^K \quad \text{such that}\quad  \d_{\beta_{ji}}(F_{J}(x_0)) = \sum_{Iu}  f_{IJu}(x_0) \xi^I \d_{\beta_{ji}}(\xi_u).
$$
Here $K$ is a multi-index.
In particular we have  
\begin{equation}\label{eq F_J}
\frac{\partial}{\partial \xi_u} (F_{J}(x_0)) = \sum_{I} f_{IJu}(x_0) \xi^I. 
\end{equation}
Note that $F_{J}(x_0)$ is defined up to a constant. However, if we assume that 
$$
\text{weight}(F_{J}(x_0)) = \text{weight}\big(\sum_{I} f_{IJu}(x_0) \xi^I\big) - (\beta_{ji}-\alpha_i),
$$
then $F_{J}(x_0)$ is unique. 

Now we need  to show that $F_{J}(x_0)$ is a smooth (or holomorphic) function in $x_0$. Consider the following equality 
$$
F_{KJ} =  \frac{\partial}{\partial \xi^K} (F_{J}),
$$
where $\frac{\partial}{\partial \xi^K}$ is the corresponding to $\xi^K$ differential operator. Using (\ref{eq F_J}), we see that the function $F_{KJ}$ is a certain iterated derivative of $\sum\limits_{I} f_{IJu} \xi^I$, and hence it is smooth (or holomorphic).

Summing up, we constructed the functions $F_J= \sum\limits_{K} F_{KJ} \xi^K$ such that 
$$
\d_{\beta_{ji}}(F_{J}) = \sum_{Iu}  f_{IJu} \xi^I \d_{\beta_{ji}}(\xi_u) \,\, \mod \mathcal J^-.
$$
We put 
$$
F:= \sum_{KJ} F_{KJ}\xi^K \eta^J.
$$
Clearly the weight of $F$ is $\D_2(\delta)$ and we have $\d_{\beta_{ji}}(F)\,\, \mod \mathcal J^- = f$. 

It is remaining to show that $\d_{\gamma_{q}}(F)=0\,\, \mod \,\, \mathcal J^-$ for $q=p+1,\ldots, s$. 
Consider the function $H:= (\d_{\gamma_{q}} \, \mod \mathcal J^-)(F)$, where $q=p+1,\ldots, s$. The weight of $H$ does not contain $\beta_{ji}$ since $\D^{\bar\Lambda}(\delta)$ is multiplicity free, and it has the form $c\alpha_i+\ldots$, where $c>0$, since our assumption that the weight $\D^{\bar\Lambda}(\delta)$ depends on $\alpha_i$ non-trivially for any $i=1,\ldots, r$. Therefore, $H$ is a non-constant function for the de Rham differential $\d_{\beta_{ji}}$. Further,
\begin{align*}
&(\d_{\beta_{ji}} \, \mod \mathcal J^-) (H)= (\d_{\beta_{ji}} \, \mod \mathcal J^-)\circ (\d_{\gamma_{q}} \, \mod \mathcal J^-)(F) = \\
-&(\d_{\gamma_{q}} \, \mod \mathcal J^-)\circ (\d_{\beta_{ji}} \, \mod \mathcal J^-) (F) = -(\d_{\gamma_{q}} \, \mod \mathcal J^-)(f)=0.
\end{align*} 
Hence $H=0$.

If we iterate our construction, we get a function 
$$
F'\in (\mathcal O_{\widetilde{\mathcal N}}/\mathcal J^-)_{\delta} = (\mathcal O_{\mathcal N})_{\delta}
$$
such that $\D^{\bar\Lambda}(F') = f'$. The proof is complete.$\Box$

\medskip

\bigskip

\bigskip

\noindent{\it Elizaveta Vishnyakova}

\noindent {Universidade Federal de Minas Gerais, Brazil}

\noindent{\emph{E-mail address:}
	\verb"VishnyakovaE@googlemail.com"}


\begin{thebibliography}{99}

\bibitem[BGG]{BruceGrGr}  {\it Bruce A., Grabowska K., Grabowski J.} Linear duals of graded bundles and higher analogues of (Lie) algebroids, J. Geom. Phys. 101 (2016), 71-99. 

\bibitem[BGR]{BruceGrab SIGMA}  {\it Bruce A., Grabowski J. and Rotkiewicz M.} Polarisation of graded bundles. SIGMA 12 (2016), 106, 30 pages

\bibitem[BCMZ]{Bursztyn}  {\it Bursztyn, H., Cattaneo, A., Mehta R. and Zambon M.} Reduction of Courant algebroids via supergeometry. Work in progress.

\bibitem[BCMZ]{BursztynAdv}  {\it Bursztyn, H., Cavalcanti, G. R. and Gualtieri, M.} Reduction of Courant algebroids and generalized complex structures, Advances in Mathematics 211(2) (2007), 726-765.

\bibitem[CM]{Fernando}  {\it del Carpio-Marek. F.} Geometric structure on degree $2$ manifolds. PhD-thesis, IMPA, Rio de Janeiro, 2015.

\bibitem[Cou]{Courant}  {\it Courant, T.} Dirac manifolds, Trans. Amer. Math. Soc. 319 (1990), 631-661.


\bibitem[GR]{GR}  {\it Grabowski, J.,
	Rotkiewicz, M.} Higher vector bundles and multi-graded symplectic manifolds, J. Geom. Phys., Volume 59, Issue 9, 2009, Pages 1285-1305.

\bibitem[GrMa]{GraciaMackenzie}  {\it Gracia-Saz A., Mackenzie K.C.H.} Duality functors for triple vector bundles, Lett. Math. Phys. 90 (2009), 175-200.

\bibitem[GrMe]{GraciaMetha}{\it Gracia-Saz, A. and Mehta, R.A.} Lie algebroid structures on double vector bundles and representation theory of Lie algebroids, Advances in Mathematics 223 (4), 2010, Pages 1236-1275.



\bibitem[Gu]{Gualteri}  {\it Gualtieri, M.} Generalized complex geometry, Ann. of Math. Volume 174 (2011), Issue 1 (2011), Pages 75-123.


\bibitem[KS]{Kosmann}  {\it Kosmann-Schwarzbach, Y.} Courant Algebroids. A Short History. SIGMA 9 (2013), 014, 8 pages.

\bibitem[JL]{JL}  {\it Jotz Lean, M.} N-manifolds of degree 2 and metric double vector bundles,  arXiv:1504.00880.

\bibitem[L]{Leites} {\it Leites D.A.} Introduction to the theory of supermanifolds. Uspekhi Mat. Nauk, 1980, Volume 35,	Issue 1(211), 3-57
	
\bibitem[LB]{Li-Bland}  {\it Li-Bland, D.S.} LA-Courant algebroids and their applications. Ph.D. Thesis, University of Toronto, 2012.

\bibitem[LS]{LS}{\it Lyakhovich S.L., Sharapov A.A.}
Characteristic classes of gauge systems
Nucl. Phys. B703 (2004) 419-453, e-Print: hep-th/0407113.

	
\bibitem[Man]{Man} {\it Manin Yu.I.} Gauge field theory and complex geometry, Grundlehren Math.
Wiss., vol. 289, Springer-Verlag, Berlin 1988, 1997.

	
\bibitem[M1]{Mackenzie}  {\it Mackenzie, K.C.H.} General theory of Lie groupoids and Lie algebroids. Volume 213 of London Mathematical Society Lecture Note Series. Cambridge: Cambridge University Press, 2005.


\bibitem[M2]{Mackenzie Crelles}  {\it Mackenzie, K.C.H.} Ehresmann doubles and Drinfel'd doubles for Lie algebroids and Lie bialgebroids. Journal f\"{u}r die reine und angewandte Mathematik, Volume 2011, Issue 658 (2011).

\bibitem[R]{Roytenberg}  {\it Roytenberg, D.} On the structure of graded symplectic supermanifolds and Courant algebroids. Contemp. Math., Vol. 315, Amer. Math. Soc., Providence, RI, (2002)
	
\bibitem[S]{Severa}  {\it \v{S}evera, P.} Some title containing the words ``homotopy'' and ``symplectic'', e.g. this one, Travaux math\'{e}matiques. Fasc. XVI (2005), Pages 121-137
	
\bibitem[Va]{Va}  {\it Vaintrob A.} Lie algebroids and homological vector fields, Rossi\u{i}skaya Akademiya
Nauk. Moskovskoe Matematicheskoe Obshchestvo. Uspekhi Matematicheskikh Nauk, 52, 2(314), 1997, 161-162.
	
\bibitem[Vit]{Luca}  {\it Vitagliano L.} Vector bundle valued differential forms on ℕQ-manifolds. Pacific Journal of Mathematics 283, 2 (2016) 449-482.


\bibitem[Vo1]{Vor}  {\it Voronov, Th.Th.} $Q$-Manifolds and Mackenzie Theory. Communications in Mathematical Physics, 2012, Volume 315, Issue 2, pp 279-310.
	
\bibitem[Vo2]{Voronov graded}  {\it Voronov, Th.Th.} Graded manifolds and Drinfeld doubles for Lie bialgebroids.  Quantization, Poisson Brackets and Beyond. Volume 315 of Contemp. Math., Providence, RI: Amer. Math. Soc., 2002, pp. 131-168


\end{thebibliography}
\end{document}